\documentclass[11pt]{amsart}
\usepackage[letterpaper,margin=1in]{geometry}                
\usepackage{amssymb,amsmath,amsfonts,stmaryrd}
\usepackage{xcolor}
\usepackage[colorlinks,
             linkcolor=blue,
             citecolor=green!70!black,
             pdfproducer={pdfLaTeX},
             pdfpagemode=None,
             bookmarksopen=true
             bookmarksnumbered=true]{hyperref}
\usepackage{graphicx}
\usepackage{multicol}

\usepackage{tikz}
\usetikzlibrary{decorations.markings,arrows,decorations.pathreplacing}
\usepackage{arydshln}

\newcommand\arXiv[1]{\href{http://arxiv.org/abs/#1}{\nolinkurl{arXiv:#1}}}
\newcommand\MRnumber[1]{\href{http://www.ams.org/mathscinet-getitem?mr=#1}{\nolinkurl{MR#1}}}
\newcommand\DOI[1]{\href{http://dx.doi.org/#1}{\nolinkurl{DOI:#1}}}
\newcommand\MAILTO[1]{\href{mailto:#1}{\nolinkurl{#1}}}

\newtheorem*{theorem}{Theorem}

\newtheorem*{proposition}{Proposition}
\newtheorem*{corollary}{Corollary}
\newtheorem*{lemma}{Lemma}
\newtheorem*{conjecture}{Conjecture}

\usepackage{dsfont}
\renewcommand\mathbb\mathds

\newcommand\bC{\mathbb C}

\newcommand\bF{\mathbb F}

\newcommand\bR{\mathbb R}

\newcommand\bZ{\mathbb Z}

\newcommand\rA{\mathrm A}
\newcommand\rB{\mathrm B}
\newcommand\rC{\mathrm C}
\newcommand\rD{\mathrm D}
\newcommand\rE{\mathrm E}
\newcommand\rF{\mathrm F}
\newcommand\rG{\mathrm G}
\newcommand\rH{\mathrm H}

\newcommand\rJ{\mathrm J}

\newcommand\rL{\mathrm L}
\newcommand\rM{\mathrm M}

\newcommand\rO{\mathrm O}
\newcommand\rP{\mathrm P}

\newcommand\rT{\mathrm T}
\newcommand\rU{\mathrm U}

\usepackage[mathscr]{euscript}

\DeclareMathOperator\homology{H}
\renewcommand\H{\homology}

\renewcommand\d{\mathrm d}

\newcommand\longto\longrightarrow
\newcommand\mono\hookrightarrow
\newcommand\epi\twoheadrightarrow
\newcommand\isom{\overset\sim\to}
\newcommand\<\langle
\renewcommand\>\rangle
\newcommand\sminus\smallsetminus

\newcommand\st{\text{ s.t.\ }}

\newcommand\id{\mathrm{id}}

\newcommand\Co{\mathrm{Co}}

\newcommand\Suz{\mathrm{Suz}}

\newcommand\Spin{\mathrm{Spin}}

\newcommand\SU{\mathrm{SU}}
\newcommand\Sp{\mathrm{Sp}}
\newcommand\SO{\mathrm{SO}}
\newcommand\PSO{\mathrm{PSO}}
\newcommand\PSU{\mathrm{PSU}}
\newcommand\PSp{\mathrm{PSp}}

\DeclareMathOperator\Aut{Aut}
\DeclareMathOperator\Sym{Sym}
\DeclareMathOperator\Alt{Alt}

\DeclareMathOperator\Fer{Fer}
\DeclareMathOperator\Cliff{Cliff}

\DeclareMathOperator\Sq{Sq}
\DeclareMathOperator\SH{SH}
\DeclareMathOperator\rSH{rSH}

\newcommand\tr{\mathrm{tr}}

\newcommand\ev{\mathrm{ev}}
\newcommand\odd{\mathrm{odd}}

\newcommand\adj{\mathrm{adj}}

\newcommand\vac{(\mathrm{vac})}

\DeclareMathOperator\diag{diag}

\newcommand\define[1]{\emph{#1}}
\newcommand\cat[1]{\textsc{#1}}

\title{Supersymmetry and the Suzuki chain}
\author{Theo Johnson-Freyd}
\email{\MAILTO{theojf@pitp.ca}}
\address{Perimeter Institute for Theoretical Physics, Waterloo, Ontario}

\begin{document}
\begin{abstract}
  We classify $N{=}1$ SVOAs with no free fermions and with bosonic subalgebra a simply connected WZW algebra which is not of type $\rE$. 
  The latter restriction makes the classification tractable; the former restriction implies that the $N{=}1$ automorphism groups of the resulting SVOAs are finite.
  We discover two infinite families and nine exceptional examples. The exceptions are all related to the Leech lattice: their automorphism groups are the larger groups in the Suzuki chain ($\Co_1$, $\Suz{:}2$, $\rG_2(4){:}2$, $\rJ_2{:}2$, $\rU_3(3){:}2$) and certain large centralizers therein ($2^{10}{:}\rM_{12}{:}2$, $\rM_{12}{:}2$, $\rU_4(3){:}D_8$, $\rM_{21}{:}2^2$). Along the way, we elucidate fermionic versions of a number of VOA operations, including simple current extensions, orbifolds, and 't Hooft anomalies.
\end{abstract}
\maketitle

\section{Introduction}

The goal of this article is to prove the following:
\begin{theorem}
  Suppose $V$ is an $N{=}1$ SVOA with no free fermions, and that $V_{\ev} = G_k$ is a  simply connected WZW algebra which is not one of $\rE_{7,2}$, $\rE_{7,1}^2$, or $\rE_{8,2}$. Then
   $V$ is on the following list:
$$\begin{array}{l|l|l|l}
 V_{\ev} & \dim V_{3/2} & c & \Aut_{N{=}1}(V) \\ \hline &&& \\[-9pt]
 \Spin(m)_3 & \frac{m(m-1)(m+4)}6 & \frac{3m(m-1)}{2(m+1)} & S_{m+1} \\[3pt]
 \Spin(m)_1^3 & m^3 & \frac{3m}2 & \begin{cases} 2^{2(m-1)}{:}(S_3 \times S_m), & m \neq 4 \\
 2^6 {:} 3S_6, & m = 4 \end{cases}
 \\[3pt]
 \Sp(2{\times}3)_2 & 84 & 7 & \rU_3(3){:}2 \\[3pt]
 \Sp(2{\times}3)_1^2 & 196 & 8 \frac25 & \rJ_2{:}2 \\[3pt]
 \SU(6)_2 & 175 & 8 \frac34 & \rM_{21}{:}2^2 \\[3pt]
 \Sp(2{\times}6)_1 & 429 & 9 \frac34 & \rG_2(4){:}2\\[3pt]
 \SU(6)_1^2 & 400 & 10 & \rU_4(3){:}D_8 \\[3pt]
 \Spin(12)_2 & 462 & 11 & \rM_{12}{:}2 \\[3pt]
 \SU(12)_1 & 924 & 11 & \Suz{:}2 \\[3pt]
 \Spin(12)_1^2 & 1024 & 12 & 2^{10}{:}\rM_{12}{:}2 \\[3pt]
 \Spin(16)_1 \times \Spin(8)_1 & 1024 & 12 & 2^8 \cdot \rO_8^+(2).2 \\[3pt]
 \Spin(24)_1 & 2048 & 12 & \Co_1 \\[3pt]
\end{array}$$
 In each case on the list, the $N{=}1$ SVOA $V$ is uniquely determined up to isomorphism by $G_k$, and the supersymmetry-preserving automorphism group of $V$ is the listed finite group $\Aut_{N{=}1}(V)$.
 \end{theorem}
 The type $\rE$ case is discussed in Section~\ref{sec.typeE}, 
 where we construct an $N{=}1$ SVOA with even subalgebra $\rE_{7,1}^2$, but do not prove its uniqueness, and explain why the other two cases seem unlikely.
\begin{conjecture}
 There is a unique-up-to-isomorphism $N{=}1$ SVOA with even subalgebra $\rE_{7,1}^2$. There does not exist an $N{=}1$ SVOA with even subalgebra $\rE_{7,2}$ or $\rE_{8,2}$.
\end{conjecture}

We now elaborate on the statement of the Theorem. Our notation for groups is the following: as in the ATLAS~\cite{ATLAS}, a colon denotes a semidirect product,
$2 = \bZ_2$ denotes a cyclic group of order $2$, 
$2^{2m}$ is the elementary abelian group of that order,
and $S_m$, $\rJ_2$, etc.\ denote specific finite groups; 
we write $\Sp(2{\times}m)$ for the group of type $\rC_m$ (which is  variously called $\Sp(m)$ and $\Sp(2m)$); 
and the small-$m$ members of the ``$\Spin(m)$'' families must be interpreted appropriately~(see \S\ref{subsec:small-m}).
We will later use the notation $nG$ for a perfect central extension of $G$ by $\bZ_n$ and $\rP G$ for the adjoint form of the simply connected group $G$; when $m$ is divisible by $4$, we will write $\SO^+(m)$ for the image of $\Spin(m)$ in the positive half-spin representation; and we will denote the central product of $G$ and $G'$, where their centres have been identified $Z = Z(G) = Z(G')$, by $G \circ G' = (G \times G')/Z$.

The acronym ``SVOA'' stands for super vertex operator algebra. Elements of an SVOA are the ``fields'' (also typically called ``vertex operators''). Axioms can be found in the standard textbooks \cite{MR1651389,MR2082709,MR2023933}. Our SVOAs $V$ enjoy many niceness axioms: $C_2$-cofiniteness; the $L_0$-action on $V$ is diagonalizable and determines a $\frac12\bZ$-grading $V = \bigoplus_{h \in \frac12\bZ} V_h$ that we call ``spin''; the vacuum vector spans the space $V_0$ of spin-$0$ fields. (The first two of these axioms define ``niceness'' in the sense of \cite{MR1614941,MR2352133}, and the third says that the SVOA is of ``CFT type'' \cite{MR2648364}.) 
All the SVOAs in this paper  are moreover rational (they have finitely many simple modules).
The space $V_1$ of spin-$1$ fields is naturally a Lie algebra,
and a version of Noether's theorem identifies it with the Lie algebra of the group $\Aut(V)$ of SVOA automorphisms of $V$, which is always a reductive algebraic group~\cite{MR2097833}.

We will furthermore require that our SVOA be \define{unitary}, i.e.\ equipped with an invariant positive definite Hermitian form. This requirement does two things. First, it provides inside each graded component $V_h$ a real subspace of self-adjoint, aka real, fields. 
Second, it makes each $V_h$ into a Hilbert space, and picks out a compact form $\Aut(V)_{\mathrm{cpt}}$ of $\Aut(V)_\bC$ acting on each $V_h$ unitarily. Indeed, $V_h$ is the complexification of its self-adjoint subspace, which is a real $\Aut(V)_{\mathrm{cpt}}$-module. Using this, we may move freely between compact and complex forms of groups, and between their unitary (perhaps real) representations and their algebraic (perhaps self-dual) representations.

  A \define{simply connected WZW algebra} is a VOA generated by its spin-$1$ fields (the ``currents'')~\cite{HenriquesWZW}. Simply connected WZW algebras are labeled by a  simply connected reductive Lie group $G$ together with a positive class $k$, called the ``level,'' in the integral cohomology group $\H^4(BG;\bZ)$. In the quasisimple case, $\H^4(BG;\bZ) \cong \bZ$, and the level is a positive integer. We will generally write the WZW algebra corresponding to $(G,k)$ as $G_k$, simplifying for example $(\rE_7)_2$ to $\rE_{7,2}$.  We will call a WZW algebra \define{simple} if the corresponding group is quasisimple.
  (All WZW algebras are ``simple'' in the VOA sense.)
   Automorphism groups of simply connected WZW algebras are always $\Aut(G_k) = \Aut(\mathfrak{g}) = \rP G{:}\operatorname{Out}(\mathfrak{g})$, where $\mathfrak{g}$ is the Lie algebra of $G$, $\rP G$ is its adjoint form, and $\operatorname{Out}(\mathfrak{g})$ is its group of Dynkin diagram automorphisms. Note in particular that although $G$ acts naturally on each $G_k$, the action has kernel the centre of $G$.
  The  SVOAs  in the Theorem are all ``non-simply connected'' WZW algebras. A non-simply connected WZW algebra is determined by a connected but not simply connected group $G$ together with a level which, in the quasisimple case, is a positive integer satisfying a divisibility criterion that depends on ($G$ and) $\pi_1(G)$; details are reviewed in~\S\ref{subsec.extension}. $\rE_{8,2}$ has some peculiar behaviour, and is forbidden by fiat in \cite{HenriquesWZW}, where VOAs containing an~$\rE_{8,2}$ factor are termed ``$\rE_{8,2}$-contaminated.''

An SVOA is called ``$N{=}1$'' if it is equipped with a supersymmetry: a real \define{superconformal vector}, which is a spin-$\frac32$ field $\tau$ (also typically called ``$G$'') satisfying the OPE
$$ \tau(z) \tau(0) \sim \frac{\frac23c}{z^3} + \frac{\nu(0)}{z}, $$
where the spin-$2$ field $\nu$ is the conformal vector (also typically called ``$L$'' or ``$T$'') and $c$ is the bosonic central charge. 
For an $N{=}1$ SVOA $V$, we will adopt the names $\Aut_{N{=}0}(V)$ and $\Aut_{N{=}1}(V)$ for the groups of automorphisms of $V$ as an SVOA and as an $N{=}1$ SVOA, respectively.

A \define{free fermion} in an SVOA is a field of spin $\frac12$. 
We will write $\Fer(m)$ for the SVOA generated by an $m$-dimensional vector space of free fermions, and call it the ``(purely) free fermion'' algebra. It is unique up to isomorphism, and $\Aut(\Fer(m)) = \rO(m)$.
The free fermions in any SVOA split off as a tensor factor~\cite{MR968813},
 although perhaps not compatibly with a chosen supersymmetry.
There are many $N{=}1$ SVOAs with free fermions: highlights include the supersymmetric lattice SVOAs~\cite{HeluaniKac2007} and the beautiful classification identifying $N{=}1$ structures on purely free-fermion algebras with semisimple Lie algebras~\cite{MR791865}.
 Noether's theorem for $N{=}1$ SVOAs identifies $\operatorname{Lie}(\Aut_{N{=}1}(V))$ with the space of free fermions in $V$, equipped with a Lie bracket derived, in the sense of \cite{MR2104437}, from the ordinary Lie bracket:
$$ [x,y]_\tau = [x,[\tau,y]].$$
In particular, if $V$ has no free fermions, then $\Aut_{N{=}1}(V)$ is a finite group.

One reason to be interested in $N{=}1$ SVOAs $V$ without free fermions, and with $V_\ev$ easy to understand (for instance, a  simply connected WZW algebra), comes from Duncan's beautiful work~\cite{MR2352133}.
In that work, Duncan discovered a specific $N{=}1$ SVOA $V^{f\natural}$ with no free fermions; it enjoys $V^{f\natural}_\ev = \Spin(24)_1$ and $\Aut_{N{=}1}(V^{f\natural}) = \Co_1$, Conway's largest sporadic simple group. 
The original motivation for this paper was to produce similar $N{=}1$ SVOAs for  the \define{Suzuki chain} groups
$$ \Suz{:}2 \supset \rG_2(4){:}2 \supset \rJ_2{:}2 \supset \rU_3(3){:}2 \supset \rL_2(7){:}2 \supset A_4{:}2 \supset A_3{:}2. $$
This is a chain of subgroups of $\Co_1$ (except that $\Suz{:}2$ lives therein only through its abelian extension $3\Suz{:}2$).
In fact, the Theorem provides a systematic construction of such SVOAs only for $\Suz{:}2,\dots,\rU_3(3){:}2$. We did not find a satisfactory SVOA representation of $\rL_2(7){:}2$, and our SVOA representations of $A_4{:}2$ and $A_3{:}2$ place them  
naturally as entries in the infinite chain $S_{m+1} = A_{m+1}{:}2$ and not in the Suzuki chain. That said, Wilson organizes $\rU_4(3){:}D_8$, which does appear in the Theorem, with the Suzuki chain \cite{MR723071}, and the $m=8$ entry in our $\Spin(m)_1^3$ family is closely related to the ``tricode group'' $2^{2+12}{:}(S_3 \times A_8)$, which Wilson also discusses in the context of the Suzuki chain.
(The passage $2^{14} \leadsto 2^{2+12}$ occurs because some symmetries are broken, and others are centrally extended, when passing from $\Spin(8)_1^3$ to $\Spin(24)_1$.)

The paper is structured as follows. 
Section~\ref{sec.WZW} classifies the SVOAs $V$ with even subalgebra a simple simply connected WZW algebra which have no free fermions but which do have a nonzero field of spin $\frac32$; we find a short list, consisting of those algebras listed in the Theorem together with  $\rE_{7,2}$, $\rE_{7,1}^2$, and $\rE_{8,2}$. 
The $N{=}1$ structures claimed in the Theorem are constructed in Section~\ref{sec.existence}. Uniqueness (up to conjugation by SVOA automorphisms, of course, with the asserted stabilizer) is shown for the $\Spin(m)_3$ and $\Spin(m)_1^3$ families in Section~\ref{sec.spinm3}, and for the exceptional cases in Section~\ref{sec.uniqueness}. For the exceptional cases we take advantage of Duncan's work on the Conway group $\Co_1$ and its module $V^{f\natural}$, and through this we connect directly to the Suzuki chain. Finally, Section~\ref{sec.typeE} addresses the type $\rE$ case.

\section{WZW algebras with a spin-$\frac32$ abelian anyon}\label{sec.WZW}

This section restricts the possible WZW algebras $G_k$ that can appear as the even parts of the $N{=}1$ SVOAs $V$ considered in the Theorem. First, in \S\ref{subsec.extension}, we ask what WZW algebras can appear as the even parts of (not necessarily supersymmetric) SVOAs with a field of spin $\frac32$ and with no free fermions; in particular, this already cuts the possibilities down to just those simple groups appearing in the Theorem. We then make a comment in \S\ref{subsec:small-m} about some unstable behaviour enjoyed by the Spin groups. In \S\ref{subsec.genus} we use ``elliptic genus'' considerations to rule out most non-simple possibilities. Finally, in \S\ref{subsec.inclusions}, we chart the WZW algebras that require further study.

\subsection{$\bZ_2$ simple current extensions} \label{subsec.extension}

Suppose $V$ is an SVOA, not necessarily supersymmetric, with a nontrivial fermionic part $V_{\odd}$. Then the bosonic subalgebra $V_{\ev}$ is the fixed points of a nontrivial $\bZ_2$-action on $V$ (namely, fermion parity). It follows that $V$ is a ``$\bZ_2$ simple current extension'' of $V_{\ev}$. If $V$ is to be an $N{=}1$ SVOA without free fermions, then $V_{\odd}$ will have conformal dimension exactly $\frac32$.

This means the following. The representations of a VOA are called \define{anyons}. An anyon is \define{abelian} if it is invertible for the braided monoidal structure on the category of anyons; abelian anyons are also called ``simple currents.'' (The braided monoidal category of anyons is a modular tensor category when the VOA is rational \cite{MR2468370}. All VOAs appearing in this paper are rational.) From the unitary representation theory of Virasoro algebras \cite{MR740343}, one sees that each irreducible anyon $M$ has a \define{conformal dimension} $h_M \geq 0$, with equality only for the trivial anyon; conformal dimension is also called ``minimal energy,'' and is by definition the smallest eigenvalue of the action of $L_0 \in V$ on the anyon~$M$.
The abelian anyons form an abelian group $A$, and $$q = (M \mapsto \exp(2\pi i h_M)) : A \to \rU(1)$$ is a quadratic function. Given a subgroup $A' \subset A$, the direct sum $\bigoplus_{M\in A'} M$, admits an SVOA structure if and only if the restriction of $q$ to $A'$ is a group homomorphism to $\{\pm 1\}$. In this case that SVOA structure is unique up to isomorphism and is called the ``$A'$ simple current extension'' of $V$~\cite{MR4050091}.

Simple current extensions of WZW algebras are studied in detail in \cite{MR1822111,HenriquesWZW}. The representation theory of the WZW algebra $G_k$ is well understood, and is easily accessed in Schellekens' computer algebra program ``Kac'' \cite{Kac-computer} and in the detailed tables compiled in \cite{MR1117679}.
(We focus in this section on the case when $\mathfrak{g} = \operatorname{Lie}(G)$ is simple. The general semisimple case is considered in \S\ref{subsec.genus}.)
The irreducible anyons are indexed by dominant integral weights $\lambda$ of $G$ such that $\langle \lambda, \alpha_{\mathrm{max}} \rangle \leq k$, where $\alpha_{\mathrm{max}}$ denotes the highest root; the fields of minimal conformal dimension in the $\lambda$th anyon  form the simple $G$-representation of highest weight $\lambda$.
The conformal dimensions of all anyons are given by an easy formula in terms of the weight lattice.
 With one notable exception, the group $A$ of abelian anyons is naturally isomorphic to the centre~$Z(G)$ of~$G$. The conformal dimension of the $G_k$-anyon corresponding to a fixed element in $Z(G)$ depends linearly on $k$ and is listed in \cite{MR1822111}:
 \begin{lemma}[\cite{MR1822111}] The conformal dimension $h_a$ of the $G_k$-anyon corresponding to $a \in Z(G)$ is:\\[-3pt]

\hfill
$ \begin{array}[b]{l|l|l}
 G & Z(G) & h_a \\ \hline && \\[-9pt]
 \SU(m) & \bZ_m & i \mapsto \frac {ki(m-i)}{2m} \\[3pt]
 \Sp(2{\times}m) & \bZ_2 & \frac {mk} 4 \\[3pt]
 \Spin(m), \, m \,\odd & \bZ_2 & \frac k2 \\[3pt]
 \Spin(m), \, m \,\ev & \bZ_4 \text{ or } (\bZ_2)^2 &  v \mapsto \frac k 2, \; s^\pm \mapsto \frac{km}{16} \\[3pt]
 \rE_6 & \bZ_3 & \frac{2k}3 \\[3pt]
 \rE_7 & \bZ_2 & \frac{3k}4
\end{array}
$ \hfill \qed
\end{lemma}
We list only the nonzero values of $h_a$. In type $\rA$, $i$ ranges over $\{1,\dots,m-1\}$, and in type $\rD$ we write $v \in Z(\Spin(m))$  (``vector'') for the nontrivial element in $\ker(\Spin(m) \to \SO(m))$ and $s^\pm$ (``spinor'') for the other two elements. The exceptional groups $\rG_2$, $\rF_4$, and $\rE_8$ have trivial centre and so are not listed. Levels for the spin groups $\Spin(m)$ with $m\leq 4$ have some unstable behaviour, described in \S\ref{subsec:small-m}, and so the formulas above in those cases must be appropriately interpreted.

 The simple current extension of $G_k$ by $A' \subset Z(G)$, if it exists, is a ``non-simply connected WZW algebra'' corresponding to the quotient group~$G/A'$~\cite{HenriquesWZW}. 
 We will follow the reasonably standard convention that $(G/A')_k = G_k / A'$ denotes the $A'$ simple current extension of $G_k$, so that a level for a non-simply connected group is determined by its pullbacks to the simply connected cover --- note that  this means deciding that for a non-simply connected group, the minimal level may not be ``level $1$'' (cf.\ \S\ref{subsec:small-m}). As will be discussed in \S\ref{subsec.inclusions}, one should not read too much into the name ``$(G/A')_k$,'' as there is no known functorial construction taking in a compact but non-simply connected group and producing a WZW algebra: a functorial construction is outlined in~\cite{HenriquesWZW}, but would require stronger results than are available relating loop groups and VOAs; the rigorous construction in~\cite{HenriquesWZW} is via a case-by-case analysis. Furthermore, it is worth emphasizing that the simply connected group $G$ acts, but usually non-faithfully, on the WZW algebra $G_k$, and the quotient group $G/A'$ may fail to act on the VOA $(G/A')_k$.

  The one exception to the rule ``$A = Z(G)$'' is $G_k = \rE_{8,2}$~\cite{MR1096120}, which has an abelian anyon (corresponding to the $3875$-dimensional $\rE_8$-module) even though $Z(\rE_8)$ is trivial.  Henriques does not consider the corresponding $\bZ_2$ extension of $\rE_{8,2}$ to be a ``WZW algebra,'' and it does not seem to have a standard name; we will abusively call it ``$\rE_{8,2}/\bZ_2$.''
For this exception, the simple current has conformal dimension $h = \frac32$.

As we have already explained, we are interested in the case when $a \in A$ has order $2$ with $h_a = \frac32$. 
Let us focus for now on the case when the group $G$ is simple.
In type $\rA$, this means that $m$ is even and $\frac{k(m/2)^2}{2m} = \frac 32$. In type $\rC$, this means $\frac{km}4 = \frac32$. In type $\rD_\ev$, we could take $a=s^\pm$ and solve $\frac{km}{16} = \frac32$. Otherwise, in all types $\rB\rD$, we can take $a = v$ and $m$ arbitrary provided $\frac k2 = \frac 32$, and for $\rE_7$ we want $\frac{3k}4 = \frac32$. 
The abelian anyon in the exceptional case $\rE_{8,2}$ has the desired conformal dimension.
All together, we find:

\begin{corollary}
The SVOAs $V$ such that $V_\ev = G_k$ is a simple simply-connected WZW algebra and $V_\odd$ has conformal dimension $\frac32$ are the following:\\[-3pt]

\hfill
$ \begin{array}[b]{l|l|l}
V & V_\ev & \text{Spin-$\frac32$ fields in }V_\odd \\[3pt] \hline & \\[-9pt]
\SO(m)_3 & \Spin(m)_3 & \text{$\Sym^3(\mathbf{m}) \ominus \mathbf{m}$, of dimension $\frac{m(m-1)(m+4)} 6$}\\[3pt]
\PSp(2{\times}3)_2 & \Sp(2{\times}3)_2 & \text{Irrep of dimension $84$} \\[3pt]
3\PSU(6)_2 & \SU(6)_2 & \text{Irrep of dimension $175$} \\[3pt]
\SO^+(12)_2 & \Spin(12)_2 & \text{Either irrep of dimension $462$}  \\[3pt]
\PSp(2{\times}6)_1 & \Sp(2{\times}6)_1 & \text{$\Alt^6(\mathbf{12}) \ominus \Alt^4(\mathbf{12})$, of dimension $429$} \\[3pt]
6\PSU(12)_1 & \SU(12)_1 & \text{$\Alt^6(\mathbf{12})$, of dimension $924$} \\[3pt]
\SO^+(24)_1 & \Spin(24)_1 & \text{Either half-spin irrep, of dimension $2048$} \\[3pt]
\rP\rE_{7,2} & \rE_{7,2} & \text{Irrep of dimension $1463$} \\[3pt]
\rE_{8,2}/\bZ_2 & \rE_{8,2} & \text{Irrep of dimension $3875$} 
\end{array}\hfill\qed
$
\end{corollary}

By ``$\mathbf{m}$'' we mean, of course, the $m$-dimensional vector representations of $\Spin(m)$, and $\mathbf{12}$ means the vector representation of $\Sp(2{\times}6)$ or $\SU(12)$. 
Their symmetric and alternating powers are not irreducible; for instance, $\Sym^3(\mathbf{m})$ splits as $\mathbf{m}$ plus an irrep that we will simply call $\Sym^3(\mathbf{m}) \ominus \mathbf{m}$ consisting of the ``traceless'' symmetric $3$-tensors. $\Sp(2{\times}6)$ has two $429$-dimensional irreps: the other one, which does not lead to a simple current extension of $\Sp(2{\times}6)_1$, is $\Alt^4(\mathbf{12}) \ominus \Alt^2(\mathbf{12})$.  Other irreps are listed by dimension. We took advantage of the famous exceptional isomorphisms of small spin groups to leave redundant entries off the list; see \S\ref{subsec:small-m}.
The name ``$\PSp$'' means the adjoint form of $\Sp$, and ``$n\PSU(m)$,'' for $n$ dividing $m$, means the perfect central extension of $\PSU(m)$ by a cyclic group of order $n$, i.e.\ $n\PSU(m) = \SU(m)/(\bZ_{m/n})$.

\subsection{Conventions for Spin groups} \label{subsec:small-m}

When $m$ is divisible by $4$, the centre of $\Spin(m)$ is a Klein-four group, and so $\Spin(m)$ has three quotients by $\bZ_2$. One of these is $\SO(m)$, defined as the image of $\Spin(m)$ in the adjoint representation. It is invariant under the outer automorphism of $\Spin(m)$. The other two are called $\SO^\pm(m)$. They are exchanged by the outer automorphism of $\Spin(m)$, and are by definition the image of $\Spin(m)$ in the half-spin representations.

The Spin groups have some well-known exceptional behaviour for small $m$. The triality automorphism of $\Spin(8)$ relates the 
vector representation with the two half-spin representations; thus the SVOAs $\SO(8)_3$, $\SO^+(8)_3$, and $\SO^-(8)_3$ are all isomorphic.
The exceptional isomorphisms 
\begin{gather*}
\Spin(3) \cong \SU(2) \cong \Sp(2{\times}1), \quad \Spin(4) \cong \SU(2)^2 \cong \Sp(2{\times}1)^2, \\ \Spin(5) \cong \Sp(2{\times}2), \quad   \Spin(6) \cong \SU(4)
\end{gather*}
mean that we have left out of the Corollary from~\S\ref{subsec.extension} the redundant entries $\Sp(2{\times}1)_6 = \SU(2)_6$ (see next paragraph),  $\Sp(2{\times}2)_3$, and $\SU(4)_3$. 

Recall that to define a simply connected WZW algebra for the group $G$ requires a choice of ``level'' $k \in \H^4(BG;\bZ)$. When $m\geq 5$, the algebra $\Spin(m)_k$ corresponds to the class $k \frac{p_1}2 \in \H^4(B\Spin(m);\bZ)$, where
 $\frac{p_1}2$, the \define{fractional first Pontryagin class}, is the positive generator of $\H^4(B\Spin(\infty);\bZ)$; it is a \define{stable} class in the sense that it is restricted from $B\Spin(\infty) = \varinjlim B\Spin(m)$. Our convention, even for $m<5$, will be to write ``$\Spin(m)_k$'' for the level given by restricting $k\frac{p_1}2$ along the standard inclusion $B\Spin(m) \subset B\Spin(\infty)$.

Unpacking this convention for $m=4$, we have for example the group $\Spin(4) \cong \Sp(2{\times}1)^2$. This group is not simple, and so not technically considered in the Corollary, but does appear in our Theorem; see \S\ref{subsec.genus}. The level restricts diagonally: the appropriate meaning of ``$\Spin(4)_3$'' is $\Sp(2{\times}1)_3^2$, with level $(3,3) \in \H^4(B\Spin(4);\bZ) = \bZ^2$.

For the  $m=3$ case of the  $\Spin(m)_3$ family,  we will always write scare-quotes around ``$\Spin(3)_3$.''  The reason for the scare-quotes is that it is not technically a ``level $3$'' algebra, at least with the usual normalization for $\H^4(B\Spin(3);\bZ)$. Indeed, the generator $\frac{p_1}2 \in \H^4(B\Spin(\infty);\bZ)$ restricts along the standard inclusion $B\Spin(3) \subset B\Spin(\infty)$ to \emph{twice} the generator of $\H^4(B\Spin(3);\bZ)$.
The generator of $\H^4(B\Spin(3);\bZ)$ does have a stable interpretation, in fact two of them: the inclusion $\SU(m) \subset \SU(m+1)$ induces an isomorphism $\H^4(B\SU(m+1);\bZ) \isom \H^4(B\SU(m);\bZ)$ when $m \geq 2$, and the generator is the \define{Chern class} $c_2$; also, for all $m$, the inclusion $\Sp(2{\times}m) \subset \Sp(2{\times}(m+1))$ induces an isomorphism on $\H^4$, and the generator of $\H^4(B\Sp(2{\times}m);\bZ)$ is  the so-called \define{quaternionic Pontryagin class}. In summary, the $m=3$ case of the $\Spin(m)_3$ family is $\text{``}\Spin(3)_3\text{''} = \SU(2)_6 = \Sp(2{\times}1)_6$, which we otherwise should have included in our list.

Finally, it is worth addressing the $m=2$ case of the $\Spin(m)_3$ family. The group $\Spin(2)$ is not simply connected, being isomorphic to $\rU(1)$, and so ruled out of out Theorem by fiat, but is important in our analysis in \S\ref{su12.uniqueness} and \S\ref{unique.su6}. 

If $G/A$ is a compact Lie group with  simply connected cover $ G$ and finite $\pi_1(G/A) = A$, then the restriction $\H^4(B(G/A);\bZ) \to \H^4(B G)$ is an injection, and it is traditional to say that $G/A$ acts with ``level $k$'' if the induced $ G$-action has level $k$. For example, writing $\rP \rE_6 = \rE_6 / \bZ_3$ for the adjoint form of $\rE_6$, the inclusion $\bZ \cong \H^4(B\rP \rE_6;\bZ) \to \H^4(B\rE_6;\bZ) \cong \bZ$ has cokernel of order $3$, and so $\rP \rE_{6,k} = (\rP \rE_6)_k$ is meaningful only when $k$ is divisible by $3$. 

This convention fails for $\rU(1)$, since the simply connected cover of $\rU(1)$ is $\bR$, and $\H^4(B\rU(1);\bZ) \to \rH^4(B\bR;\bZ)$ is not injective. Rather, ``$\rU(1)_k$'' refers to the lattice SVOA corresponding to the lattice $\sqrt{k}\bZ \subset \bR$; it is a bosonic VOA only when $k$ is even.
The reason for this convention corresponds to deciding that the vector representation of $\rU(m)$, which restricts to a level-$1$ representation of $\SU(m)$, should have ``level $1$'' on all of $\rU(m)$.
 In particular,  ``$\rU(1)$ at level $2$'' refers to the generator $-c_1^2 \in \H^4(B\rU(1))$, and ``level~1'' refers to a supercohomology class that could be called ``$-\frac{c_1^2}2$.'' 
 Supercohomology will be discussed further in \S\ref{orbifolds}; see also \cite[\S5.3--5.4]{MR3978827} for definitions and relations to invertible phases of matter, and \cite[\S1.4]{JFT} for a supercohomological interpretation of the fractional Pontryagin class $\frac{p_1}2$.
 
All together, the fractional Pontryagin class $\frac{p_1}2$, which restricts along $\Spin(3) = \SU(2)$ to $2c_2$, restricts further along the maximal torus $\Spin(2) \subset \Spin(3)$ to $4(-\frac{c_1^2}2)$. Thus $3\frac{p_1}2 \mapsto 12(-\frac{c_1^2}2)$, and so ``$\Spin(2)_3$'' is the bosonic VOA $\rU(1)_{12}$. But we come full circle: its simple current extension ``$\SO(2)_3$'' ends up being the SVOA $\rU(1)_3$, which is consistent with the isomorphism $\rU(1) \cong \SO(2)$.

\subsection{Elliptic genus constraints} \label{subsec.genus}

An SVOA $V$ with nontrivial fermionic part has two types of anyons. 

A \define{Neveu--Schwarz (NS) sector anyon} is simply a ``vertex module'' for $V$, understood internal to the category of supervector spaces: it is a supervector space $M$ with a vertex algebra action of $V$; the locality axiom describing how odd elements of $V$ act on odd elements of $M$ takes into account the Koszul sign rule. The vacuum module $V$ itself is an NS-sector anyon, and NS-sector anyons form a braided monoidal ``supercategory'' (a monoidal category enriched in supervector spaces, with a braiding whose axiomatics take into account the Koszul sign rule).

The other type of anyons are called \define{Ramond (R) sector}. In terms of vertex (super) modules, an R-sector anyon for $V$ is a ``twisted module,'' where the ``twisting'' is the canonical parity-reversal automorphism $(-1)^f$. The R-sector anyons also form a supercategory, but it is not monoidal. Rather, the R-sector is a ``module category'' for the NS-sector. Any supercategory has an underlying bosonic category given by forgetting the odd morphisms. (For instance, the underlying bosonic category of $\cat{SVec}$ itself is $\cat{Vec}^2$.) The direct sum of underlying bosonic categories of the NS- and R-sectors for $V$ is the category of (bosonic) representations of the even subalgebra $V_\ev$. A $V_\ev$-anyon $M_\ev$ produces a $V$-anyon $M$ whose even subspace is $M_\ev$ and whose odd subspace is $M_\odd = M_\ev \otimes_{V_\ev} V_\odd$ (where ``$\otimes_{V_\ev}$'' denotes fusion of $V_\ev$-anyons). The $V$-anyon $M$ is in the NS- or R-sector according to the \define{charge} of $M_\ev$, controlled by the braiding of $V_\odd$ with $M_\ev$. When $V_\ev = G_k \neq E_{8,2}$, this charge is precisely the action on $M_\ev$ of the central element of $G$ corresponding to $V_\odd$.

A (super) $V$-module $M$, whether NS- or R-sector, has furthermore two different ``characters,'' i.e.\ graded dimensions. One, the \define{NS-character}, is the graded dimension
$$ \chi_{\mathrm{NS}}(M) = \tr_M\bigl(q^{L_0 - c/24}\bigr)$$
 of the underlying non-super vector space of $M$. The \define{R-character}, on the other hand, is the graded superdimension 
$$ \chi_{\mathrm{R}}(M) = \operatorname{str}_M\bigl(q^{L_0 - c/24}\bigr) = \tr_M\bigl( (-1)^f q^{L_0 - c/24}\bigr).$$
Zhu's famous modularity result  \cite{MR1317233} in the super case asserts that the ``NS-NS'' characters, i.e.\ the NS characters of the NS-sector anyons, form a vector valued modular form for $\Gamma_0(2)$, whereas the ``R-R'' characters form a vector valued modular form for the whole modular group $\mathrm{SL}_2(\bZ)$ \cite{MR2175996,MR2681777,MR3077918,MR3205090}.

Suppose now that the SVOA $V$ is equipped with an $N{=}1$ superconformal vector $\tau$. As with any odd operator, the Fourier expansions of $\tau(z)$ on NS- and R-sector anyons have different gradings. As a result, NS-sector anyons become representations of the ``$N{=}1$ Neveu--Schwarz algebra'' whereas R-sector anyons become representations of the ``$N{=}1$ Ramond algebra.'' In the latter algebra, but not the former, the shifted energy operator $L_0-c/24$ has an odd square root (traditionally called ``$G_0$''). This operator pairs bosonic and fermionic fields in the same R-sector module, other than the \define{ground states} with $L_0 = c/24$, and so those fields cancel out of the R-R character, and we recover the famous boson-fermion cancellation:

\begin{lemma}
  If $V$ is an $N{=}1$ SVOA, then for every R-sector anyon $M$, its R-character evaluates to an integer (i.e.\ it has no $q$-dependence), equal to the signed count of ground states in $M$. In particular, if $M_\ev$ is a nontrivially-charged $V_\ev$-anyon, then $M_\ev$ and $M_\ev \otimes_{V_\ev} V_\odd$ must have the same graded dimensions up to a constant. \qed
\end{lemma}

We will use this Lemma to rule out most semisimple but non-simple groups from consideration in our Theorem.

Suppose that $G_{(1)},G_{(2)},\dots$ are simple simply connected groups, and that $G = G_{(1)} \times G_{(2)} \times \dots$ is equipped with the level $k = (k_{(1)},k_{(2)},\dots)$, and that $c = (c_{(1)},c_{(2)},\dots) \in G$ is a central element of order $2$ such that $V = (G / \{1,c\})_k$ admits an $N{=}1$ structure. Then each $c_{(i)} \in G_{(i)}$ must be nontrivial: else $\tau$ would have trivial OPE with $(G_{(i)})_{k_{(i)}} \subset V$ and so could not generate the conformal vector. Write $M_{(i)}$ for the abelian $(G_{(i)})_{k_{(i)}}$-anyon corresponding to $c_{(i)}$, so that $V_\odd = \bigotimes_i M_{(i)}$. The total conformal dimension of $V_\odd$, equal to $\frac32$ by assumption, is the sum of the conformal dimensions of the $M_{(i)}$s. Inspecting the Lemma from \S\ref{subsec.extension}, we see find the following choices for conformal dimension $<3/2$:
$$ \begin{array}{l|l|l}
h_a & (G_{(i)})_{k_{(i)}} & \dim(M_{(i)}) \\ \hline 
1/4 & \Sp(2{\times}1)_1 & 2 \\
\hline 1/2 & \Spin(m)_1 & m \\
\hline 3/4 & \Sp(2{\times}1)_3 & 4 \\
 & \Sp(2{\times}3)_1 & 14 \\
 & \SU(6)_1 & 20 \\
 & \Spin(12)_1 & 32 \\
 & \rE_{7,1} & 56 \\
\hline 1 & \Spin(m)_2 & \frac{(m+2)(m-1)}2 \\
 & \Sp(2{\times}4)_1 & 42 \\
 & \Spin(16)_1 & 128 \\
\hline 5/4 & \Sp(2{\times}1)_5 & 6 \\
 & \Sp(2{\times}5)_1 & 132 \\
 & \Spin(20)_1 & 512
\end{array}$$
Note that $\Sp(2{\times}2)_k = \Spin(5)_k$, $\SU(4)_k = \Spin(6)_k$, and $\SU(2)_{2k} = \Sp(2{\times}1)_{2k} = \text{``}\Spin(3)_{k}\text{''}$ do not require separate entries.

For example, $V_\ev = \Sp(2{\times}3)_1 \times \Spin(12)_1$ has an anyon $V_\odd = M_{(1)} \otimes M_{(2)}$ of conformal dimension~$\frac32$. Is it possible for the corresponding SVOA $V$ to admit an $N{=}1$ structure? No. Indeed, writing $\vac$ for vacuum module of either $\Sp(2{\times}3)_1$ or $\Spin(12)_1$, we find that the $V_\ev$-modules $M_{(1)} \boxtimes \vac$ and $\vac \boxtimes M_{(2)}$ are nontrivially charged and exchanged by fusion with $V_\odd$ and hence merge into an R-sector $V$-anyon, but that their characters are
$$ \chi\bigl(M_{(1)} \boxtimes \vac\bigr) = 14 \, q^{0.325} + \dots, \quad \chi\bigl(\vac \boxtimes M_{(2)}\bigr) = 32 \, q^{0.325} + \dots.$$
The coefficients are simply the dimensions of the minimal-spin subspaces of the anyons $M_{(i)}$, and the power of $q$ is $\frac34 - \frac c{24}$ (a quantity that Schellekens' computer algebra program ``Kac'' calls the ``modular anomaly''). Since $14\neq 32$, we have found a non-constant R-R character, and so this SVOA $V$ does not admit an $N{=}1$ structure.
The same conclusion holds whenever $G_k = (G_{(1)})_{k_{(1)}} \times (G_{(2)})_{k_{(2)}}$ with the $M_{(i)}$s of conformal dimensions $(\frac34,\frac34)$ or $(\frac14,\frac54)$: we must have $(G_{(1)})_{k_{(1)}} = (G_{(2)})_{k_{(2)}}$.

In the $\Spin(m)_1 \times \Spin(m')_2$ and $\Spin(m)_1 \times \Sp(2{\times}4)_2$ cases, the $V_\ev$-anyons $M_{(1)} \boxtimes \vac$ and $\vac \boxtimes M_{(2)}$ merge to an NS-sector $V$-anyon, where $V_\odd = M_{(1)} \boxtimes M_{(2)}$. To build an R-sector $V$-anyon, one can take start with the (not necessarily irreducible) $\Spin(m)_1$-anyon $N$ of conformal dimension $h = \frac m {16}$ 
 whose minimal-spin fields are the full spinor representation (of dimension $2^{m/2}$); then $N \boxtimes \vac$ merges with $N \boxtimes M_{(2)}$ to form an R-sector anyon with nonconstant R-R character.
  
Similar arguments also rule out more complicated combinations with three or more simple factors. All together, we see that, in order for $G_k = (G_{(1)})_{k_{(1)}} \times (G_{(2)})_{k_{(2)}} \times \dots$ to have a $\bZ_2$-extension with $N{=}1$ supersymmetry and no free fermions, a necessary condition is that all the simple factors $(G_{(i)})_{k_{(i)}}$ must be isomorphic, with one exception: $\Spin(16)_1 \times \Spin(8)_1$ extended by its anyon of conformal dimension $\frac32$ passes the test in the previous Lemma.

\begin{corollary}
  The SVOAs $V$ with no free fermions, $V_\ev = G_k$ a simply connected WZW algebra, $V_\odd$ of conformal dimension $\frac32$, and for which all R-R characters are constants are those listed in the Corollary in \S\ref{subsec.extension} (including $\Spin(4)_3 = \Sp(2{\times}1)_3^2$) and 
  $$
   \Spin(m)_1^3\quad \text{(including $\Spin(3)_1^3 = \Sp(2{\times}1)_2^3$ and $\Spin(4)_1^3 = \Sp(2{\times}1)_1^6$)}, $$
\quad\hfill   $\Sp(2{\times}3)_1^2,\quad \SU(6)_1^2, \quad \Spin(12)_1^2, \quad \Spin(16)_1 \times \Spin(8)_1, \quad \text{and} \quad \rE_{7,1}^2$. \hfill \qed
\end{corollary}

\subsection{Inclusions of non-simply connected WZW algebras} \label{subsec.inclusions}

The (bosonic) simply connected WZW algebras $V_\ev = G_k$ listed in the Corollaries in \S\ref{subsec.extension} and \S\ref{subsec.genus} are related by the following inclusions:
$$ \begin{tikzpicture}[xscale=1.1]
  \path 
    (0,1) node (A16) {\!``$\Spin(3)_3$''\!}
    +(0,-.5) node[rotate=90] {$=$}
    +(0,-1) node (C16) {$\Sp(2{\times}1)_6$}
    ++(2,0) node (D23) {$\Spin(4)_{3}$}
    +(0,-.5) node[rotate=90] {$=$}
    +(0,-1) node (C132) {$\Sp(2{\times}1)_3^2$}
    ++(2,0) node (B23) {$\Spin(5)_3$}
    +(0,-.5) node[rotate=90] {$=$}
    +(0,-1) node (C23) {$\Sp(2{\times}2)_3$}
    ++(2,0) node (D33) {$\Spin(6)_3$}
    +(0,-.5) node[rotate=90] {$=$}
    +(0,-1) node (A33) {$\SU(4)_3$}
    ++(2,0) node (B33) {$\Spin(7)_3$}
    ++(2,0) node (D43) {$\Spin(8)_3$}
    +(0,-.5) node[rotate=90] {$=$} +(0,-.5) node[anchor=west] {\tiny triality}
    +(0,-1) node (D43t) {$\Spin(8)_3$}
    ++(2,0) node (B43) {$\Spin(9)_3$}
    ++(1.5,0) node (dots) {$\dots$}
  ;
  \path 
    (0,2) node (A123) {\!``$\Spin(3)_1^3$''\!} 
    +(0,.5) node[rotate=90] {$=$}
    +(0,1) node (C123) {$\Sp(2{\times}1)_2^3$}
    ++(2,0) node (D213) {$\Spin(4)_1^3$}
    +(0,.5) node[rotate=90] {$=$}
    +(0,1) node (C116) {$\Sp(2{\times}1)_1^6$}
    ++(2,0) node (B213) {$\Spin(5)_1^3$}
    +(0,.5) node[rotate=90] {$=$}
    +(0,1) node (C213) {$\Sp(2{\times}2)_1^3$}
    ++(2,0) node (D313) {$\Spin(6)_1^3$}
    +(0,.5) node[rotate=90] {$=$}
    +(0,1) node (A313) {$\SU(4)_1^3$}
    ++(2,0) node (B313) {$\Spin(7)_1^3$}
    ++(2,0) node (D413) {$\Spin(8)_1^3$}
    +(0,.5) node[rotate=90] {$=$} +(0,.5) node[anchor=west] {\tiny triality}
    +(0,1) node (D413t) {$\Spin(8)_1^3$}
    ++(2,0) node (B413) {$\Spin(9)_1^3$}
    ++(1.5,0) node (dots1) {$\dots$}
  ;
  \path
    (0,4) node (C32) {$\Sp(2{\times}3)_2$}
    ++(2,1) node (A52) {$\SU(6)_2$}
    ++(2,-1) node (C61) {$\Sp(2{\times}6)_1$}
    ++(0,2) node (D62) {$\Spin(12)_2$}
    ++(2,-1) node (A111) {$\SU(12)_1$}
    ++(4,-1) node (D8D4) {$\Spin(16)_1 \times \Spin(8)_1$}
    ++(0,2) node (D121) {$\Spin(24)_1$}
  ;
  \path
    (D62)
    ++(0,1) node (E72) {$\rE_{7,2}$}
    ++(0,1) node (E82) {$\rE_{8,2}$}
  ;
  \path
    (C32) ++(2,0) node (C312) {$\Sp(2{\times}3)_1^2$}
    (A52) ++(2,0) node (A512) {$\SU(6)_1^2$}
    (D62) ++(2,0) node (D612) {$\Spin(12)_1^2$}
    (E72) ++(2,0) node (E712) {$\rE_{7,1}^2$}
  ;
  \draw[right hook - >] (A16) -- (D23);
  \draw[right hook - >] (A123) -- (D213);
  \draw[right hook - >] (A16) -- (A123);
  \draw[right hook - >] (D23) -- (B23);
  \draw[right hook - >] (B23) -- (D33);
  \draw[right hook - >] (D33) -- (B33);
  \draw[right hook - >] (B33) -- (D43);
  \draw[right hook - >] (D43) -- (B43);
  \draw[right hook - >] (B43) -- (dots);
  \draw[right hook - >] (D23) -- (D213);
  \draw[right hook - >] (B23) -- (B213);
  \draw[right hook - >] (D33) -- (D313);
  \draw[right hook - >] (B33) -- (B313);
  \draw[right hook - >] (D43) -- (D413);
  \draw[right hook - >] (B43) -- (B413);
  \draw[right hook - >] (D213) -- (B213);
  \draw[right hook - >] (B213) -- (D313);
  \draw[right hook - >] (D313) -- (B313);
  \draw[right hook - >] (B313) -- (D413);
  \draw[right hook - >] (D413) -- (B413);
  \draw[right hook - >] (B413) -- (dots1);
  \draw[right hook - >] (C16) -- (C132);
  \draw[right hook - >] (C132) -- (C23);
  \draw[right hook - >] (C23) -- (A33);
  \draw[right hook - >] (A33) -- (D43t);
  \draw[right hook - >] (C123) -- (C116);
  \draw[right hook - >] (C116) -- (C213);
  \draw[right hook - >] (C213) -- (A313);
  \draw[right hook - >] (A313) -- (D413t);
  \draw[right hook - >] (C123) -- (C32);
  \draw[right hook - >] (C116) -- (C312);
  \draw[right hook - >] (C213) -- (C61);
  \draw[right hook - >] (A313) -- (A111);
  \draw[right hook - >] (C32) -- (C312);
  \draw[right hook - >] (C312) -- (C61);
  \draw[right hook - >] (C32) -- (A52);
  \draw[right hook - >] (C312) -- (A512);
  \draw[right hook - >] (C61) -- (A111);
  \draw[right hook - >] (A52) -- (C61);
  \draw[right hook - >] (A52) -- (A512);
  \draw[right hook - >] (A512) -- (A111);
  \draw[right hook - >] (A52) -- (D62);
  \draw[right hook - >] (A512) -- (D612);
  \draw[right hook - >] (A111) -- (D121);
  \draw[right hook - >] (D62) -- (A111);
  \draw[right hook - >] (D62) -- (D612);
  \draw[right hook - >] (D612) -- (D121);
  \draw[right hook - >] (D413t) -- (D8D4);
  \draw[right hook - >] (D8D4) -- (D121);
  \draw[right hook - >] (D62) -- (E72);
  \draw[right hook - >] (D612) -- (E712);
  \draw[right hook - >] (E72) -- (E712);
  \draw[right hook - >] (E72) -- (E82);
\end{tikzpicture} $$
Indeed, let $G' \subset G$ be an inclusion of  simply connected Lie groups. Then $G'_{k'} \subset G_k$ exactly when $k'$ is the image of $k$ under the restriction map $\H^4(BG;\bZ) \to \H^4(BG';\bZ)$. (When $G'$ and $G$ are simple, these cohomology groups are both isomorphic to $\bZ$, and the restriction map is multiplication by a positive integer called the \define{Dynkin index} of the inclusion.)
This is how the above diagram was charted.

Except in special situations, the inclusion $G'_{k'} \subset G_k$ is not \define{conformal} --- it does not intertwine conformal vectors. Rather, there is a \define{coset} VOA $G_k / G'_{k'}$ defined to consist of those fields in $G_k$ which commute with $G'_{k'}$. Writing $\nu$ and $\nu'$ for the conformal vectors of $G_{k}$ and $G'_{k'}$ respectively, the conformal vector for the coset is $\nu'' = \nu - \nu'$, so that the total inclusion $G'_{k'} \boxtimes (G_k / G'_{k'}) \subset G_k$ is conformal. The coset of a conformal inclusion is the trivial VOA $\bC$. Often but not always a sub-VOA is equal to its double coset; for example, this fails for conformal inclusions, but holds in many examples related to level-rank duality. Complete lists of conformal inclusions of simply connected WZW algebras are available in \cite{MR867023,MR867243}, and the relation to level-rank duality is discussed in \cite{MR3039775}.

One should expect that the non-simply connected WZW algebras also depend functorially on the corresponding groups, but this is not manifest from the construction and does not seem to be known in general (but see \cite[Section~2]{HenriquesWZW} for a heuristic construction of non-simply connected WZW algebras which is manifestly functorial). If we did have such functoriality, then it would follow that all of the above bosonic inclusions, other than $\rE_{7,2} \subset \rE_{8,2}$, extend to the $\bZ_2$-extensions: the reader is invited to check that in all cases, the given $\bZ_2$-quotients of the simply connected groups do map appropriately. Since we do not have a general functorialty result, we will instead check that for the SVOAs from  \S\ref{subsec.extension}, the bosonic inclusions extend.
We check this as follows. 

Let $G' \subset G$ be an inclusion of simple simply connected Lie groups from the above chart, 
so that we are trying to establish an inclusion $G'_{k'}/\bZ_2 \subset G_k/\bZ_2$ of SVOAs. We do have inclusions $G'_{k'} \subset G_k \subset G_k/\bZ_2$. Write $M = (G_k/\bZ_2)_\odd$ for the $G_k$-anyon consisting of the fermionic part of $G_k/\bZ_2$; it has conformal dimension $h_M = \frac32$, and its space $M_{3/2}$ of fields of spin $\frac32$ is the simple $G$-module listed, in the simple-$G$ case, in the Corollary in \S\ref{subsec.extension}. Now decompose $M$ over $G'_{k'} \times (G_k / G'_{k'})$. It decomposes as
$$ M = \bigoplus A_i \boxtimes B_i $$
where the $A_i$s are simple $G'_{k'}$-anyons and the $B_i$s are simple $(G_k / G'_{k'})$-anyons. They have complementary conformal dimensions: $h_{A_i} + h_{B_i} = h_M = \frac32$. Since $G'_{k'}$ is  simply connected, the fields of minimal spin in any simple anyon $A$ are a simple module for the compact group $G'$, and $A$ is determined by this module. In the case at hand, these simple modules are precisely the modules appearing as direct summands inside the restriction $M_{3/2}|_{G'}$ of $M_{3/2}$ to $G'$.

Write $M' = (G'_{k'}/\bZ_2)_\odd$ for the fermionic part of $G'_{k'}/\bZ_2$, and $M'_{3/2}$ for its space of spin-$\frac32$ fields. By the remarks in the previous paragraph, $M'$ appears as one of the $A_i$s if and only if $M'_{3/2}$ appears as a direct summand of $M_{3/2}|_{G'}$. If it does, then the corresponding $B_i$ must be the vacuum $(G_k / G'_{k'})$-module, since it must have conformal dimension $0$.
But then
\begin{multline*}
 G'_{k'}/\bZ_2 \boxtimes (G_k / G'_{k'}) = (G'_{k'} \oplus M') \boxtimes (G_k / G'_{k'}) \\ = \bigl(G'_{k'} \boxtimes (G_k / G'_{k'})\bigr) \oplus \bigl(M' \boxtimes (G_k / G'_{k'})\bigr) \subset G_k \oplus M = G_k/\bZ_2
\end{multline*}
and the inclusion is verified. And sure enough:
\begin{lemma}
  For all inclusions $G' \subset G$ in the above diagram except for $\rE_{7,2} \subset \rE_{8,2}$, $M_{3/2}|_{G'}$ contains $M'_{3/2}$ as a direct summand, and so the inclusion extends to an SVOA inclusion $G'_{k'}/\bZ_2 \subset G_k/\bZ_2$.
\end{lemma}

It is worth warning that, although many of the inclusions $G'_{k'}/\bZ_2 \subset G_k/\bZ_2$ are compatible with $N{=}1$ superconformal structures, not all of them are.  We will describe the modules $M_{3/2}$ in more detail in Section~\ref{sec.existence}, and we will
 use some but not all of the resulting inclusions in Sections~\ref{sec.spinm3} and~\ref{sec.uniqueness}.

\begin{proof} For the inclusions $G_2 \subset G_1^2$, write $N$ for the abelian $G_1$-anyon of conformal dimension $\frac34$. Then $M_{3/2} = N^{\boxtimes 2}$ restricts over the diagonal inclusion $G \subset G^2$ to $N^{\otimes 2}$, and its highest-weight submodule is $M'_{3/2}$. For the remaining inclusions, we work from bottom to top and right to left:
  \begin{description}
    \item[$\Spin(m)_3 \subset \Spin(m+1)_3$] The module $M_{3/2}$ is the space $\Sym^3(\mathbf{m+1}) \ominus (\mathbf{m+1})$ of traceless symmetric 3-tensors in $m+1$ variables. It includes the space $M'_{3/2} = \Sym^3(\mathbf{m}) \ominus \mathbf{m}$ of traceless 3-tensors in $m$ variables.
    \item[$\Spin(m)_3 \subset \Spin(m)_1^3$] The module $M_{3/2}$ is the ``outer'' product $\mathbf{m}^{\boxtimes 3}$ of the vector representations of the three copies of $\Spin(m)$. It restricts along the diagonal inclusion $G' = \Spin(m) \subset \Spin(m)^3 = G$ to $\mathbf{m}^{\otimes 3}$, which contains the irrep $\Sym^3(\mathbf{m}) \ominus \mathbf{m}$ as a direct summand.
    \item[$\Spin(m)_1^3 \subset \Spin(m+1)_1^3$] $(\mathbf{m+1})^{\boxtimes 3} \supset \mathbf{m}^{\boxtimes 3}$. 
    \item[$\Spin(8)_1^3 \subset \Spin(16)_1 \times \Spin(8)_1$] $M_{3/2}$ is the tensor product of positive half-spin modules $\mathbf{128}_+ \boxtimes \mathbf{8}_+$. It decomposes over $\Spin(8)^3$ as $(\mathbf{8}_+ \boxtimes \mathbf{8}_+ \boxtimes \mathbf{8}_+) \oplus (\mathbf{8}_- \boxtimes \mathbf{8}_- \boxtimes \mathbf{8}_+)$. The first summand becomes $M'_{3/2} = \mathbf{8}^{\boxtimes 3}$ after applying the triality automorphism of $\Spin(8)$.
    \item[$\Spin(16)_1 \times \Spin(8)_1 \subset \Spin(24)_1$] $M_{3/2}$ is the positive half-spin module $\mathbf{2048}_+$. It decomposes over $\Spin(16) \times \Spin(8)$ as $(\mathbf{128}_+ \boxtimes \mathbf{8}_+) \oplus (\mathbf{128}_- \boxtimes \mathbf{8}_-)$. The first summand is $M'_{3/2}$.
    \item[$\SU(4)_1^3 \subset \SU(12)_1$] $M_{3/2} = \Alt^6(\mathbf{12})$. 
    Write $\mathbf{4}^{\oplus 3} = \mathbf{12}|_{\SU(4)^3}$ for the direct sum of the three vector representations. Then $M_{3/2}$
     restricts over $\SU(4)^3$ as $\Alt^6(\mathbf{4}^{\boxplus 3}) \supset \Alt^2(\mathbf{4})^{\boxtimes 3}$.
    \item[$\Sp(2{\times}2)_1^3 \subset \Sp(2{\times}6)_1$]  $M_{3/2} = \mathbf{429} = \Alt^6(\mathbf{12}) \ominus \Alt^4(\mathbf{12})$. The vector representation $\mathbf{12}$ restricts along the diagonal inclusion $\Sp(2{\times}2)^3 \subset \Sp(2{\times}6)$ to $\mathbf{4}^{\boxplus 3}$, and so $\mathbf{429}$ contains $M'_{3/2} = \bigl(\Alt^2(\mathbf 4) \ominus \mathbf 1\bigr)^{\boxtimes 3}$.
    \item[$\Sp(2{\times}1)_1^6 \subset \Sp(2{\times}3)_1^2$] The $\Sp(2{\times}3)$-irrep $\mathbf{14}_- = \Alt^3(\mathbf6) \ominus \mathbf6$ restricts along $\Sp(2{\times}1)^3 \subset \Sp(2{\times}3)$ to $\Alt^3(\mathbf2^{\boxplus 3}) \ominus (\mathbf2^{\boxplus 3}) \supset \mathbf2^{\boxtimes 3}$, and so $M_{3/2} = \mathbf{14}_-^{\boxtimes 2}$ contains $M'_{3/2} = \mathbf2^{\boxtimes 6}$.
    \item[$\Sp(2{\times}1)_2^3 \subset \Sp(2{\times}3)_2$] $M_{3/2} = \mathbf{84} = \Sym^2(\mathbf{14}_-) \ominus \mathbf{21}$, where $\mathbf{21}$ is the adjoint representation of $\Sp(2{\times}3)$. As above, $\mathbf{14}_-$ restricts along $\Sp(2{\times}1)^3 \subset \Sp(2{\times}3)$ to contain $\mathbf2^{\boxtimes 3}$, and so $\Sym^2(\mathbf{14}_-)$ contains $M'_{3/2} = \Sym^2(\mathbf2)^{\boxtimes 3}$, whereas $\mathbf{21}|_{\Sp(2{\times}1)^3} = \mathbf3^{\boxplus 3}$. 
    \item[$\SU(12)_1 \subset \Spin(24)_1$] $M_{3/2} = \mathbf{2048}_+$ splits over $\SU(12)$ as $\Alt^0(\mathbf{12}) \oplus \Alt^2(\mathbf{12}) \oplus \Alt^4(\mathbf{12}) \oplus \Alt^6(\mathbf{12}) \oplus \Alt^8(\mathbf{12}) \oplus \Alt^{10}(\mathbf{12}) \oplus \Alt^{12}(\mathbf{12})$, and so contains $M'_{3/2} = \Alt^6(\mathbf{12})$ as a direct summand.
    \item[$\Sp(2{\times}6)_1 \subset \SU(12)_1$] $M_{3/2} = \Alt^6(\mathbf{12})$ contains $\Alt^6(\mathbf{12}) \ominus \Alt^4(\mathbf{12}) = M'_{3/2}$.
    \item[$\Sp(2{\times}3)_1^2 \subset \Sp(2{\times}6)_1$] $M_{3/2} = \Alt^6(\mathbf{12}) \ominus \Alt^4(\mathbf{12})$ restricts to $\Alt^6(\mathbf{6}^{\boxplus 2}) \ominus \Alt^4(\mathbf{6}^{\boxplus 2}) \supset \bigl(\Alt^3(\mathbf{6}) \ominus \mathbf{6}\bigr)^{\boxtimes 2} = M'_{3/2}$.
    \item[$\Spin(12)_1^2 \subset \Spin(24)_1$] $M_{3/2}|_{\Spin(12)^2} = \mathbf{2048}_+|_{\Spin(12)^2} = \mathbf{32}_+^{\boxtimes 2} \oplus \mathbf{32}_-^{\boxtimes 2} \supset \mathbf{32}_+^{\boxtimes 2} = M'_{3/2}$.
    \item[$\SU(6)_1^2 \subset \SU(12)_1$] $M_{3/2}|_{\SU(6)^2} = \Alt^6(\mathbf{6}^{\boxplus 2}) \supset \Alt^3(\mathbf{6})^{\boxtimes 2} = M'_{3/2}$.
    \item[$\Sp(2{\times}3)_1^2 \subset \SU(6)_1^2$] $M_{3/2} = \Alt^3(\mathbf{6})^{\boxtimes 2}$ contains $M'_{3/2} = \bigl(\Alt^3(\mathbf{6}) \ominus \mathbf{6}\bigr)^{\boxtimes 2}$.
    \item[$\Spin(12)_2 \subset \SU(12)_1$] $M_{3/2} = \mathbf{924}$ splits as $\mathbf{462}_+ \oplus \mathbf{462}_-$, where $M'_{3/2} = \mathbf{462}_+$.
    \item[$\SU(6)_2 \subset \Sp(2{\times}6)_1$] $\Alt^6(\mathbf{12})$ contains $\Sym^2(\Alt^3(\mathbf{6}))$, and the desired inclusion $M'_{3/2} \subset M_{3/2}$ arises as the traceless parts of these modules.    
    \item[$\SU(6)_2 \subset \Spin(12)_2$] $\mathbf{462}_+$ is a simple submodule of $\Sym^2(\mathbf{32}_+) = \mathbf{462}_+ \oplus \mathfrak{so}(12)$, where $\mathbf{32}_+$ is the positive half-spin module of $\Spin(12)$. It splits over $\SU(6)$ as $\mathbf{32}_+ = \Alt^1(\mathbf{6}) \oplus \Alt^3(\mathbf{6}) \oplus \Alt^5(\mathbf{6})$, and so $\Sym^2(\mathbf{32}_+)$ contains $\Sym^2(\Alt^3(\mathbf{6})) = \mathbf{175}_+ \oplus \mathfrak{su}(6)$.
    \item[$\Sp(2{\times}3)_2 \subset \SU(6)_2$] The inclusion $M'_{3/2} \subset M_{3/2}$ comes from applying $\Sym^2$ to the inclusion $\Alt^3(\mathbf6) \ominus \mathbf6 \subset \Alt^3(\mathbf6)$ and restricting to traceless parts.
    \item[$\Spin(12)_1^2 \subset \rE_{7,1}^2$] The standard representation $\mathbf{56}$ of $\rE_7$ splits over $\Spin(12)$ as $\mathbf{32}_+ \oplus 2 \otimes \mathbf{12}$, and so $M_{3/2} = \mathbf{56}^{\boxtimes 2}$ contains $M'_{3/2} = \mathbf{32}_+^{\boxtimes 2}$.
    \item[$\Spin(12)_2 \subset \rE_{7,2}$] Applying $\Sym^2$ to the restriction $\mathbf{56}|_{\Spin(12)} \supset \mathbf{32}_+$, we see that $\Sym^2(\mathbf{56}) = \mathbf{1463} \oplus \mathbf{133}$ contains $\Sym^2(\mathbf{32}_+) = \mathbf{462}_+ \oplus \mathbf{66}$, where $\mathbf{1463} = M_{3/2}$, $\mathbf{462}_+ = M'_{3/2}$, $\mathbf{133} = \mathfrak{e}_7$, and $\mathbf{66} = \mathfrak{so}(12)$.
    \item[$\rE_{7,2} \subset \rE_{8,2}$] 
    According to the GAP package ``SLA'' \cite{SLA-GAP},
    $M_{3/2} = \mathbf{3875}$ splits over $\rE_7$ as $\mathbf{1} \oplus \mathbf{1539} \oplus 2 \otimes \mathbf{56} \oplus 2 \otimes \mathbf{912} \oplus 3 \otimes \mathbf{133}$, and so does not contain $M'_{3/2} = \mathbf{1463}$. \qedhere
  \end{description}
\end{proof}

\section{Constructing $N{=}1$ structures}\label{sec.existence}

In this section we construct the $N{=}1$ structures listed in the Theorem. The underlying SVOAs~$V$ are those listed in \S\ref{subsec.extension} and \S\ref{subsec.genus}.
 Our strategy is the following. Given such a $V$, we will cleverly choose a finite subgroup $S \subset \Aut_{N{=}0}(V)$. It is worth remarking that $\Aut(G_k) = \Aut(\mathfrak{g}) = G^\adj{:}\operatorname{Out}(G)$, where $G^\adj$ is the adjoint form of $G$ and $\operatorname{Out}(G)$ is the group of Dynkin diagram automorphisms of $G$, and the colon denotes a semidirect product. The extension from $V_\ev = G_k$ to $V = V_\ev \oplus V_\odd$ involves passing to a (possibly trivial) double cover of $G^\adj$, and may break some of the $\operatorname{Out}(G)$ symmetry.
 We will choose the subgroup $S \subset \Aut_{N{=}0}(V)$ so that it fixes a unique (up to scalar) non-null spin-$\frac32$ field $\tau$, but such that 
 the adjoint representation $\mathfrak{g}$ remains simple upon restriction to $S$. The space of spin-$2$ fields in a simply connected WZW algebra $G_k$ is a submodule of $\Sym^2(\mathfrak{g})\oplus \mathfrak{g}$, with equality except in very low level. (The $\Sym^2(\mathfrak{g})$ summand consists of bilinears in the Kac--Moody generators, and the $\mathfrak{g}$-summand consists of the derivatives of the Kac--Moody generators.) So simplicity of $\mathfrak{g}$ will imply that the conformal vector $\nu$ is the only (up to scalar) $S$-fixed spin-$2$ field.
 
The generic self-OPE of $\tau$ is
$$ \tau(z)\,\tau(0) \sim \frac{c'}{z^3} + \frac{X(0)}{z}$$
for some spin-$2$ field $X$ ---
locality rules out a term of the form $\frac{Y(0)}{z^2}$, hence its omission.
There is a unique $\Aut_{N{=}0}(V)$-invariant nondegenerate inner product on the spin-$\frac32$ fields, which sets $\|\tau\|^2 = c'$, and so provided $\tau$ is not lightlike for this inner product, we can rescale $\tau$ to set $c' = \frac23 c$.
Since $S$ is finite and the spin-$\frac32$ fields have a nondenegerate inner product, the $S$-module of spin-$\frac32$ fields has a natural real form for which that inner product is positive definite (namely, the real subspace on which the bilinear inner product agrees with the Hermitian form coming from unitarity of the SVOA). If $S$ has a complex fixed point, then it has a real fixed point, and so if $\tau$ was the unique fixed point, it must be real and cannot be lightlike.
Finally, $X$ is spin-$2$ and $S$-fixed and so proportional to $\nu$, and the Jacobi identity will force the correct normalization \cite[Lemma 5.9]{MR1651389}.

We now proceed with the examples. We will follow the ATLAS~\cite{ATLAS} for names for finite groups: for instance, ``$2$'' means the cyclic group $\bZ_2$, $\rU_3(3)$ means the simple subquotient of the third unitary group over $\bF_3$,  $\Suz$ is Suzuki's sporadic group, and a colon denotes a split extension.
We rely on GAP \cite{GAP} (and its CTblLib library \cite{CTblLib}, which uses the ATLAS naming conventions) for all character table computations; groups in the ATLAS were accessed through the GAP package AtlasRep \cite{AtlasRep-Gap}.

\subsection{Existence for $\Spin(m)_3$}\label{exist.spinm}

When $m\neq8$,
the automorphism group of the bosonic WZW algebra $V_\ev = \Spin(m)_3$ is the projective orthogonal group $\mathrm{PO}(m) = \rO(m)/\{\pm1\}$, isomorphic to $\PSO(m){:}2$ when $m$ is even and to $\SO(m)$ when $m$ is odd. This lifts to the full orthogonal group $\rO(m)$ upon including the simple current $V_\odd$.
When $m=8$, $\Aut(V_\ev) = \PSO(8){:}S_3$ includes the triality automorphism, but $\Aut_{N{=}0}(V)$ does not and remains $\rO(8)$.
 Choose $S = A_{m+1}$ the alternating group, embedded into $\rO(m)$ via the simple submodule of the permutation representation. The adjoint representation $\mathfrak{so}(m) = \Alt^2(\mathbf{m})$ remains simple upon restriction to $S$, and so the conformal vector $\nu$ is the unique $S$-fixed spin-$2$ field.

We claim that $A_{m+1}$ has a unique fixed point in $\Sym^3(\mathbf{m}) \ominus (\mathbf{m})$, or equivalently in $\Sym^3(\mathbf{m})$. Indeed, we can witness a fixed point explicitly: write $\mathbf{m}$ as the span of elements $e_0,\dots,e_m$, permuted by $A_{m+1}$, modulo $\sum e_i = 0$; then $\sum e_i^3$ is a nonzero $S_{m+1}$-fixed element of $\Sym^3(\mathbf{m})$. On the other hand, a Young diagram computation confirms that $\hom_{A_{m+1}}(\mathbf{m},\Sym^2(\mathbf{m}))$ is one-dimensional, and so $\Sym^3(\mathbf{m})$ has at most one fixed point.

We remark that this $\tau$ is in fact preserved by the whole symmetric group $S_{m+1} = A_{m+1}{:}2$. Actually, $S_{m+1}$ maps to $\rO(m)$ in two different ways, related by the sign representation, and $\tau$ is fixed for one of these and \define{antifixed} (i.e.\ acted on by the sign representation) for the other.

As we remarked already in \S\ref{subsec:small-m}, special cases include $\Spin(6)_3 = \SU(4)_3$, $\Spin(5)_3 = \Sp(2{\times}4)_3$, $\Spin(4)_3 = \Sp(2{\times}1)_3^2 = \SU(2)_3^2$, and $\text{``}\Spin(3)_3\text{''} = \Sp(2{\times}1)_6 = \SU(2)_6$. In \S\ref{su12.uniqueness} and \S\ref{unique.su6}, we will use the non-simply-connected special case $\text{``}\Spin(2)_3\text{''} = \rU(1)_{12}$.

\subsection{Existence for $\Spin(m)_1^3$}\label{exist.spinm.3}

The automorphism group of $V_\ev = \Spin(m)_1^3$ contains (and is equal to, when $m\neq 4,8$) the wreath product $\rP\rO(m)\wr S_3 = \rP\rO(m)^3{:}S_3$ --- the $S_3$ permutes the three simple factors --- and lifts to the double cover $2.\rP\rO(m)^3{:}S_3 = \rO(m)^{\circ 3}{:}S_3$ upon inclusion of the simple current $V_\odd$. (If $G$ and $G'$ have identical centres $Z = Z(G) = Z(G')$, then we will write $G \circ G'$ to denote their \define{central product} $(G \times G') / Z$.
As in \S\ref{exist.spinm}, when $m=8$, $\Aut(V_\ev)$ also includes the triality automorphism, but $\Aut_{N{=}0}(V)$ does not.)
 The space of spin-$\frac32$ fields is $\mathbf{m}^{\boxtimes 3}$, the outer product of three copies of the vector representation of $\rO(m)$. Choose an orthonormal basis $e_1,\dots,e_m$ for $\mathbf{m}$; then a basis for $\mathbf{m}^{\boxtimes 3}$ is $e_{\vec\imath} = e_{i_1} \boxtimes e_{i_2} \boxtimes e_{i_3}$, where $\vec\imath = (i_1,i_2,i_3) \in \{1,\dots,m\}^3$. 

Embed $S = 2^{2(m-1)}{:}(S_3 \times S_m)$ into $\rO(m)^{\circ 3}{:}S_3$ as follows. The $S_3$ embeds along the $S_3$ permuting the three factors. The $S_m$ embeds diagonally inside $\rO(m) \subset \rO(m)^{\circ 3}$ via the permutation representation. Finally, write $2^{2m} = (2^2)^m$, and declare that the $j$th copy of the Klein-4 group $2^2$ embeds into $\rO(m)^3$ so as to switch the signs of two out of the three ``$j$th'' coordinates; this copy of $2^{2m} \subset \rO(m)^3$ projects to $2^{2(m-1)} \subset \rO(m)^{\circ 3}$.
 (For instance, the first nontrivial element in the first copy of $2^2 \subset (2^2)^m$ maps to $\bigl(\diag(-1,1,1,\dots), \diag(-1,1,1,\dots), \id \bigr) \in \rO(m)^3$.)  
Then $\mathbf{m}^{\boxtimes 3}$ has a unique $S$-fixed vector, namely
$$ \tau = \sum_i e_{i,i,i}.$$
On the other hand, the adjoint representation $\mathfrak{so}(m)^3$ of $\rO(m)^{\circ 3}{:}S_3$ remains simple upon restriction to $S$.

\subsection{Existence for $\Sp(2{\times}3)_2$}\label{exist.sp3}
  
  The automorphism group of $V_\ev = \Sp(2{\times}3)_2$ is the adjoint form $\PSp(2{\times}3)$. This extends to $\Aut_{N{=}0}(V) = 2 \times \PSp(2{\times}3)$ when $V_\odd$ is included --- the central extension splits because the level is even. 
  The group $S = \rU_3(3)$ has a unique six-dimensional irrep. It is quaternionic (i.e.\ its Frobenius--Schur indicator is $-1$), and so determines an injection $S \subset \Sp(2{\times}3)$. Since $S$ has no centre, this injection descends to an injection $S \subset \PSp(2{\times}3)$. The adjoint representation $\mathbf{21} = \mathfrak{sp}(3)$ of $\PSp(2{\times}3)$ remains simple upon restriction to $S$.
  
   Denote the three fundamental representations of $\Sp(2{\times}3)$ by $\mathbf{14}_+$, $\mathbf{14}_-$, and $\mathbf{6}$. The third is the standard representation, and the other two are
   $$ \mathbf{14}_+ = \Alt^2(\mathbf{6}) \ominus \mathbf{1}, \qquad \mathbf{14}_- = \Alt^3(\mathbf{6}) \ominus \mathbf{6}.$$
   $\mathbf{14}_+$ is real and descends to $\PSp(2{\times}3)$ whereas $\mathbf{14}_-$ is quaternionic and acted on nontrivially by the centre of $\Sp(2{\times}3)$. 
  The $84$-dimensional space of spin-$\frac32$ fields appears inside
  $$ \Sym^2(\mathbf{14}_-) = \mathbf{21} \oplus \mathbf{84}.$$

The restriction of $\mathbf{14}_+$ to $S$ is irreducible, but the restriction of $\mathbf{14}_-$ breaks as $\mathbf{7} \oplus \overline{\mathbf{7}}$ for a dual pair of  seven-dimensional complex irreps. ($S$ also has a real seven-dimensional irrep.)  Since the restriction of $\mathbf{21}$ to $S$ is simple, it follows that $\Sym^2(\mathbf{14}_-)$, and hence $\mathbf{84}$, has a unique $S$-fixed point, providing us with our desired spin-$\frac32$ field~$\tau$.
  
  We remark that $\PSp(2{\times}3) \times 2$ contains two conjugacy classes of $\rU_3(3){:}2$ subgroups, whose actions on $\mathbf{84}$ differ by the sign representation. Indeed, 
  write $f : \rU_3(3){:}2 \mono \PSp(2{\times}3)$ for the unique conjugacy class of embeddings therein; then the other copy of $\rU_3(3){:}2 \subset \PSp(2{\times}3) \times 2$ is given by the map $(f, \pi)$, where $\pi : \rU_3(3){:}2 \to 2$ is the projection. It follows that $\tau$ is fixed by one (and not the other) of the two copies of $\rU_3(3){:}2 \subset \PSp(2{\times}3) \times 2$. 
  The other copy \define{antifixes} $\tau$ in the sense that it acts on $\tau$ via the unique nontrivial one-dimensional ``sign'' representation of $\rU_3(3){:}2$ (namely, the composition of the projection $\pi$ with the standard ``sign'' representation $2 \isom \rO(1)$).
  
\subsection{Existence for $\Sp(2{\times}3)_1^2$}\label{exist.sp3.2}

The automorphism group of $V_\ev = \Sp(2{\times}3)_1^2$ is $\PSp(2{\times}3) \wr 2 = \PSp(2{\times}3)^2 {:} 2$, extending to $\Sp(2{\times}3)^{\circ 2} {:} 2$ upon including $V_\odd$. 
Choose $S = \rJ_2$. Its double cover $2\rJ_2$ has two symplectic six-dimensional irreps, exchanged by the outer automorphism, and together they provide a map $2\rJ_2 \subset \Sp(2{\times}3)^2$ covering an inclusion $\rJ_2 \subset \Sp(2{\times}3)^{\circ 2} \subset \Sp(2{\times}3)^{\circ 2} {:} 2$.

The space of spin-$\frac32$ fields in $V$ is $\mathbf{14}_-^{\boxtimes 2}$, where $\mathbf{14}_- = \Alt^3(\mathbf{6}) \ominus \mathbf{6}$ is one of the fundamental $\Sp(2{\times}3)$-irreps, namely the one with nontrivial central character. 
The two choices of six-dimensional irrep of $2\rJ_2$ lead to the same 14-dimensional irrep $\mathbf{14}_-$, and a character table calculation confirms that $\mathbf{14}_-^{\boxtimes 2}$ has a unique fixed point upon restriction to $\rJ_2$.

The adjoint representation $\mathfrak{sp}(2{\times}3)^2$ is not simple upon restricting to $\rJ_2$, but it becomes simple if we include the outer automorphism and work with $\rJ_2{:}2$. Because of the double cover, there are two inclusions $\rJ_2{:}2 \subset \Sp(2{\times}3)^{\circ 2} {:} 2$, differing only by the sign of the action of the ``${:}2$.'' By uniqueness, the $\rJ_2$-fixed point in $\mathbf{14}_-^{\boxtimes 2}$ is fixed by one of the copies of $\rJ_2{:}2$, and antifixed by the other.

\subsection{Existence for $\SU(6)_2$}\label{exist.su6}
  
  The automorphism group of $V_\ev = \SU(6)_2$ is $\Aut(\mathfrak{su}(6)) = \PSU(6){:}2$. This extends to $\Aut_{N{=}0}(V) = \PSU(6){:}2 \times 2$. Write $\mathbf{6}$ for the standard representation of $\SU(6)$, and $\mathbf{20} = \Alt^3(\mathbf{6})$ for the third fundamental representation. The $175$-dimensional space of spin-$\frac32$ fields arises as
  $$\Sym^2(\mathbf{20}) = \mathbf{175} \oplus \mathbf{35},$$
  where $\mathbf{35} = \mathfrak{su}(6)$ is the adjoint  representation.
  
  The Mathieu group $\rM_{21}$ is not sporadic, being isomorphic to $\mathrm{PSL}_3(\bF_4) = \rL_3(4)$. It embeds into $\PSU(6)$ via the 6-dimensional representation of its 6-fold cover $6\rM_{21}$ (which GAP's character table library only knows under the name ``$6.\rL_3(4)$,'' but we will use the Mathieu name). The adjoint representation $\mathbf{35}$ remains simple when restricted to $\rM_{21}$, and so the conformal vector $\nu$ is the unique $\rM_{21}$-fixed spin-$2$ field. On the other hand, $\mathbf{20}$ breaks over $2\rM_{21}$ as $\mathbf{10} \oplus \overline{\mathbf{10}}$, where $\mathbf{10}$ and $\overline{\mathbf{10}}$ are a dual pair of $10$-dimensional complex irreps of $2\rM_{21}$, and so 
    $\Sym^2(\mathbf{20})$, and hence $\mathbf{175}$, has a unique fixed point when restricted to $\rM_{21}$, picking out the superconformal vector $\tau$.
    
    We remark that, although $\rM_{21} = \rL_3(4)$ is not sporadic, it is exceptional, having very large Schur multiplier ($3 \times 4^2$) and outer automorphism group ($D_{12} = 2 \times S_3$). The outer automorphism of order $3$ permutes the three double covers.
   In particular, there is no (interesting) group ``$2.\rM_{21}.3$,'' and so no ``$6.\rM_{21}.3$'' which could embed into $\SU(6)$. 
    (GAP can compute, in a few seconds, that $\Aut(6.\rM_{21}) = \rM_{21}{:}2^2$.)

  There are three groups of shape $\rM_{21}{:}2$, coming from the three conjugacy classes of order-$2$ element in $\operatorname{Out}(\rM_{21}) = D_{12}$. The one corresponding to the central element in $D_{12}$ is called $\rM_{12}{:}2\mathrm{a}$. It has a 6-fold central extension, and GAP is able to  build the group $6\rM_{21}{:}2{\mathrm{a}}$ from the $6$-dimensional matrix representation over $\bF_{25}$ listed in the ATLAS, convert it into a permutation group, and calculate its character table. It does have a $6$-dimensional complex irrep, and so embeds into $\SU(6)$ covering an embedding $\rM_{21}{:}2{\mathrm{a}} \subset \PSU(6)$. The other two groups, $\rM_{21}{:}2\mathrm{b}$ and $\rM_{21}{:}2\mathrm{c}$, can be built as normal subgroups of $\Aut(6\rM_{21}) = \rM_{21}{:}2^2$. One can see (for instance by calculating Schur multipliers using Holt's program ``Cohomolo'') that the outer automorphisms $2\mathrm{b}$ and $2\mathrm{c}$ exchange $\mathbf{6}$ and $\overline{\mathbf{6}}$. In particular, $\rM_{21}{:}2^2$ is not a subgroup of $\PSU(6)$, but is a subgroup of $\PSU(6){:}2$. There are, therefore, four embeddings $\rM_{21}{:}2^2 \subset \PSU(6){:}2 \times 2$. These differ by signs in their actions on the spin-$\frac32$ fields. It follows that one of these embeddings, but not the others, fixes $\tau$.
    
  \subsection{Existence for $\SU(6)_1^2$}\label{exist.su6.2} 
  
The automorphism group of $V_\ev = \SU(6)_1^2$ is $(\PSU(6){:}2)^2{:}2$, extending to the diagonal double cover $2.(\PSU(6){:}2)^2.2$ upon including $V_\odd$. The finite simple group called ``$\rU_4(3)$'' in the ATLAS has a large Schur multiplier~($4\times3^2$) and a large outer automorphism group~($D_8$). Of its 6-fold covers, the cover called ``$6_1.\rU_4(3)$'' in GAP's character table libraries has two dual six-dimensional irreps, exchanged by an outer automorphism. These provide a diagonal inclusion $\rU_4(3) \subset 2.\PSU(6)^2 \subset 2.(\PSU(6){:}2)^2.2$. (The diagonal double cover of $\PSU(6)^2$ restricts trivially along the diagonal $\rU_4(3) \subset \PSU(6)^2$.)
  
  Writing $\mathbf{20} = \Alt^3(\mathbf{6})$ for third fundamental irrep of $\SU(6)$, the space of spin-$\frac32$ fields in $V$ is $\mathbf{20}^{\boxtimes 2}$. The central character of $\mathbf{20}$ is just a sign, and so $\mathbf{20}$ pulls back from $2\PSU(6)$; its restriction to $2\rU_4(3)$ is simple (and in fact the unique  $2\rU_4(3)$-irrep of dimension $20$). A character table calculation confirms that $\mathbf{20}^{\otimes 2}$ has a unique $\rU_4(3)$-fixed point $\tau$.
  
  The adjoint representation $\mathfrak{su}(6)^{\boxplus 2} = \mathbf{35}^{\boxplus 2}$ of $2.(\PSU(6){:}2)^2.2$ is not simple when restricted to $\rU_4(3)$, but it becomes simple upon including the outer automorphism called ``$2_3$'' in GAP's character table libraries. There are two ways to extend the inclusion $\rU_4(3) \subset 2.(\PSU(6){:}2)^2.2$ to $\rU_4(2){:}2_3$, differing by a sign; one of them fixes $\tau$ (and the other antifixes $\tau$). Thus we may take $S = \rU_4(3){:}2_3$ to complete the proof: $\tau$ and $\nu$ are the unique $S$-fixed fields of spin $\frac32$ and $2$, respectively.
  
  Using the ``${:}2$s'' in $2.(\PSU(6){:}2)^2.2$, we may in fact embed the full automorphism group $\rU_4(3){:}D_8$, in various ways differing by some signs, exactly one of which fixes $\tau$.
       
  \subsection{Existence for $\Spin(12)_2$}\label{exist.spin12}
  
  $V_\ev = \Spin(12)_2$ has two anyons of conformal dimension $\frac32$, and the choice breaks
  $\Aut(V_\ev) = \Aut(\mathfrak{so}(12)) = \mathrm{PO}(12)$  to $\PSO(12)$, which then extends to $\Aut_{N{=}0}(V) = \PSO(12) \times 2$. 
  
  The group $2\rM_{12}$ has a unique $12$-dimensional irrep.
  It is real, and the central element acts by~$-1$. This provides an inclusion $2\rM_{12} \mono \SO(12)$ covering $\rM_{12} \mono \PSO(12) \subset \PSO(12) \times 2$. The adjoint representation $\mathbf{66}$ of $\PSO(12)$ remains simple upon restriction to $\rM_{12}$, and so the conformal vector $\nu$ is the unique $\rM_{12}$-fixed spin-$2$ field. 
  
  Write $\mathbf{32}_+$ and $\mathbf{32}_-$ for the two irreducible spinor representations of $\Spin(12)$; both are quaternionic, and $\mathbf{32}_+$ is in fact a representation of the double cover $\SO^\pm(12)$ of $\PSO(12)$. 
  Recall that double covers of $G$ are classified by $\H^2(BG; \bZ_2)$ --- for instance, 
  the three double covers $\SO(12)$, $\SO^+(12)$, and $\SO^-(12)$ correspond to the three nontrivial classes in $\H^2(B\PSO(12); \bZ_2) = (\bZ_2)^2$. Since $\H^2(B\rM_{12}; \bZ_2) = \bZ_2$ and the double cover $\SO(12) \to \PSO(12)$ restricts nontrivially to $\rM_{12}$, it follows that one of the double covers $\SO^+(12) \to \PSO(12)$ and $\SO^-(12)\to \PSO(12)$ --- the former, say --- restricts trivially over $\rM_{12}$ and the other restricts to the nontrivial double cover. Said another way, $\mathbf{32}_+$ restricts to a representation of $\rM_{12}$ whereas $\mathbf{32}_-$ restricts to a representation of $2\rM_{12}$ in which the centre acts nontrivially. In fact, $\mathbf{32}_-$ remains simple upon restriction to $2\rM_{12}$ --- it is the unique quaternionic irrep of $2\rM_{12}$ --- whereas $\mathbf{32}_+$ breaks over $\rM_{12}$ as $\mathbf{16} \oplus \overline{\mathbf{16}}$, where $\mathbf{16}$ and $\overline{\mathbf{16}}$ are the two dual complex irreps of $\rM_{12}$.
  
  The two $462$-dimensional representations of $\PSO(12)$ appear as
  $$ \Sym^2(\mathbf{32}_+) = \mathbf{462}_+ \oplus \mathbf{66}, \qquad \Sym^2(\mathbf{32}_-) = \mathbf{462}_- \oplus \mathbf{66}.$$
  Let us choose our SVOA $V$ to be the simple current extension of $V_\ev = \Spin(12)_2$ for which the spin-$\frac32$ fields in $V_\odd$ form the representation $\mathbf{462}_+$. Then, upon restriction to $\rM_{12}$, we have
  $$ \mathbf{462}_+ |_{\rM_{12}} = \Sym^2(\mathbf{16} \oplus \overline{\mathbf{16}}) \ominus \mathbf{66}, $$
  and this has a unique $\rM_{12}$-fixed point (since $\mathbf{16}$ and $\overline{\mathbf{16}}$ are dual irreps and $\mathbf{66}$ is irreducible). 
  This provides an $\rM_{12}$-fixed spin-$\frac32$ field $\tau$ in $V$.
  
  For comparison, $\mathbf{462}_-$ does not have any $\rM_{12}$-fixed points.
  The reader may worry that we seem to have made no choices --- we started with the unique $12$-dimensional irrep of $2\rM_{12}$ --- and somehow broke the symmetry between $\mathbf{462}_+$ and $\mathbf{462}_-$. In fact, the $12$-dimensional irrep $2\rM_{12}$ is unique in the sense that there is a unique conjugacy class of irreducible maps $2\rM_{12} \to \rO(12)$, but there are two conjugacy classes of irreducible maps $2\rM_{12} \to \SO(12)$, exchanged by the outer automorphism of the target. A choice of one of these is what breaks the symmetry between $\mathbf{462}_+$ and $\mathbf{462}_-$.
  
  We remark that $\PSO(12) \times 2$ contains two copies of $\rM_{12}{:}2$, differing by a sign in how they act  on the spin-$\frac32$ fields, just like in \S\ref{exist.sp3}. It follows that the superconformal vector is fixed by one $\rM_{12}{:}2$ and antifixed by the other.

\subsection{Existence for $\Spin(12)_1^2$}\label{exist.spin12.2} 

The automorphism group of $V_\ev = \Spin(12)_1^2$ is $\rP\rO(12)\wr 2 = \rP\rO(12)^2{:}2$. The space of spin-$\frac32$ fields in $V_\odd$ is $\mathbf{32}_+^{\boxtimes 2}$. This choice breaks and lifts each $\rP\rO(12)$ to a copy of $\SO^+(12)$, so that $\Aut_{N{=}0}(V) = \SO^+(12)^{\circ 2}{:}2$. (As above, the $\circ$ denotes a central product.) As in \S\ref{exist.spin12}, $2\rM_{12}$ has a unique isomorphism class of $12$-dimensional irrep $2\rM_{12} \subset \rO(12)$, which splits into two conjugacy classes of embeddings $2\rM_{12} \subset \SO(12)$. In \S\ref{exist.spin12} we used the one for which the corresponding map $2\rM_{12} \to \SO^+(12)$ factored through $\rM_{12}$; we now use the other one, giving an injection $2\rM_{12} \subset \SO^+(12)$. Its diagonal $2\rM_{12} \subset \SO^+(12)^2$ covers an injection $\rM_{12} \subset \SO^+(12)^{\circ 2} \subset \SO^+(12)^{\circ 2}{:}2$. The representation $\mathbf{32}^+$ of $\SO^+(12)$ remains simple (and quaternionic) upon restriction to this $2\rM_{12}$. It follows that $\mathbf{32}_+^{\boxtimes 2}$ has a unique $\rM_{12}$-fixed point $\tau$. We remark that the $12$-dimensional irrep of $2\rM_{12}$ extends, in two ways differing by a sign, to $2\rM_{12}.2$, and so the inclusion $\rM_{12} \subset \SO^+(12)^{\circ 2}{:}2$ extends, again in two ways differing by a sign, to an inclusion of $\rM_{12}{:}2$. One of these embeddings fixes $\tau$, and the other antifixes it.

The adjoint representation $\mathbf{66}^{\boxplus 2} = \mathfrak{so}(12)^{\boxplus 2}$ of $\SO^+(12)^{\circ 2}{:}2$ is not simple when restricted to $\rM_{12}$, but becomes simple upon inclusion of the outer automorphism. All together, we see that taking $S = \rM_{12}{:}2 \subset \SO^+(12)^{\circ 2}{:}2$ provides unique spin-$\frac32$ and spin-$2$ fields $\tau$ and $\nu$.

In \S\ref{unique.spin1212} we will give a different description of the superconformal structure, and see that $\Aut_{N{=}1}(V)$ is not just $S = \rM_{12}{:}2$ but in fact $2^{10}{:}\rM_{12}{:}2$.

\subsection{Existence for $\Sp(2{\times}6)_1$}\label{sp6.exist}
  
  The automorphism group $\PSp(2{\times}6)$ of $V_\ev = \Sp(2{\times}6)_1$ extends to $\Aut(V) = \PSp(2{\times}6) \times 2$ upon including $V_\odd$ because 
  $\mathbf{429} = \Alt^6(\mathbf{12}) \ominus \Alt^4(\mathbf{12})$ descends to $\PSp(2{\times}6)$, and so $V_\odd$ is acted on trivially by the centre of $\Sp(2{\times}6)$. 
  We choose $S = \rG_2(4)$, embedded into $\PSp(2{\times}6)$ via the 12-dimensional quaternionic irrep $2\rG_2(4) \subset \Sp(2{\times}6)$. The adjoint representation $\mathbf{78} = \Sym^2(\mathbf{12})$ remains simple upon restriction to $S$. On the other hand, a character table computation quickly confirms that $\Alt^6(\mathbf{12}) \ominus \Alt^4(\mathbf{12})$ has a unique $S$-fixed point.  Thus we have unique $\rG_2(4)$-fixed fields $\nu$ and $\tau$ as desired.
  
  We remark that $\PSp(2{\times}6) \times 2$ contains two copies of $\rG_2(4){:}2$, differing by a sign by their action on the spin-$\frac32$ fields. (One copy is inside $\PSp(2{\times}6)$ and the other is not. The two copies are perfect analogues of the two copies of
  $\rU_3(3){:}2$ inside $\PSp(2{\times}3) \times 2$ described in
  \S\ref{exist.sp3}.) It follows that one of these copies but not the other fixes $\tau$.
      
  \subsection{Existence for $\SU(12)_1$}\label{su12.existence}
  
  The automorphism group of $V_\ev = \SU(12)_1$ is $\PSU(12){:}2$, extending to a nontrivial double cover $\Aut_{N{=}0}(V) = 2\PSU(12){:}2$ upon including $V_\odd$. Suzuki's sporadic group $S = \Suz$ embeds into $2\PSU(12) = \SU(12)/\bZ_6$, since its Schur multiplier $6\Suz$ has a $12$-dimensional irrep. The adjoint representation $\mathfrak{su}(12) = \mathbf{143} = (\mathbf{12} \otimes \overline{\mathbf{12}}) \ominus \mathbf{1}$ remains simple upon restriction to $\Suz$, but a character table computation quickly confirms that $\Alt^6(\mathbf{12})$ has a unique $\Suz$-fixed point. Thus we have unique $\Suz$-fixed fields $\nu$ and $\tau$ as desired.
  
  We remark that $2\PSU(12){:}2$ contains two copies of $\Suz{:}2$, differing by a sign in how they act on the spin-$\frac32$ fields. It follows that one of these copies fixes $\tau$ and the other antifixes $\tau$.
  
  \subsection{Existence for $\Spin(16)_1 \times \Spin(8)_1$}\label{spin16spin8.existence}
  
  Let $V'$ denote the $N{=}1$ SVOA with even subalgebra $V'_{\ev} = \Spin(8)_1^3$ from \S\ref{exist.spinm}. The inclusion $V' \subset V$ constructed in \S\ref{subsec.inclusions}  is conformal, and so the superconformal vector for $V'$ is also a superconformal.
  
  The adjoint representation $\mathfrak{so}(16) \boxplus \mathfrak{so}(8)$ is not simple as an $\Aut_{N=0}(V)$-module, and so any group $S \subset \Aut_{N=0}(V)$ will fix at least two fields of spin-$2$. Thus the methods of this section do not quite apply. But a construction in the flavour of this section will be given in \S\ref{unique.spin16spin8}.
    
  \subsection{Existence for $\Spin(24)_1$}\label{spin24.existence}
  
  Let $V'$ denote the $N{=}1$ SVOA with even subalgebra $V'_{\ev} = \Spin(16)_1 \times \Spin(8)_1$ from \S\ref{spin16spin8.existence}. The inclusion $V' \subset V$ constructed in \S\ref{subsec.inclusions}  is conformal, and so the superconformal vector for $V'$ is also a superconformal.
  
  We may also show existence directly, using the methods of this section. $\Aut(\Spin(24)_1) = \mathrm{PO}(24)$, breaking and lifting to $\SO^+(24)$ upon extending from $V_\ev$ to $V = V_\ev \oplus V_\odd$, where $\SO^+(24)$ denotes the image of $\Spin(24)$ in the positive spinor representation $\mathbf{2048}_+$. Take $S = \Co_1$. There are two conjugacy classes of embeddings $S \subset \PSO(24)$, merging inside $\mathrm{PO}(24)$, and we will choose the one such that the double cover $\SO^+(24) \to \PSO(24)$ splits when restricted to $S$. (For the other embedding, $\SO^+(24) \to \PSO(24)$ restricts to the nontrivial double cover of $S$ but $\SO^-(24) \to \PSO(24)$ splits. We saw the same phenomenon happen for $\rM_{12} \subset \PSO(12)$.) The adjoint representation $\mathfrak{so}(24) = \mathbf{276}$ remains simple when restricted to $S$, but $\mathbf{2048}_+$ splits at $\mathbf{1} \oplus \mathbf{276} \oplus \mathbf{1771}$. Thus we have unique $\Co_1$-fixed fields $\nu$ and $\tau$ as desired.
   
   The $N{=}1$ extension of $V_{\ev} = \Spin(24)_1$ is due to \cite{JohnDuncan-thesis,MR2352133}, where it is called $V^{f\natural}$ and shown to be unique and to have automorphism group $\Co_1$. The construction in \cite{MR2352133} presents $V^{f\natural}$ directly as an extension of $\Spin(24)_1$. The construction in \cite{JohnDuncan-thesis} uses $V'$, which is built as a canonical $N{=}1$ subalgebra of the ``supersymmetric $\rE_8$ lattice'' SVOA. This explains the appearance of $\rE_8$ in \S\ref{unique.spin16spin8}.

\section{$\Spin(m)$ cases} \label{sec.spinm3}

In this section we analyze in detail the $N{=}1$ extensions of $\Spin(m)_3$ and $\Spin(m)_1^3$ constructed in~\S\ref{exist.spinm} and~\S\ref{exist.spinm.3}. Henceforth we will use the name ``$\SO(m)_3$'' for the simple current extension of $\Spin(m)_3$ and ``$\Spin(m)^3_1/\bZ_2$'' for the simple current extension of $\Spin(m)_1^3$; see~\S\ref{subsec:small-m} for discussion of these naming choices.
We first show (\S\ref{sec.spinm3auto}) that symmetric group $S_{m+1}$ from \S\ref{exist.spinm} exhausts all of $\Aut_{N{=}1}(\SO(m)_3)$. Then we show (\S\ref{subsec.spinm3unique}) that, up to the action of $\Aut_{N{=}0}(\SO(m)_3) = \rO(m)$, there is a unique conformal vector $\tau$. Finally (\S\ref{unique.spinm1.3}), we show that the $N{=}1$ structure on $\Spin(m)^3_1/\bZ_2$ is unique, and calculate its automorphism group.

\subsection{Automorphism group for $\Spin(m)_3$}  \label{sec.spinm3auto}
\begin{proposition}
  $\Aut_{N{=}1}(\SO(m)_3) = S_{m+1}$ for the $N{=}1$ structure constructed in \S\ref{exist.spinm}.
\end{proposition}

\begin{proof}
  Recall that the space of spin-$\frac32$ fields in $\SO(m)_3$ is the $\rO(m)$-irrep $\Sym^3(\mathbf{m}) \ominus \mathbf{m}$, where $\mathbf{m}$ is the vector representation of $\rO(m)$; the subspace of self-adjoint fields is formed by interpreting $\mathbf{m}$ as the real vector representation  $\bR^m$ of the compact form of $\rO(m)$.
  Note that $\Sym^3(\mathbf{m}) \ominus \mathbf{m} \subset \Sym^3(\mathbf{m})$, and identify $\Sym^3(\mathbf{m})$ with the space of cubic functions on $\bR^m$.
   As in \S\ref{exist.spinm}, give $\bR^{m+1}$ its standard coordinates, and embed $\bR^m \subset \bR^{m+1}$ as the set of vectors $x = (x_0,\dots,x_{m+1})$ such that $\sum_i x_i = 0$. In this parameterization, the superconformal vector $\tau$ corresponded to the function
   $$ \tau(x) = \sum_{i=0}^m x_i^3.$$
   We saw already that this function is in fact in the submodule $\Sym^3(\mathbf{m}) \ominus \mathbf{m} \subset \Sym^3(\mathbf{m})$, since $\mathbf{m}$ has no $S_{m+1}$-fixed points, whereas $\tau$ is manifestly $S_{m+1}$-fixed. To prove the proposition, we must show simply that there are no further symmetries stabilizing $\tau$.
   
   We will strongly maximize $\tau(x)$ subject to the constraints
   $$ c_1(x) = \sum x_i = 0, \qquad c_2(x) = \sum x_i^2 = m(m+1).$$
   By ``strongly maximize,'' we mean to look for points which are maxima with negative-definite Hessian (as opposed to negative-semidefinite).
   The first constraint just forces $x$ into our $\bR^m \subset \bR^{m+1}$, and the second constraint cuts down to an $(m-1)$-dimensional sphere $S^{m-1} \subset \bR^m$.
   We claim that there are precisely $m+1$ maxima, namely the vectors 
   $$ \bigl(-1,-1,\dots,-1,m,-1,\dots,-1\bigr).$$
   This would complete the proof, since any automorphism of $\tau$ must permute these $m+1$ maxima, and the only such elements of $\rO(m)$ are the elements of $S_{m+1}$.
   
   When maximizing a function subject to constraints, one does not have that the derivative of the function vanishes, but rather that it is in the span of the derivatives of the constraints:
   $$\d\tau(x) = \bigl(3 x_0^2,\dots, 3x_m^2\bigr) = a \,\d c_1 + b \,\d c_2 = a\bigl(1,\dots,1\bigr) + b\bigl(2x_0,\dots,2x_m\bigr),$$
   i.e.\ there are real numbers $a,b \in \bR$ such that $3x_i^2 = a + 2bx_i$ for all $0 \leq i\leq m$. Summing these equations gives:
   $$ 3c_2(x) = \sum 3x_i^2 = (m+1)a + 2b \sum x_i = (m+1)a + 2b c_1(x).$$
   After imposing the constraints $c_2(x) = m(m+1)$ and $c_1(x) = 0$, we find that $a = 3m$.
   
   If we instead multiply the equation $3x_i^2 = a + 2bx_i$ by $x_i$ and then sum,  we learn that $3\tau(x) = ac_1(x) + 2b c_2(x) = 2b m(m+1)$. Since $\tau$ is an odd function under $x \mapsto -x$, its maxima are  positive, and so $b$ is positive.
   
   Substituting $a = 2m$, we have:
   $$ 3x_i^2 = 3m + 2bx_i, \qquad \text{i.e.\ } x_i = \frac b 3 \pm \sqrt{ \left( \frac b 3 \right)^2 + m }.$$
   In particular, no $x_i$ is zero at any critical point.
   
   The matrix of second derivatives is
   $$ H_{ij} = \frac{\partial^2 \tau}{\partial x_i \partial x_j} = 6x_i \delta_{ij}.$$
   This matrix is nondegenerate since no $x_i$ is zero. If $x$ is to strongly-maximize $\tau$ subject to the constraints, then this matrix must restrict to a negative-definite matrix on the tangent space $$\rT_x S^{m-1} = \left\{v \in \bR^{m+1} \st \sum v_i = \sum 2 x_i v_i = 0\right\}.$$
   Consider the (nontangent) vectors $(1,1,\dots,1)$ and $(\frac1{x_0},\frac1{x_1},\dots,\frac1{x_m})$. The first of these is null for the inner product determined by $H$: its $H$-norm is $\sum 6x_i(1)(1) = 6c_1(x) = 0$. Furthermore, the $H$-dot-product of these two vectors is $\sum 6x_i(1)(\frac1{x_i}) = 6(m+1)$, and so they are linearly independent. Finally, both vectors are $H$-orthogonal to the tangent space $\rT_x S^{m-1}$, since $H(1,\dots,1) = (6x_0,\dots,6x_m) = 3\,\d c_2$ and $H(\frac1{x_0},\dots,\frac1{x_m}) = (6,\dots,6) = 6\,\d c_1$. Thus, in terms of the $H$-norm, we have split $\bR^{m+1}$ as an orthgonal direct sum of an $(m-1)$-dimensional negative-definite space $\rT_x S^{m-1}$ and a two-dimensional space which contains a null vector. Nondegeneracy of $H$ then shows that that $2$-dimensional space must have signature $(1,1)$. All together, we find that $H$ has signature $(1,m)$, which is to say exactly one of the $x_i$ is positive and the other $m$ are negative.
   
   But $x_i = \frac b 3 \pm \sqrt{ ( \frac b 3 )^2 + m }$ with $b$ positive, and so $x_i$ is positive or negative depending on the choice of $\pm$. Thus we find: exactly one of the $\pm$s is $+$ and the other $m$ $\pm$s are $-$. 
   Summing the $x_i$s then gives
   $$ 0 = c_1(x) = (m+1) \frac b 3 - (m-1) \sqrt{ \left( \frac b 3 \right)^2 + m },$$
   which solves to 
   $$ \frac b3  = \frac{m-1}2,$$
   and so
   $$ x_i = \frac{m-1}2 \pm \sqrt{ \left( \frac{m-1}2\right)^2 + m} = \frac{m-1}2 \pm \sqrt{ \left(\frac{m+1}2\right)^2} = m \text{ or} -1.$$
   Thus the (strong) maxima are precisely the ones claimed.
\end{proof}

\subsection{Uniqueness for $\Spin(m)_3$} \label{subsec.spinm3unique}

In \S\ref{subsec.inclusions} we observed that the inclusion $\Spin(m-1)_3 \subset \Spin(m)_3$ of simply connected WZW algebras extends to an inclusion $\SO(m-1)_3 \subset \SO(m)_3$. This inclusion is compatible with the superconformal structures constructed in \S\ref{exist.spinm} as follows. Continue to identify the $m$-dimensional vector representation $\mathbf{m}$ of $\SO(m)$ with the subspace of $\bR^{m+1}$ cut out by the equation $\sum_{i=0}^m x_i = 0$, where $x_0,\dots,x_m$ are the coordinates on $\bR^{m+1}$. Then the vector representation of $\SO(m-1)$ is just the subspace for which $x_m = 0$. Write $\tau = \sum_{i=0}^m x_i^3$ for the superconformal vector for $\SO(m)_3$, and $\tau' = \sum_{i=0}^{m-1} x_i^3$ for the superconformal vector for $\SO(m-1)_3$. Then the difference $\tau'' = \tau - \tau' = x_m^3$ lives in the coset algebra $\SO(m)_3 / \SO(m-1)_3$, and is in fact a superconformal vector therein. More generally, we will say that an inclusion $(V',\tau') \subset (V,\tau)$ of $N{=}1$ SVOAs is \define{supersymmetric} if $\tau'' = \tau-\tau' \in V/V'$, and \define{superconformal} if it is both supersymmetric and conformal (i.e.\ if $\tau'' = 0$).

Our strategy to prove the uniqueness of the $N{=}1$ structure on $\SO(m)_3$ will be to show that, for any superconformal vector $\tau$, there is a supersymmetric inclusion from $\SO(m-1)_3$, whose superconformal vector $\tau'$ is unique by induction.

\begin{proposition}
  The superconformal vectors for $\SO(m)_3$ form a single orbit under the action by $\Aut_{N{=}0}(\SO(m)_3) = \rO(m)$.
\end{proposition}

\begin{proof}
  A superconformal vector $\tau \in \SO(m)_3$ is in particular a traceless symmetric three-tensor, and so defines, as in \S\ref{sec.spinm3auto}, a  cubic function, which we will also call $\tau$, on $\mathbf{m} = \bR^m$. 
  Let $S^{m-1} \subset \bR^m$ denote the unit sphere, and choose $e_0 \in S$ to maximize $\tau|_{S^{m-1}}$. Complete $e_0$ to an orthonormal basis $e_0,e_1,\dots,e_{m-1}$ of $\bR^m$. We will show that the inclusion $\SO(m-1)_3 \subset \SO(m)_3$ along the subspace spanned by $e_1,\dots,e_{m-1}$ is  supersymmetric, for some superconformal vector $\tau'$ on $\SO(m-1)_3$.
  
  We may coordinatize a neighbourhood of $e_0 \in S^{m-1}$ by vectors $\vec y = (y_1,\dots,y_{m-1})$ by projection to $S^{m-1}$: 
  $$ \frac{e_0 + \sum_{i>0} y_i e_i}{(1 + \| y\|^2)^{1/2}} \in S^{m-1}.$$
  Let us Taylor-expand $\tau|_{S^{m-1}}$ near $e_0$. Since $\tau$ is homogeneous cubic, we have
  \begin{multline*}
    \tau\left(\frac{e_0 +  y_i e_i}{(1 + \| y\|^2)^{1/2}} \right) = \frac{\tau\left(e_0 +  y_i e_i\right)}{(1 + \delta_{ij}y_iy_j)^{3/2}} \\ = \left( \tau^{(0)} + \tau^{(1)}_i y_i + \frac12 \tau^{(2)}_{ij} y_i y_j + \frac16 \tau^{(3)}_{ijk} y_i y_j y_k \right) 
    \left( 1 - \frac32 \delta_{ij} y_i y_j + O(y^4) \right)
    \\ = \tau^{(0)} + \tau^{(1)}_i y_i + \left( \frac12 \tau^{(2)}_{ij} - \tau^{(0)}\frac32 \delta_{ij} \right) y_i y_j + \frac16 \tau^{(3)}_{ijk} y_i y_j y_k + O(y^4). \end{multline*}
  To save space we have dropped the summation signs: a repeated index in a monomial is to be summed over, with all indices ranging over $\{1,\dots,m-1\}$. By ``$O(y^4)$'' we mean of course terms vanishing to fourth order and higher in $\vec y$.
  
  By assumption, $e_0$ is a maximum of $\tau|_{S^{m-1}}$. This implies:
  $$ \tau^{(0)}>0, \qquad \tau^{(1)} = 0, \qquad\text{and}\qquad \frac12\left(\tau^{(2)} - 3 \tau^{(0)}\right) \leq 0.$$
  By ``$\tau^{(2)} - 3 \tau^{(0)} \leq 0$'' we mean that the symmetric matrix $\tau^{(2)}_{ij} - 3 \tau^{(0)}\delta_{ij}$ is negative semidefinite.
  
  We furthermore know that $\tau$ is a superconformal vector. Thinking of $\tau$ as a symmetric 3-tensor on all of $\bR^m$, the assertion that $\tau$ is a superconformal vector is equivalent to the assertion that a certain ``partial trace'' of $\tau \otimes \tau$ returns some fixed multiple of the identity map on $\Alt^2(\bR^m)$; after rescaling $\tau$, we may assume that multiple to be $\pm1$, and checking the examples from \S\ref{exist.spinm} shows it to be negative. So as not to confuse with the indices $i,j,k \in \{1,\dots,m-1\}$, let us use indices $p,q,r,s,t$ to range over $\{0,\dots,m-1\}$. For instance, the symmetric 3-tensor $\tau_{prq}$ is related to the function $\tau(-)$ by $\tau(y_0 e_0 + y_i e_i) = \frac1{3!} \tau_{prq} y_p y_q y_r$. In coordinates, the equation to be a superconformal vector is then:
  $$ \tau_{rps} \tau_{rqt} - \tau_{rpt} \tau_{rqs} = -(\delta_{ps}\delta_{qt} - \delta_{pt}\delta_{qs}).$$
  
  Consider this equation when $p=s=0$, $q=i$, and $t=j$, where $i,j \neq 0$. In terms of the previous Taylor expansion, we find:
  $$ 6\tau^{(0)} \tau_{ij}^{(2)} - \tau_{ki}^{(2)} \tau_{kj}^{(2)} = -\delta_{ij}$$
  The factor of $6 = 3!$ comes from the combinatorics of Taylor expansion: $\tau(y_0e_0) = \frac1{3!} \tau_{000} y_0^3 = \tau^{(0)}y_0^3$. If we were working in some other basis, then there would also be $\tau^{(1)}\tau^{(1)}$ and $\tau^{(1)}\tau^{(3)}$ terms, but we have assumed $\tau^{(1)}=0$.
  Completing the square gives
  $$ 1 + \bigl(3\tau^{(0)}\bigr)^2 = \bigl( \tau^{(2)} - 3\tau^{(0)} \bigr)^2,$$
  where $1$ and $\tau^{(0)}$ stand for the corresponding scalar matrices, and the right-hand side means matrix multiplication. Thus, for fixed $\tau^{(0)}$, the eigenvalues of $\tau^{(2)}$ are $3\tau^{(0)} \pm \sqrt{1 + (3\tau^{(0)})^2}$.
  
  But $\tau^{(2)} - 3\tau^{(0)}$ is negative semidefinite! It follows that $\tau^{(2)}$ is a scalar matrix:
  $$ \tau^{(2)} = 3\tau^{(0)} - \sqrt{1 + (3\tau^{(0)})^2}.$$
  In particular, $\tau''(y_0,\vec y) = \tau^{(0)}y_0^3 + \frac12 y_0 \tau^{(2)}_{ij} y_i y_j$ is $\SO(m-1)$-invariant, and so an element of the coset $\SO(m)_3/\SO(m-1)_3$. On the other hand, $\tau' = \tau^{(3)}$ is traceless since $\tau$ was, hence in $\SO(m-1)_3$, and it is not hard to check that it is a superconformal vector. All together, we have written
  $$ \tau = \tau' + \tau'', \qquad \tau' \in \SO(m-1)_3, \quad \tau'' \in \SO(m)_3/\SO(m-1)_3.$$
  
  This completes the proof of uniqueness of the superconformal vector $\tau$: $\tau'$ is unique by induction,  the space of $\SO(m-1)$-fixed spin-$\frac32$ fields is one-dimensional (since $\Sym^3(\mathbf{m}) \ominus \mathbf{m}$ splits over $\SO(m-1)$ as $\bigl(\Sym^3(\mathbf{m-1}) \ominus (\mathbf{m-1})\bigr) \oplus \bigl(\Sym^2(\mathbf{m-1})\ominus\mathbf 1\bigr) \oplus (\mathbf{m-1}) \oplus \mathbf 1$), and $\tau''$ is necessarily a superconformal vector, setting its normalization. But we may also confirm the uniqueness directly. We have not used that $\tau$ is traceless (except to confirm that $\tau'$ was traceless). The $0$-component of its trace is
  $$ 0 = \tau_{0pp} = \tau_{000} + \tau_{0ii} = 6\tau^{(0)} + \tau^{(2)}_{ii} = 6\tau^{(0)} + (m-1)\left(3\tau^{(0)} - \sqrt{1 + (3\tau^{(0)})^2}\right).$$
  This unpacks to a quadratic equation for $\tau^{(0)}$, only one of whose solutions is positive. 
\end{proof}

\subsection{Uniqueness and automorphisms for $\Spin(m)_1^3$} \label{unique.spinm1.3}

In \S\ref{exist.spinm.3} we wrote down the superconformal vector $\tau = \sum e_{i,i,i} \in \mathbf{m}^{\boxtimes 3}$ for the $\bZ_2$-extension $V = \Spin(m)_1^3/\bZ_2$ of $\Spin(m)_1^3$. 

\begin{proposition}
  The $N{=}1$ superconformal vector $\tau$ on $V$ is unique up to the action of $\Aut_{N{=}0}(V)$.
\end{proposition}

\begin{proof}
  Consider the conformal embedding $\Spin(m)_1^3 \subset \Spin(3m)_1$. We claim that it extends to an embedding of SVOAs $V \subset \SO(3m)_1$. Indeed, let us write $M$ for the abelian $\Spin(m)_1$-anyon of conformal dimension $1/2$. Its minimal-spin fields form the vector representation $\mathbf{m}$ of $\Spin(m)$, and the $\bZ_2$-extension $\SO(m)_1 = \Spin(m)_1 \oplus M$ is precisely the free fermion algebra $\Fer(m)$. Free fermion algebras tensor well: $\SO(m)_1^3 = \Fer(m)^3 = \Fer(3m) = \SO(3m)_1$. All together $\SO(m)_1^3$ is a $\bZ_2^3$-extension of $\Spin(m)_1^3$, and the $\bZ_2$-extension $V$ that we care about is a subextension. In terms of fields, the spin-$\frac32$ fields in $V$ form the vector space $\mathbf{m}^{\boxtimes 3}$, whereas the spin-$\frac32$ fields in $\SO(3m)_1$ form $\Alt^3(\mathbf{m}^{\boxplus 3})$, and the embedding $V \subset \SO(3m)_1$ is the standard inclusion $\mathbf{m}^{\boxtimes 3} \subset \Alt^3(\mathbf{m}^{\boxplus 3})$.
  
  Thus any superconformal vector $\tau$ on $V$ provides a superconformal vector on $\SO(3m)_1$. The superconformal vectors on a free fermion algebra like $\SO(3m)_1 = \Fer(3m)$ have a beautiful classification \cite{MR791865}: they are in bijection with Lie algebra structures on $\mathbf{3m} = \bR^{3m}$ whose Killing form is (minus) the Euclidean inner product; up to isomorphism, they are in bijection with $3m$-dimensional semisimple Lie algebras; the bijection simply interprets the superconformal vector $\tau \in \Alt^3(\mathbf{3m})$ as the Lie bracket. In the case at hand, we start with $\tau \in \mathbf{m}^{\boxtimes 3} \subset \Alt^3(\mathbf{m}^{\boxtimes 3})$, which is to say that the only nontrivial brackets are of the form $$\bigl[(\text{element from $\mathbf{m}'$}), (\text{element from $\mathbf{m}''$})\bigr] = (\text{element from $\mathbf{m}'''$})$$
  where $\mathbf{m}'$, $\mathbf{m}''$, and $\mathbf{m}'''$ are the three ``coordinate'' $\mathbf{m}$s in some order.
  
  In particular, each ``coordinate'' $\mathbf{m} \subset \mathbf{3m}$ is an abelian subalgebra of our to-be-determined semisimple Lie algebra. Over the complex numbers, it is not true that every abelian subalgebra extends to a Cartan subalgebra: the $N^2$-dimensional space of block matrices of the form $\bigl( \begin{smallmatrix} 0 & B \\ 0 & 0 \end{smallmatrix}\bigr)$ inside $\mathfrak{sl}(2N)$ is the standard counterexample. But it is true for compact groups, which is to say for Lie algebras over $\bR$ with a negative-definite Killing form. It follows that our $3m$-dimensional semisimple Lie algebra has rank at least~$m$. The only possibility is $\mathfrak{su}(2)^m$.
  
  The ``coordinate'' $\mathfrak{su}(2)$s inside $\mathfrak{su}(2)^m$ are canonically determined: they are the $\mathfrak{su}(2)$ subalgebras with maximal centralizer $\mathfrak{su}(2)^{m-1}$. Each of these coordinate $\mathfrak{su}(2)$s intersects each coordinate $\mathbf{m}$ in exactly one dimension (which can be resolved further to two vectors by using the metric). Thus the superconformal vector ends up equipping each $\mathbf{m}$ with a basis-up-to-sign $\{e_1,\dots,e_m\}$, ordered up to simultaneously reordering the three bases, and in this basis 
  $  \tau = \sum e_{i,i,i} = \sum e_{i} \boxtimes e_i \boxtimes e_i.$ 
\end{proof}

\begin{proposition}
When $m\neq 4$, the stabilizer of $\tau$ is $\Aut_{N{=}1}(V) = 2^{2(m-1)} {:} (S_3 \times S_m)$, where $S_3$ and $S_m$ act by the tensor product of the standard $2$- and $(m-1)$-dimensional representations (the quotient of the permutation representation by the fixed vector).
\end{proposition}

\begin{proof} 
When $m \neq 4,8$, the bosonic subalgebra $V_\ev = \Spin(m)_1^3$ has automorphism group $\Aut(V_\ev) = \mathrm{PO}(m)^3 {:} S_3$, lifting to $\Aut_{N{=}0}(V) = 2.\mathrm{PO}(m)^3 {:} S_3 = \rO(m)^{\circ 3} {:} S_3$, where $\mathrm{PO}(m)$ is the projective orthogonal group (isomorphic to $\SO(m)$ when $m$ is odd) and $\rO(m)^{\circ 3}$ denotes the central product $\rO(m)^3 / 2^2$.  When $m = 8$, $\Aut(V_\ev)$ is larger than $\mathrm{PO}(m)^3 {:} S_3$ due to the triality automorphism of $\Spin(8)$, but the extension to $V$ breaks these extra automorphisms and we still have $\Aut_{N{=}0}(V) = \rO(m)^{\circ 3} {:} S_3$. The extension $V \subset \SO(3m)_1$ does not break any of these symmetries,
but lifts the central product $\rO(m)^{\circ 3}$ to the product $\rO(m)^3$.
 
  What are the symmetries inside $\rO(m)^3{:}S_3 \subset \rO(3m)$ that preserve $\tau \in \SO(3m)_1$? The basis $\{e_i\}$ suffers a sign ambiguity: for each $i \in \{1,\dots,m\}$, there is a Klein-4 group $2^2$ acting by switching the signs of any two out of the three $e_i$s, and so $2^{2m} \subset \Aut_{N{=}1}(V)$. Any other symmetry must rearrange the basis. We may arbitrarily permute the set $\{1,\dots,m\}$, provided we apply the same permutation to all three bases $\{e_i\}$, providing a group $S_m \subset \Aut_{N{=}1}(V)$. Together with $2^{2m}$, we find that $2^{2m}{:}S_m$ is the full group of symmetries which do not permute the three $\mathbf{m}$s. 
  What about permutations $S_3$ permuting the three $\mathbf{m}$s? 
  The even permutations $A_3$ are harmless, but the odd ones are slightly subtle, at least in our description inside $\Fer(3m) = \SO(3m)_1$: due to the fermionic nature of the fields $e_i$, one should declare, for example, that the transposition $(12)$ acts on each each triple of $e_i$s by the signed permutation $\left( \begin{smallmatrix} & -1 & \\ -1 & & \\ & & -1\end{smallmatrix} \right)$ rather than the plain permutation $\left( \begin{smallmatrix} & 1 & \\ 1 & & \\ & & 1\end{smallmatrix} \right)$. These signs are irrelevant for the purposes of understanding $\mathbf{m}^{\boxtimes 3}$ as an abstract $S_3$-module: they arise simply from the embedding $\mathbf{m}^{\boxtimes 3} \subset \Alt^3(\mathbf{m}^{\boxplus 3})$.
  
  Thus we find that the stabilizer of $\tau$ under the $\rO(m)^3{:}S_3$ action is $2^{2m} {:} (S_3 \times S_m)$. Its image in $\rO(m)^{\circ 3}{:}S_3 = (\rO(m)^3{:}S_3) / 2^2$ is the group $2^{2(m-1)}{:}(S_3 \times S_m)$ described in the Proposition.
\end{proof}

\begin{proposition}
When $m = 4$, the automorphism group is $\Aut_{N{=}1}(V) = 2^6 {:} 3S_6$.
\end{proposition}  

The group $2^6{:}3S_6$ is the ``sextet group'' arising as a maximal subgroup of Mathieu's group $\rM_{24}$ and surveyed for example in Chapter~11 of \cite{MR1662447}. The extension $3S_6$ is perfect but not central, and it acts faithfully on $2^6$. This group has an index-two subgroup $2^6{:}3A_6$, and the extension $3A_6$ is central, perfect, and exceptional.
  
\begin{proof}
When $m = 4$, there is an exceptional isomorphism $\Spin(4)_1 = \SU(2)_1^2$. This does not affect $\Aut(\Spin(4)_1) = \mathrm{PO}(4) = \mathrm{PSU}(2) {:} 2$, but it does contribute extra automorphisms to $\Aut(\Spin(4)_1^3) = \Aut(\SU(2)_1^6)$ that permute the six copies of $\SU(2)_1$ independently. These permutations are compatible with the extension from $V_\ev = \Spin(4)_1^3$ to $V$, so that
  $$ \Aut_{N{=}0}(V) = \SU(2)^{\circ 6}{:} S_6.$$
  The further extension to $\SO(12)_1$ breaks these extra symmetries, and so $\Aut_{N{=}1}(V)$ may be larger than the group $2^6 {:} (S_4 \times S_3)$ predicted by the general case. 
  
  There is, manifestly, a map $\Aut_{N{=}1}(V) \to S_6$ given by looking just at the action on the Dynkin diagram of $\SU(2)_1^6$. 
  We will first calculate the kernel of this map. This kernel doesn't care that $S_6$ is larger than would be predicted from the general-$m$ case. For even $m>4$, in place of $S_6$ we would find $2^3 {:} S_3$ many Dynkin diagram automorphisms, and we would want to understand the map $2^{2(m-1)}{:}(S_m \times S_3) \to 2^3 {:} S_3$. The $S_3$s match, and so equivalently we want to understand the map $2^{2(m-1)}{:}S_m \to 2^3$. As in the proof of the previous Proposition, we may lift from $\rO(m)^{\circ 3} {:} S_3$ 
   to $2^{2m}{:}S_m \subset \rO(m)^3 {:} S_3 \subset \rO(3m)$. Each element of $2^{2m}{:}S_m$ consists of a triple of changes of basis of $\bR^m$, and the map to $2^3$ records whether each change of basis is oriented. An odd element of $S_m$ changes the orientation of $\bR^{3m}$, whereas $2^{2m}{:}A_m$ preserves the orientation of $\bR^{3m}$, and so the kernel is inside $2^{2m}{:}A_m$ and manifestly contains $A_m$. Finally, on the $2^{2m}$ part, the map to $2^3$ is the sum map $2^{2m} \to 2^2 \subset 2^3$, and so has kernel $2^{2(m-1)}$. All together, we find that $\ker \bigl( 2^{2m}{:} (S_m \times S_3) \to 2^3 {:} S_3\bigr) = 2^{2(m-1)}{:}A_m$, and since this automatically contains the kernel of $2^{2m}{:} (S_m \times S_3) \to 2^{2(m-1)}{:} (S_m \times S_3)$, we find
   $$ \ker \bigl( \Aut_{N{=}1}(V) \to \Aut(\text{Dynkin diagram})\bigr) = 2^{2(m-2)}{:} A_m,$$
   for any $m$.
   
   Specializing to $m=4$, we thus have a kernel of shape $2^4{:} A_4$, and the full automorphism group extends this by some subgroup of $S_6$. Unlike the larger alternating groups, $A_4$ is not simple: it has shape $2^2{:}3$. If we were working ``upstairs'' in $\rO(4)^3{:}S_3$ rather than $\rO(4)^{\circ 3}{:}S_3$, then we would have $2^6 {:} 2^2 {:} 3$, and the action of $2^2$ on $2^6$ would be nontrivial. But in the quotient $2^4 {:} 2^2 {:} 3$, the $2^2$ acts trivially on $2^4$, and so
   $$ \ker \bigl( \Aut_{N{=}1}(V) \to S_6\bigr) = 2^6{:}3, \text{ i.e. } \Aut_{N{=}1}(V) = (2^6{:}3).(\text{subgroup of }S_6)$$
   where $3$ acts by three copies of the two-dimensional simple representation on $2^2$. Of course, the subgroup of $S_6$ is precisely the image of $\Aut_{N{=}1}(V) \to S_6$. By the general-$m$ case, we know that this image contains $2 \times S_3$, where the $S_3$ permutes the six vertices of the Dynkin diagram in pairs, and the $2$ acts diagonally on the three pairs. In particular, the $2$ acts by an odd permutation, so that the image of $\Aut_{N{=}1}(V) \to S_6$ is not contained in $A_6$.
   
   To complete the proof (i.e.\ to identify this image, and to identify the extension), it suffices to write down an $N{=}1$ superconformal vector $\tau' \in V$ with manifest symmetry $2^6{:}3S_6$, since the first Proposition in this section shows that any two superconformal vectors are in the same orbit under $\rO(4)^{\circ3}{:}S_3 \subset \SU(2)^{\circ 6}{:}S_6$. In unpublished work \cite{HarveyMoore-unpub} described in \cite{GregLecture-Freedfest,GregLecture-Readfest}, Harvey and Moore have used the hexacode (c.f.\ Chapter~11 of \cite{MR1662447}) to construct an explicit $\tau'$ which is symmetric under (at least) $2^6{:}3A_6$. (Their motivation was to understand the full superconformal field theory studied in \cite{GTVW}, which has our $N{=}1$ SVOA $V = \SU(2)_1^6/\bZ_2$  as its chiral algebra.)
    In brief, one may find a subgroup $S = 2^6{:}3A_6 \subset \SU(2)^{\circ 6}{:}S_6$ which satisfies the conditions used in \S\ref{sec.existence}: it preserves a unique (up to scalar) $\tau' \in V_{3/2}$, but the adjoint representation $\mathfrak{su}(2)^6$ remains simple when restricted to $S$, and so $\tau'$ must be an $N{=}1$ superconformal vector.
   This shows that the image of $\Aut_{N{=}1}(V) \to S_6$ contains $A_6$. But we remarked already that the image is not contained in $S_6$, and so must be all of~$S_6$.
   
   One can also see that the symmetry group enhances from $2^6{:}3A_6$ to $2^6{:}3S_6$ as follows. There are two embeddings $2^6{:}3S_6 \subset \SU(2)^{\circ 6}{:}S_6$ containing the subgroup $2^6{:}3A_6$, which differ by the ``sign'' map from $2^6{:}3S_6$ to the centre of $\SU(2)^{\circ 6}{:}S_6$.
    Since $\tau'$ is the unique $2^6{:}3A_6$-vector in $V_{3/2}$, either extension to $2^6{:}3S_6$ must fix $\tau'$ up to a phase. It follows that one of the extensions fixes $\tau'$ and the other \define{antifixes} $\tau'$, meaning that it acts on $\tau'$ by the sign representation.
\end{proof}

\section{Leech lattice groups} \label{sec.uniqueness}

This section completes the proof of the Theorem by showing that for each of the WZW algebras $V_\ev = \Sp(2{\times}3)_2$, $\Sp(2{\times}3)_1^2$, $\SU(6)_2$, $\Sp(2{\times}6)_1$, $\SU(6)_1^2$, $\Spin(12)_2$, $\SU(12)_1$, $\Spin(12)_1^2$, and $\Spin(24)_1$, the superconformal vector $\tau \in V = V_\ev \oplus V_\odd$  found in Section~\ref{sec.existence} is unique up to the action of $\Aut_{N=0}(V)$, with the claimed stabilizer $\Aut_{N=1}(V)$.

The $V_\ev = \Spin(24)_1$ case is already known, and is one of the main results of \cite{MR2352133}. That paper uses the name $V^{f\natural}$ for the SVOA we have been calling $\SO^+(24)_1$, and we will use that name for the remainder of this Section.

\begin{theorem}[\cite{MR2352133}]
  Let $V^{f\natural} = \SO^+(24)_1$ denote the simple current extension of $V^{f\natural}_\ev = \Spin(24)_1$ through either anyon of conformal dimension $\frac32$. Then all superconformal vectors for $V^{f\natural}$ form a single orbit for $\Aut_{N{=}0}(V^{f\natural}) = \SO^+(24)$, with stabilizer $\Aut_{N{=}1}(V^{f\natural}) = \Co_1$. \qed
\end{theorem}

Our strategy will be to derive uniqueness in all other cases from the uniqueness of $V^{f\natural}$. We showed in \S\ref{subsec.inclusions} at the level of non-supersymmetric SVOAs that each of the algebras $V$ in question embeds into $V^{f\natural}$. Choose some SVOA $W$ such that the embedding $V \subset V^{f\natural}$ extends to a conformal embedding $V \otimes W \subset V^{f\natural}$ (we will choose embeddings from the lists compiled in \cite{MR867023,MR867243}), and suppose furthermore that $W$ is equipped with some $N{=}1$ superconformal vector $\tau_W \in W$. Then any superconformal vector $\tau_V \in V$ will produce a superconformal vector $\tau = \tau_V + \tau_W \in V^{f\natural}$, which is unique with known stabilizer by Duncan's theorem. Automorphisms of $V$ and $W$, unless they are broken by the extension, will produce automorphisms of $V^{f\natural}$. Assuming we understand $\Aut_{N{=}1}(W) \subset \Aut_{N{=}1}(V^{f\natural}) = \Co_1$,  we will have good control over $\Aut_{N{=}1}(V)$ and hence over the $N{=}1$ structure on $V$.

This is exactly what we will do in the ``Suzuki chain'' cases $V_\ev = \SU(12)_1$, $\Sp(2{\times}6)_1$, $\Sp(2{\times}3)_1^2$, $\Sp(2{\times}3)_2$, and $\SU(6)_1^2$, where we will take $W$ to be $\SO(2)_3$, $\SO(3)_3$, $\SO(4)_3$, $\SO(5)_3$, and $\SO(2)_3^2$, respectively, each time equipped with the $N{=}1$ structure constructed in \S\ref{exist.spinm}. The exceptional isomorphisms of even subalgebras $\text{``}\Spin(2)_3\text{''} = \rU(1)_{12}$, $\text{``}\Spin(3)_3\text{''} = \Sp(2{\times}1)_6$, $\Spin(4)_3 = \Sp(2{\times}1)_3^2$, $\Spin(5)_3 = \Sp(2{\times}2)_3$, and $(\SO(2)_3^2)_\ev = \rU(1)_{6}^2$ illustrate in each cases that $W$ is the ``level-rank dual'' of $V$. The automorphism group of $W$ is in each case a symmetric group (\S\ref{sec.spinm3auto}), broken to an alternating group by the extension to $V^{f\natural}$. The automorphism group of $V$ will then be the normalizer-mod-normal of an alternating group in $\Co_1$.

In this way we will recover the larger groups in the \define{Suzuki chain} 
$$ \Suz{:}2 \supset \rG_2(4){:}2 \supset \rJ_2{:}2 \supset \rU_3(3){:}2 \supset \rL_2(7){:}2 \supset A_4{:}2 \supset A_3{:}2 $$
of subgroups (except that $\Suz{:}2$ projectively embeds) of $\Co_1$. 
The Suzuki chain was first constructed by Thompson and named for its largest member; 
details of the construction first appeared in Wilson's article \cite[\S2.2]{MR723071} classifying the maximal subgroups of $\Co_1$.  
Thompson's construction goes as follows. There is a unique conjugacy class of inclusions $A_9 \subset \Co_1$. The nontrivial double cover $2\Co_1$ pulls back to the nontrivial double cover $2A_9$, and acts on the $24$-dimensional representation through the inclusion $2A_9 \subset \SO^+(8) \cong \SO(8) \subset \SO(24)$, 
where $2A_9 \subset \SO^+(8)$  lifts the standard inclusion $A_9 \subset \SO(8)$, and $\SO(8) \subset \SO(24)$ is diagonal.  Now let $m \in \{2,\dots,8\}$ and consider the standard inclusion $A_{m+1} \subset A_9$. Because of the appearance of triality, we find that $A_{m+1} \subset \SO(m)$ maps into $\SO^+(24)$ through an exceptional isomorphism $\Spin(m) \cong (\dots) \subset \SO(24)$.
Now 
write $NA_{m+1}$ for its normalizer in $\Co_1$. The $m$th entry in the Suzuki chain arises as the normalizer-mod-normal $NA_{m+1}/A_{m+1}$ of this subgroup, and its derived subgroup as the centralizer-mod-center.

Wilson identifies four more ``Suzuki chain'' groups: for $m = 2,3,4$, the centralizer $CA_{m+1} \subset A_{m+1}$ contains a second copy of $A_{m+1}$, and when $m=2$ there is yet a third copy of $A_3 \subset C(A_3^2)$, and so one can form the normalizers-mod-normals
$$ \frac{NA_3^2}{A_3^2} = \rU_4(3){:}D_8, \quad \frac{NA_4^2}{A_4^2} = 2 \times 2^4.A_5.2, \quad \frac{NA_5^2}{A_5^2} = 2 \times D_{10}.2, \quad \frac{NA_3^3}{A_3^3} = 3^4{:}2S_4^2.$$
The first of these appears in our Theorem. The normalizer-mod-normal construction of the others, and also of the smaller Suzuki-chain groups $\rL_2(7){:}2$, $A_4{:}2$ and $A_3{:}2$, correspond to other interesting sub-VOAs of $V^{f\natural}$, but they do not arise in our Theorem. For instance $NA_4^2/A_4^2$ acts naturally on the coset of a $\PSp(2{\times}1)_6^2 \subset V^{f\natural}$; this coset is the $\bZ_4$ simple current extension of $\Spin(6)_4 = \SU(4)_4$, with even subalgebra the non-simply connected WZW algebra $\SO(6)_4$.

For the other three SVOAs in our Theorem, we study the conformal embeddings $\SU(6)_2 \times \rU(1)_{12} \subset \rU(6)_2 \subset \Sp(2{\times}6)_1$, $\Spin(12)_2 \subset \SU(12)_1$, and $\Spin(12)_1^2 \subset \Spin(24)_1$. Each of these is selected by an order-2 automorphism of the larger VOA, explaining why the corresponding centralizers appear in the Theorem.

We now provide details.

\subsection{Uniqueness for $\SU(12)_1$} \label{su12.uniqueness}
  
  Let $V$ denote the simple current extension of $V_\ev = \SU(12)_1$, and set $W = \SO(2)_3 = \rU(1)_3$, with even subalgebra $W_\ev = \rU(1)_{12}$. There is a conformal embedding $V_\ev \otimes W_\ev = \SU(12)_1 \otimes \rU(1)_{12} \subset \Spin(24)_1 \subset \SO(24)_1 = V^{f\natural}$, familiar from level-rank duality \cite{MR867023,MR867243}. We claim that this extends to a conformal embedding $V \otimes W \subset V^\natural$. There are two things to check. First, the even subalgebra $(V\otimes W)_\ev$ is not just $V_\ev \otimes W_\ev = \SU(12)_1 \otimes \rU(1)_{12}$, but rather $(V_\ev \otimes W_\ev) \oplus (V_\odd \otimes W_\odd)$, and so we need to see that $\SU(12)_1 \otimes \rU(1)_{12}$ extends inside $\Spin(24)_1 = V^{f\natural}_\ev$ to $ V_\odd \otimes W_\odd$. For both $\SU(12)_1$ and $\rU(1)_{12}$, all anyons are abelian and the groups for each are naturally isomorphic to $\bZ_{12}$. For $i \in \bZ_{12}$, let us write $\SU(12)_1(i)$ or $\rU(1)_{12}(i)$ for the corresponding anyon. In this notation, 
   $V_\odd = \SU(12)_1(6)$ and $W_\odd= \rU(1)_{12}(6)$.
   Then $\Spin(24)_1$ decomposes over $\SU(12)_1 \otimes \rU(1)_{12}$ as
   $$ V^{f\natural}_\ev = \Spin(24)_1 = \bigoplus_{i \text{ even}} \SU(12)_1(i) \otimes \rU(1)_{12}(i), $$
   which definitely contains $\SU(12)_1(6) \otimes \rU(1)_{12}(6)$. (Of course, $i$ is considered modulo $12$.) Second, we need to confirm that $V^{f\natural}_\odd$ contains $(V \otimes W)_\odd = (V_\odd \otimes W_\ev) \oplus (V_\ev \otimes W_\odd)$. One might expect that $V_\odd$ is the above sum but with odd $i$s, but this is false --- that sum does define a $\Spin(24)_1$-anyon, but it has conformal dimension $\frac12$ (and the extension of $\Spin(24)_1$ through it is the free fermion algebra $\Fer(24) = \SO(24)_1$ and not $V^{f\natural} = \SO^+(24)_1$). The other two $\Spin(24)_1$-anyons decompose as
   $$ \bigoplus_{i \text{ even}} \SU(12)_1(i) \otimes \rU(1)_{12}(i + 6), \qquad \bigoplus_{i \text{ odd}} \SU(12)_1(i) \otimes \rU(1)_{12}(i + 6).$$
   In fact, these anyons are exchanged by an outer automorphism of $\Spin(24)_1$, which does not preserve the chosen $\SU(12)_1 \otimes \rU(1)_{12}$ subalgebra. The copy of $V^{f\natural}$ that we want is the one for which $V^{f\natural}_\odd$ is the first of these two. This manifestly contains $(V_\odd \otimes W_\ev) \oplus (V_\ev \otimes W_\odd) = (\SU(12)_1(6) \otimes \rU(1)_{12}(0)) \oplus (\SU(12)_1(0) \otimes \rU(1)_{12}(6))$.
   
   Suppose now that $\tau_V$ is any superconformal vector on $V$, and let $\tau_W$ denote the canonical superconformal vector for $W$ constructed in \S\ref{exist.spinm}. Consider the $N{=}1$ algebra $V \otimes W$ with superconformal vector $\tau = \tau_V + \tau_W$. It carries a manifest action by $\Aut_{N{=}1}(V, \tau_V) \times \Aut_{N{=}1}(W) = \operatorname{Stab}_{2\PSU(12){:}2}(\tau_V) \times S_3$. Not every automorphism in $\Aut_{N{=}0}(V \otimes W) = 2\PSU(12){:}2 \times \SO(2){:}2$ lifts to an automorphism of $V^{f\natural}$, and the lift may not be quite canonical. 
   Rather, the group $2\PSU(12)$ centrally extends to $6\PSU(12) = \SU(12)/\bZ_2$, the cartesian product becomes a central product, and the diagonal ``$2$'' extends but the individual ones do not, so that
   $$ \SO^+(24) \supset \frac{6\PSU(12) \times \SO(2)}{\bZ_3} : 2.$$
   In particular, all of $\Aut_{N{=}1}(W) = S_3$ lifts to automorphisms of $V^{f\natural}$ preserving the supersymmetry. The lifts of $\bZ_3 \subset S_3$ commute with all of $V$ (and in particular they fix $\tau_V$), whereas the reflections in $S_3$ act nontrivially on $V$ (by the ``${:}2$'' inside $\Aut_{N{=}0}(V) = 2\PSU(12){:}2$).
      
   The group $\Aut_{N{=}1}(V^{f\natural}) = \Co_1$ contains four conjugacy classes of elements of order $3$. We know which one lifts $\bZ_3 \subset \Aut_{N{=}1}(W)$, because the $\Spin(2) \subset \Spin(24)$ in question acts on the 24-dimensional vector representation without fixed points, and this holds for only one of the four conjugacy classes.
   The centralizer of this $\bZ_3$ inside $\Co_1$ is isomorphic to $3\Suz$, and $3\Suz \subset 6\PSU(12)$ acts trivially on $W$ and hence fixes $\tau_W$. Of course, $3\Suz \subset \Co_1$ also fixes $\tau = \tau_V + \tau_W$, and so $3\Suz$ fixes $\tau$. The central ``$3$'' acts nontrivially on $V^{f\natural}$ but trivially on its subalgebra~$V$, and we find that $\Aut_{N{=}1}(V, \tau_V) \supset \Suz$.
   There is a unique conjugacy class of embeddings $\Suz \subset \Aut_{N{=}0}(V) = 2\PSU(12){:}2$, and we saw in \S\ref{su12.existence} that $\Suz$ preserves a unique superconformal vector~$\tau_V$, and the automorphism group of~$\tau_V$ contains $\Suz{:}2$. Furthermore, the full group $\Aut_{N{=}1}(V)$, although it doesn't have to centralize $\bZ_3 \subset \Co_1$, must certainly normalize it. The normalizer of this $\bZ_3$ is precisely $\Suz{:}2$.

\subsection{Uniqueness for $\Sp(2{\times}6)_1$}\label{unique.sp61}

We set $V = \PSp(2{\times}6)_1$, with $V_\ev = \Sp(2{\times}6)_1$, and choose $W = \text{``}\SO(3)_3\text{''} = \PSp(2{\times}1)_6$, with $W_\ev = \Sp(2{\times}1)_6$. We claim that there is a conformal embedding $V \otimes W \subset V^{f\natural}$. The conformal embedding $V_\ev \otimes W_\ev = \Sp(2{\times}6)_1 \otimes \Sp(2{\times}1)_6 \subset \Spin(24)_1 =  V^{f\natural}_\ev$ is well known from level-rank duality.
As in \S\ref{su12.uniqueness}, we may show the stronger SVOA claim by decomposing $V^{f\natural}$ over $V_\ev \otimes W_\ev$. Again one must be a little careful: the usual explanation of level-rank duality, for this rank and level, involves embedding dual pairs of WZW algebras into $\Fer(24)$, which is a simple current extension of $\Fer(24)_\ev = \Spin(24)_1$ but not the one we want.  Regardless, the usual story already lets one quickly see that $(V\otimes W)_\ev \subset \Spin(24)_1$, and then one must be careful with the odd part. The spin-$\frac32$ fields of $V^{f\natural}_\odd$ form the $2048$-dimensional positive spinor representation of $\Spin(24)$, and so the claim amounts to observing that the restriction $\mathbf{2048}|_{\Sp(2{\times}6) \times \Sp(2{\times}1)}$ contains $\mathbf{429} \otimes \mathbf{1}$ and $\mathbf{1} \otimes \mathbf{7}$ as direct summands.

Almost all of $\Aut_{N{=}0}(V) \times \Aut_{N{=}0}(W) = (\PSp(2{\times}6) \times 2) \times (\PSp(2{\times}1) \times 2)$ lifts to $\Aut_{N{=}0}(V^{f\natural}) = \SO^+(24)$: the individual ``${\times}2$''s do not, but the diagonal does. (The usual embedding for level-rank duality is $\Sp(2{\times}6) \circ \Sp(2{\times}1) \subset \SO(24)$, but we are using instead $\SO^+(24)$.)

Now let $\tau_V$ be an arbitrary superconformal vector on $V$ and take $\tau_W$ to be the $S_4$-fixed superconformal vector on $W$ constructed in \S\ref{exist.spinm}. Then $\tau = \tau_V + \tau_W$ is a superconformal vector on~$V^{f\natural}$. 
It is stabilized by the subgroup of $\Aut_{N{=}1}(V,\tau_V) \times \Aut_{N{=}1}(W,\tau_W) = \Aut_{N{=}1}(V,\tau_V) \times S_4$ which extends to $V^{f\natural}$. 
The full stabilizer of $\tau$ is, per Duncan's theorem, a copy of $\Co_1$, and the $A_4 \subset S_4$ subgroup is the unique conjugacy class mapping into $\SO^+(24)$  through the identified embedding of $\Sp(2{\times}1)$.

The centralizer of this $A_4$ inside $\Co_1$ is $\rG_2(4)$ living, as a subgroup of $\SO^+(24)$,  inside the $\Sp(2{\times}6)$ centralizing the $\Sp(2{\times}1) \supset A_4$. It follows that $\rG_2(4)$ fixes all of $W$, and so in particular fixes $\tau_W$. But it also fixed $\tau$, and so it fixes $\tau_V$. On the other hand, there is a unique conjugacy class of $\rG_2(4)$s inside of $\Sp(2{\times}6)$, and 
we saw in \S\ref{sp6.exist} that any choice of $\rG_2(4)$ fixes a unique superconformal vector, with automorphism group extending to $\rG_2(4){:}2$. But $\Aut_{N{=}1}(V,\tau_V)$ must normalize $A_4 \subset \Co_1$, and so $\Aut_{N{=}1}(V) = \rG_2(4){:}2$.

\subsection{Uniqueness for $\Sp(2{\times}3)_1^2$} \label{unique.sp312}

Let $V$ denote the $\bZ_2$ simple current extension of $V_\ev = \Sp(2{\times}3)_1^2$, and choose $W = \SO(4)_3$ with bosonic subalgebra $W_\ev = \Sp(2{\times}1)_3^2$. These bosonic subalgebras are level-rank dual inside $\Spin(24)_1$, and by decomposing $V^{f\natural}$ over $V_\ev \otimes W_\ev$, one finds that $V \otimes W \subset V^{f\natural}$.
Of $\Aut_{N{=}0}(V) \times \Aut_{N{=}0}(W) = \Sp(2{\times}3)^{\circ 2}{:}2 \times \Sp(2{\times}1)^{\circ 2}{:}2$, the part that lifts to $\Aut_{N{=}0}(V^{f\natural}) = \SO^+(24)$ is $(\Sp(2{\times}3)^{\circ 2} \times \Sp(2{\times}1)^{\circ 2}){:}2$.

 Equip $W$ with its $S_5$-invariant superconformal vector $\tau_W$, and let $\tau_V$ be an arbitrary superconformal vector for $V$. Then $\tau = \tau_V + \tau_W$ is a superconformal vector on $V^{f\natural}$, and so stabilized by a $\Co_1$. We may identify which $A_5 \subset \Co_1$ is the lift of $A_5 \subset \Aut_{N{=}1}(W)$ by factoring it through $\Sp(2{\times}1)^{\circ 2} = \SO(4) \subset \SO^+(24)$. We find that the centralizer in $\Co_1$ of this $A_5$ is a copy of $\rJ_2 \subset \Sp(2{\times}3)^{\circ 2}$, and since it commutes with all of $W$, it certainly preserves $\tau_W$ and hence $\tau_V = \tau - \tau_W$. Thus $\tau_V$ is the superconformal vector constructed in \S\ref{exist.sp3.2}. The full automorphism group $\Aut_{N{=}1}(W)$ must normalize $A_5 \subset \Co_1$, and so is exactly the identified $\rJ_2{:}2$.

\subsection{Uniqueness for $\Sp(2{\times}3)_2$}\label{unique.sp3}

Let $V = \PSp(2{\times}3)_2$ and $W = \SO(5)_3 = \PSp(2{\times}2)_3$, with even subalgebras $V_\ev = \Sp(2{\times}3)_2$ and $W_\ev = \Sp(2{\times}2)_3$. Again we may construct a conformal embedding $V \otimes W \subset V^{f\natural}$, extending the level-rank duality $\Sp(2{\times}3)_2 \otimes \Sp(2{\times}2)_3 \subset \Spin(24)_1$. Let $\tau_V$ denote an arbitrary superconformal structure on $V$, and $\tau_W$ the superconformal structure on $W$ constructed in \S\ref{exist.spinm} with automorphism group $S_{6}$. Look at the subgroup $A_6 \subset S_6$. Again by studying how much of $\Aut_{N{=}0}(V) \times \Aut_{N{=}0}(W)$ lifts to $V^{f\natural}$, we find that $\Aut(\tau_V)$ is precisely the normalizer of $A_6 \subset \Co_1$, modulo $A_6$ of course. The conjugacy class of the embedding $A_6 \subset \Co_1$ is uniquely determined by requesting that the corresponding embedding into $\PSO(24)$ factors through $\PSp(2{\times}2)$. The normalizer-mod-normal of this $A_6$ is precisely $\rU_3(3){:}2$, which, by~\S\ref{exist.sp3}, uniquely determines the superconformal vector.

\subsection{Uniqueness for $\SU(6)_1^2$} \label{unique.su612}

Set $V = \SU(6)_1^2 / \bZ_2$, with even subalgebra $V_\ev = \SU(6)_1^2$, and set $W = \rU(1)_3^2 = \text{``}\SO(2)_3^2\text{''}$. This $W$ is a lattice SVOA for the lattice $\sqrt{3}\bZ^2$. The even sublattice of the $\bZ^2$-lattice is a  copy of the $\sqrt{2}\bZ^2$ lattice, and so $W_\ev \cong \rU(1)_6^2$.

We will embed $V \otimes W \subset V^{f\natural}$. At the level of Lie algebras this corresponds to the embedding $\mathfrak{su}(6)^2 \times \mathfrak{u}(1)^2 \subset \mathfrak{u}(6)^2 \subset \mathfrak{u}(12) \subset \mathfrak{so}(24)$. Following this, we observe that $\SU(6)_1 \otimes \rU(1)_6 \subset (\rU(6)_1)_\ev = \Spin(12)_1$, where we have used the isomorphism $\rU(n)_1 = \Fer(2n)$ to identify $(\rU(n)_1)_\ev = \Spin(2n)_1$. Let us write $\SU(6)_1(i)$ and $\rU(1)_6(i)$, with $i \in \bZ_6$, for the anyons for these algebras, so that for example $V_\odd = \SU(6)_1(3) \otimes \SU(6)_1(3)$ and $W_\odd = \rU(1)_6(3) \otimes \rU(1)_6(3)$. Then we find that $\Spin(12)_1 = (\rU(6)_1)_\ev$ decomposes over $\SU(6)_1 \otimes \rU(1)_6$ as
$$ \Spin(12)_1 = \bigoplus_{i\text{ even}}\SU(6)_1(i) \otimes \rU(6)_1(i).$$
The nontrivial $\Spin(12)_1$-anyons $\Spin(12)_1(v)$ and $\Spin(12)_1(s^\pm)$ (of conformal dimensions $\frac12,\frac34,\frac34$, respectively) decompose as
\begin{align*}
 \Spin(12)_1(v) & = \bigoplus_{i\text{ odd}} \SU(6)_1(i) \otimes \rU(6)_1(i), \\
 \Spin(12)_1(s^+) & = \bigoplus_{i\text{ even}} \SU(6)_1(i) \otimes \rU(6)_1(i+3), \\ 
 \Spin(12)_1(s^-) & = \bigoplus_{i\text{ odd}} \SU(6)_1(i) \otimes \rU(6)_1(i+3).
\end{align*}
Finally, $V^{f\natural}$ decomposes over $\Spin(12)_1^2$ as
$$ V^{f\natural} = \bigoplus_{k \in \{0,v,s^+,s^-\}} \Spin(12)_1(k) \otimes \Spin(12)_1(k),$$
where in the sum we have written $\Spin(12)_1(0)$ for the vacuum anyon. (For comparison, $\Fer(24)$ does not contain the $s^\pm$ summands and instead contains $\Spin(12)_1(0) \oplus \Spin(12)_1(v)$ and $\Spin(12)_1(v) \oplus \Spin(12)_1(0)$.)
All together, we find that $V^{f\natural}$ does contain a copy of
$$ V \otimes W = \bigoplus_{i,j \in \{0,3\}} \SU(6)_1(i) \otimes \rU(6)_1(j),$$
as desired.

Now equip $W$ with its superconformal vector $\tau_W$, which is itself the sum of two copies of the superconformal vector for $\rU(1)_3$, and let $\tau_V$ be an arbitrary superconformal vector for $V$. Then $\tau = \tau_V + \tau_W$ is a superconformal vector for $V \otimes W$ and hence for $V^{f\natural}$, and so is stabilized inside $\Aut_{N{=}0}(V^{f\natural})$ by a copy of $\Co_1$. Not all of
 $\Aut_{N{=}0}(V) = (2\PSU(6))^{\circ2}{:}2$ and
$\Aut_{N{=}0}(W) = \rO(2)^{\circ 2}{:}2$ extend to $V^{f\natural}$, and some parts extend but only in ways that act on the other tensorand. But the connected subgroups do extend and centralize each other. In particular, $A_3^2 = \Aut_{N{=}1}(W) \cap \SO(2)^{\circ 2}$ lifts to $V^{f\natural}$, where its centralizer is a copy of $3^2{\cdot}\rU_4(3)$. Conversely, $3^2{\cdot}\rU_4(3)$ stabilizes $\tau$ and also $\tau_W$ (since it stabilizes all of $W$), and hence $\tau_V$. The central $3^2 \subset 3^2{\cdot}\rU_4(3)$ acts trivially on $V$, and so we find $\Aut_{N{=}1}(V) \supset \rU_4(3)$.

This implies that $\tau_V$ is the superconformal vector constructed in \S\ref{exist.su6.2}, and so in particular $\Aut_{N{=}1}(V) \supset \rU_4(3){:}D_8$. To identify its full automorphism group, we observe that $\Aut_{N{=}1}(V)$ must normalize $A_3^2 \subset \Aut_{N{=}1}(W)$, and so $\Aut_{N{=}1}(V) \subset \rU_4(3){:}D_8$.

\subsection{Uniqueness for $\Spin(16)_1 \times \Spin(8)_1$} \label{unique.spin16spin8}

Let $V = (\Spin(16)_1 \times \Spin(8)_1)/\bZ_2$. The automorphism group of $V_\ev$ is $\Aut(\Spin(16)_1) \times \Aut(\Spin(8)_1) = \PSO(16){:}2 \times \PSO(8){:}S_3$. The space of fields of spin $\frac32$ in $V_\odd$ is $\mathbf{128}_+ \boxtimes \mathbf{8}_+$, where $\mathbf{128}_+$ and $\mathbf{8}_+$ denote the positive half-spin representations of $\Spin(16)$ and $\Spin(8)$, respectively. Thus the choice of $V_\odd$ and breaks and extends $\PSO(8){:}S_3$ to $\SO^+(8){:}2 \cong \rO(8)$, so that
  $$ \Aut_{N{=}0}(V) = \SO^+(16) \circ \rO(8) = (\SO^+(16) \times \rO(8))/\bZ_2.$$
  
  The WZW embedding $\Spin(16)_1 \times \Spin(8)_1 \subset \Spin(24)$ is conformal; the subalgebra is selected not as a coset, but as the fixed points for a $\bZ_2$ action. The embedding extends to an SVOA embedding $V \subset V^{f\natural}$ (\S\ref{subsec.inclusions}), again selected as the fixed subalgebra for the action of some $\bZ_2 \subset \Aut_{N{=}0}(V^{f\natural}) = \SO^+(24)$. The extension $V \subset V^{f\natural}$ breaks and extends $\Aut_{N{=}0}(V)$ to $\Spin(16) \circ \Spin(8) = (\Spin(16) \times \Spin(8)) / \bZ_2^2 \subset \SO^+(24)$.
  
Let $\tau$ denote any superconformal vector in $\Spin(12)_1^2/\bZ_2$. Then $\tau$ determines a superconformal vector in $V^{f\natural}$ fixed by an order-$2$ element $g$, and conversely $g$ determines $\bigl(V,\tau\bigr)$ as an $N{=}1$ SVOA.
  
  This element $g$ lifts with order $2$ to $\Spin(24)$ and acts on the $24$-dimensional representation with spectrum $1^{16} (-1)^8$. $\Co_1$ has a unique conjugacy class of such elements, called ``$2\mathrm{a}$'' in GAP. This verifies the uniqueness of $\tau$ up to the $\Aut_{N{=}0}(V)$-action.
  
  The centralizer of $g$ in $\Co_1$ has shape $2^{1+8}_+\cdot \rO_8^+(2)$ in the ATLAS notation \cite{ATLAS}. Quotienting by the central ``$2$'' to achieve a faithful action on $V$, we find that
  $$\Aut_{N{=}1}(V) \cap (\SO^+(16) \circ \SO(8)) = 2^8 \cdot \rO_8^+(2).$$
  Thus either $\Aut_{N{=}1}(V) = 2^8 \cdot \rO_8^+(2)$ or $\Aut_{N{=}1}(V)$ contains $2^8 \cdot \rO_8^+(2)$ with index $2$. We claim that the latter occurs, and that $\Aut_{N{=}1}(V) = 2^8 \cdot \rO_8^+(2){:}2$.
  
  To prove this, we will describe the embedding $2^8 \cdot \rO_8^+(2){:}2 \subset \SO^+(16) \circ \rO(8)$ by starting with the Lie group $\rE_8$. This description is essentially equivalent to the construction from \cite{JohnDuncan-thesis}, reviewed in \S6 of \cite{MR2352133}, which produces $V^{f\natural}$ as the canonical orbifold of the ``supersymmetric $\rE_8$ torus.''
  
  Let $T \subset \rE_8$ denote the maximal torus, and $W$ the Weyl group, which we think of as a subgroup of $\rO(8)$. The normalizer of $T$ is a nonsplit extension $T\cdot W$. 
  The group $W$ has shape $2 \cdot \rO_8^+(2){:}2$ in the ATLAS notation \cite{ATLAS}.
  Let $z \in W$ denote the central element. Since it acts on $T$ as $t \mapsto t^{-1}$, its lifts to $T \cdot W$ are all conjugate, and all have order $2$. Fix such a lift $\tilde z \in T \cdot W$. Then the centralizer of $\tilde z$ in $T \cdot W$ is a group of shape $2^8 \cdot W$, where $2^8 \subset T$ is the $2$-torsion subgroup.
  The centralizer of $\tilde z$ in $\rE_8$ is isomorphic to $\SO^+(16)$. 
  Indeed, a version of this construction works for any simple Lie algebra, and selects the split form of the corresponding Lie group.
  All together we find a map $2^8 \cdot W \subset \SO^+(16)$, and of course $W \subset \rO(8)$ automatically. The combined map $2^8 \cdot W \subset \SO^+(16) \times \rO(8)$ covers a map $2^8 \cdot \rO_8^+(2){:}2 \subset \SO^+(16) \circ \rO(8)$.
  
  The adjoint representation $\mathfrak{e}_8$ splits over $\SO^+(16) \subset \rE_8$ as $\mathfrak{so}(16) \oplus \mathbf{128}_+$, where the splitting is by central characters: $\tilde z$ acts trivially on $\mathfrak{so}(16)$ and with eigenvalue $-1$ on $\mathbf{128}_+$. Let $\Delta^+$ denote the set of positive roots of $\rE_8$, with respect to some decomposition into positive and negative roots. Then $\mathfrak{e}_8$ decomposes over the maximal torus $T \subset \rE_8$ as $\mathfrak{e}_8 = \mathfrak{t} \oplus \bigoplus_{r \in \Delta^+} \mathbf{1}_r \oplus \mathbf{1}_{-r}$, where $\mathbf{1}_{\pm r}$ is the weight space of weight $\pm r$, and $\mathfrak{t} = \mathbf{8}_0$ is the Cartan subalgebra. The $\tilde{z}$-action exchanges $\mathbf{1}_r$ with $\mathbf{1}_{-r}$, and acts with eigenvalue $-1$ on $\mathfrak{t}$. So $\mathfrak{so}(16)$ has a basis consisting of the vectors $(1 + \tilde z)e_r$, where $r \in \Delta^+$ and $e_r$ is a basis vector of $\mathbf{1}_r$. For comparison, $\mathbf{128}_+$ decomposes as $\mathfrak{t} \oplus \bigoplus_r (1 - \tilde z)e_r$. (Note in particular that the Cartan subalgebra $\mathfrak{t}$ does not meet $\mathfrak{so}(16)$. The inclusion $\SO^+(16) \subset \rE_8$ takes a maximal torus to a maximal torus, but not to the maximal torus $T$ that we are using. Indeed, $2^8 \subset \SO^+(16)$ is not a subgroup of any maximal torus thereof: its double cover inside $\Spin(16)$ is an extraspecial group $2^{1+8}_+$.)
  
  These decompositions diagonalize the $2^8$ action, which acts on $\mathbf{1}_r$ with character ``$r \mod 2$.'' Since roots are primitive vectors, $2^8$ acts nontrivially on $(1 \pm \tilde z)e_r$, and the only $2^8$-fixed subspace of $\mathfrak{128}_+$ is $\mathfrak{t}$, which is isomorphic to $\mathbf{8}$ as a $(2^8 \cdot W)$-module. Note that $\mathbf{8}$ is simple upon restriction to the even subgroup $2\rO_8^+(2) = W \cap \SO(8)$. Thus $2^8 \cdot \rO_8^+(2)$ has a unique fixed point inside $\mathbf{128}_+ \otimes \mathbf{8}$, which by the earlier remarks must be the superconformal vector $\tau$.
  
  But the same argument shows that $\tau$ is in fact fixed by all of $2^8 \cdot \rO_8^+(2){:}2$, and as we remarked above, this is the largest that $\Aut_{N{=}1}(V)$ could be. We remark that the inclusion $2^8 \cdot \rO_8^+(2){:}2 \subset \SO^+(16) \circ \rO(8)$ can be ``twisted'' by the nontrivial map 
$2^8 \cdot \rO_8^+(2){:}2 \to 2 \subset \SO^+(16) \circ \rO(8)$, where the second map selects the centre, thereby producing a second nonconjugate copy of $2^8 \cdot \rO_8^+(2){:}2$ inside $\SO^+(16) \circ \rO(8)$. This latter copy ``antifixes'' $\tau$.

\subsection{Uniqueness for $\Spin(12)_1^2$} \label{unique.spin1212}

Let $V = \Spin(12)_1^2/\bZ_2$, with even subalgebra $\Spin(12)_1^2$.
The WZW embedding $\Spin(12)_1^2 \subset \Spin(24)$ is conformal; the subalgebra is selected not as a coset, but as the fixed points for a $\bZ_2$ action. The embedding extends to an SVOA embedding $V \subset V^{f\natural}$ (\S\ref{subsec.inclusions}), again selected as the fixed subalgebra for the action of some $\bZ_2 \subset \Aut_{N{=}0}(V^{f\natural}) = \SO^+(24)$, and all of $\Aut_{N{=}0}(V) = \SO^+(12)^{\circ 2}{:}2$ extends (via a double cover) to $V^{f\natural}$.

Let $\tau$ denote any superconformal vector in $\Spin(12)_1^2/\bZ_2$. Then $\tau$ determines a superconformal vector in $V^{f\natural}$ fixed by an order-$2$ element $g$, and conversely $g$ determines $\bigl(V,\tau\bigr)$ as an $N{=}1$ SVOA.

This element $g$ lifts with order $2$ to $\Spin(24)$ and acts on the $24$-dimensional representation with spectrum $1^{12} (-1)^{12}$. $\Co_1$ has a unique conjugacy class of such elements, called ``$2\mathrm{c}$'' in GAP. This verifies the uniqueness of $\tau$ up to the $\Aut_{N{=}0}(V)$-action. The centralizer of $g$ in $\Co_1$ is $2^{11}{:}\rM_{12}{:}2$. Quotienting by the central ``$2$'' to achieve a faithful action on $V$, we find that $\Aut_{N{=}1}(V) = 2^{10}{:}\rM_{12}{:}2$, and that $\tau$ is the superconformal vector constructed in \S\ref{exist.spin12.2}.

\subsection{Uniqueness for $\Spin(12)_2$}\label{unique.spin12}

The (nonconformal) embedding $\Spin(12)_2 \subset \Spin(24)_1$ is unusual among WZW embeddings. Its level-rank dual is $\text{``}\Spin(2)_{12}\text{''} \subset \rU(1)_{12}$, but that has as its coset all of $\SU(12)_1$, and not just its subalgebra $\Spin(12)_2$. The reason for this unusual behaviour is that the embedding $\Spin(12)_2 \subset \SU(12)_1$ is already conformal, and so has trivial coset.  $\Spin(12)_2$ is the fixed points of an outer $\bZ_2$-action on $\SU(12)_1$.

We showed in \S\ref{subsec.inclusions} that the VOA inclusion $\Spin(12)_2 \subset \SU(12)_1$ extends to an SVOA inclusion $\SO^+(12)_2 \subset 6\PSU(12)_1$. Let us repeat the argument: it suffices to study the spin-$\frac32$ fields, and $\Alt^6(\mathbf{12}) = \mathbf{924}$ splits over $\Spin(12)$ as $\mathbf{462}_+ \oplus \mathbf{462}_-$ (with notation as in \S\ref{exist.spin12} and \S\ref{su12.existence}).

Choose an arbitrary superconformal vector $\tau$ on $\SO^+(12)_2$. Its image in $6\PSU(12)_1$ is also a superconformal vector, and so is unique up to $\SU(12)$-automorphisms, with stabilizer $\Suz{:}2$ therein, by \S\ref{su12.uniqueness}. Since $\tau \in \SO^+(12)_2$, it is fixed by the chosen outer automorphism. There are two conjugacy classes of order-$2$ elements in $\Suz{:}2$ covering the outer $2$. They can be distinguished as follows: one of them lifts to an order-$2$ element inside $2\Suz{:}2$ and the other lifts with order $4$. The one that cuts out $\Spin(12) \subset \SU(12)$ lifts with order $2$. The centralizer of any such element is a group of shape $2\rM_{12}{:}2$. It acts on $\SO^+(12)_2$ through $\rM_{12}{:}2$, and preserves $\tau$, and so we find that $\Aut_{N{=}1}(\SO^+(12)_2,\tau) \cong \rM_{12}{:}2$. But $\rM_{12}{:}2$ embeds into $\Aut_{N{=}0}(V)$ in a unique-up-to-conjugation way, preserving a (unique) superconformal vector by \S\ref{exist.spin12}.

\subsection{Uniqueness for $\SU(6)_2$}\label{unique.su6}

Let $V = \SU(6)_2/\bZ_2$, with $V_\ev = \SU(6)_2$, and set $W = \SO(2)_3 = \rU(1)_3$, with $W_\ev = \rU(1)_{12}$. Arguing as above, it is not hard to construct a conformal embedding $V \otimes W \subset \PSp(2{\times}6)_1$. Write $\tau_W$ for the unique superconformal vector in $W$; then $\Aut_{N{=}1}(W) = S_3$. Choose any superconformal vector $\tau_V$ for $V$. Then $\tau_V + \tau_W = \tau$ is a superconformal vector for 
$\PSp(2{\times}6)_1$,
 and so its automorphism group therein is $\rG_2(4){:}2$ by \S\ref{unique.sp61}. Consider the $\bZ_3 \subset S_3 = \Aut_{N{=}1}(W)$. It is compatible with the extension, and lifts to an order-$3$ element in $\rG_2(4)$. There are two (conjugacy classes of) such elements: the one coming from our $\rU(1)$ is the unique one acting in the $12$-dimensional representation without fixed points. Its normalizer in $\rG_2(4){:}2$ has shape $3\rM_{21}{:}2^2$, and so $\Aut_{N{=}1}(V,\tau_V) = \rM_{21}{:}2^2$. (The normal ``$3$'' acts nontrivially on $\Sp(2{\times}6)_1$, but trivially on the subalgebra $V \otimes W$.) On the other hand, we saw in \S\ref{exist.su6} that $\rM_{21}$ has a unique conjugacy class of embeddings into $\Aut_{N{=}0}(V)$ and 
 preserves a unique superconformal vector.

\section{Type $\rE$ cases} \label{sec.typeE}

We now discuss the Type $\rE$ cases. In \S\ref{orbifolds} and \S\ref{exist.e71.2} we construct three different $N{=}1$ structures for $\rE_{7,1}^2$; the first two are constructed by very similar methods, and the last by following the approach from Section~\ref{sec.existence}. We conjecture that these three $N{=}1$ structures are in the same orbit under the action of $\Aut_{N{=}0}(\rE_{7,1}^2/\bZ_2) = \rE_7^{\circ 2}{:}2$, but this is not obvious from the construction, and we do not attempt a proof.
We then explain in \S\ref{subsec.e72} and \S\ref{subsec.failureE} why $N{=}1$ structures for $\rE_{7,2}$ and $\rE_{8,2}$ appear unlikely.

\subsection{Existence for $\rE_{7,1}^2$ from orbifolding} \label{orbifolds}

The goal of this section is to construct an $N{=}1$ structure on the SVOA $V = \rE_{7,1}^2/\bZ_2$, the $\bZ_2$ simple current extension of $V_\ev = \rE_{7,1}^2$; we will in fact construct two $N{=}1$ structures, and it is not obvious whether they are isomorphic. Our first step will be to give a ``free fermionic'' construction of $V$ just as an SVOA. This is reasonable because $V$ is \define{holomorphic}: like $\rE_{8,1}$ and $\SO^+(24)_1 = V^{f\natural}$, $V$ has no nontrivial anyons. Indeed, $\rE_{7,1}$ has only two anyons --- the vacuum and an abelian anyon of conformal dimension $\frac32$ corresponding to the centre of $\rE_7$ --- and so the anyons of $\rE_{7,1}^2$ form a Klein-4 group. The $\bZ_2$-extension cuts down the number of anyons by a factor of $2^2$. (When there are nonabelian anyons, they must be counted weighted by the squares of their quantum dimensions. For the general statement, see \cite{MR1936496}.)

To construct $V$ as an SVOA, we will start with a different holomorphic SVOA of the same central charge, namely the free fermion algebra $\Fer(28) = \SO(28)_1$, and orbifold (aka gauge) a symmetry. Orbifolding has been part of the VOA story from the beginning: 
it was pioneered by \cite{MR996026}, who used it in their construction of the moonshine module. The first step is to choose an action of a finite group on the to-be-orbifolded SVOA. We will choose the  elementary abelian group $\bZ_2^3$ acting on $\Fer(28)$ as follows: write $28 = 4 \times 7$, and act by four copies of the seven-dimensional ``regular minus trivial'' representation of $\bZ_2^3$. (This seven-dimensional representation is essentially Hamming's error correcting code $\mathrm{Ham}(7,4)$. One choice for the generators of $\bZ_2^3$ is to have them act on $\bR^7$ via the diagonal matrices $\diag(-1,-1,-1,-1,1,1,1)$, $\diag(-1,-1,1,1,-1,-1,1)$, and $\diag(-1,1,-1,1,-1,1,-1)$.)

The second step is to pass to the $\bZ_2^3$-fixed sub-SVOA $\Fer(28)^{\bZ_2^3} \subset \Fer(28)$. (This fixed subalgebra is sometimes called the ``chiral orbifold'' in the VOA literature, whereas what we will construct is called variously the ``full orbifold'' or ``twisted orbifold.'') It is expected that for any rational SVOA and any finite group action, the fixed sub-SVOA is again rational; such a result is not known in general, but is proven for solvable groups in \cite{CarnahanMiyamoto}. In particular, $\Fer(28)^{\bZ_2^3}$ is rational. Once rationality is known, the main result of \cite{MR1923177} identifies the category of anyons of $\Fer(28)^{\bZ_2^3}$ as the Drinfel'd centre of a superfusion category $\cat{SVec}^\alpha[\bZ_2^3]$ arising as a twisted form of the finite group $\bZ_2^3$. (The papers \cite{MR1923177,MR1936496,CarnahanMiyamoto} only address the bosonic case, but they imply the corresponding fermionic statements with only a small amount of extra work.)

In particular, the anyons of $\Fer(28)^{\bZ_2^3}$ are graded by the group $\bZ_2^3$. (If $\bZ_2^3$ were nonabelian, then the anyons would be graded merely by the conjugacy classes in $\bZ_2^3$.) Suppose that we can choose one abelian anyon of each grading, such that the choices are closed under fusion.  The \define{orbifold} $\Fer(28) \sslash \bZ_2^3$, if it exists, is the simple current extension of $\Fer(28)^{\bZ_2^3}$ by these abelian anyons. In the context of orbifolds, the (generators of) the extending anyons are called \define{twist fields}. The double slash in the notation reminds that orbifolding is a two-step process: restrict to a subalgebra, and then choose an extension.

The existence of the orbifold is controlled by the \define{'t Hooft anomaly} of the action of $\bZ_2^3$ on $\Fer(28)$ (and the choices are controlled by choices of how to trivialize this anomaly). The 't Hooft anomaly appears as the ``twisting'' $\alpha$ in the superfusion category $\cat{SVec}^\alpha[\bZ_2^3]$, analogous to a twisted group algebra.
By definition, $\cat{SVec}^\alpha[\bZ_2^3]$ has simple objects $M_g$ indexed by $g \in \bZ_2^3$, and these should fuse compatibly with the group multiplication.
 In the bosonic case, the twisting $\alpha$ consists merely of associator data (aka F-matrices, aka 6j-symbols) for the fusion category $\cat{Vec}^\alpha[\bZ_2^3]$, and defines a class in ordinary group cohomology $\H^3(\bZ_2^3;\rU(1))$. In the fermionic case the anomaly is more complicated (cf.\ \cite[Section~4]{GJFIII}), and has three {layers}. In a supercategory, there are two types of simple objects: \define{ordinary} objects, with endomorphism algebra $\bC$, and \define{Majorana} objects, with endomorphism algebra $\Cliff(1)$.
The bottom \define{Majorana layer} of $\alpha$ records which objects in $\cat{SVec}^\alpha[\bZ_2^3]$ are ordinary and which are Majorana; it is a class $\alpha^{(1)} \in \H^1(\bZ_2^3;\bZ_2)$. Next is the \define{Gu--Wen layer} $\alpha^{(2)}$. When all objects are ordinary (i.e.\ when the Majorana layer vanishes), the Gu--Wen layer records whether the isomorphism $M_g \otimes M_h \cong M_{gh}$ is even or odd, and defines a class in $\H^2(\bZ_2^3;\bZ_2)$. (When $\alpha^{(1)}$ does not vanish, $\alpha^{(2)}$ is not a cocycle, but rather solves $\d \alpha^{(2)} = \Sq^2\alpha^{(1)}$.)
Finally there is the \define{Dijkgraaf--Witten layer} $\alpha^{(3)}$ prescribing the associator, which is a 3-cochain on $\bZ_2^3$ with values in $\rU(1)$. When $\alpha^{(1)}$ and $\alpha^{(2)}$ are both trivialized, $\alpha^{(3)}$ is a cocycle for ordinary cohomology. In general, it solves $\d\alpha^{(3)} = (-1)^{\Sq^2\alpha^{(2)}}(\dots)$, where the $(\dots)$ term depends on $\alpha^{(1)}$; an exact formula does not seem to appear in the literature, but could in principle be extracted from \cite[\S5.4]{MR3978827}.
All together, the three layers $\alpha^{(1)},\alpha^{(2)},\alpha^{(3)}$ compile into a class $\alpha \in \SH^3(\bZ_2^3)$ in a generalized cohomology theory called \define{supercohomology}. A supercohomology class with vanishing Majorana layer is said to live in \define{restricted supercohomology}~$\rSH^3(\bZ_2^3)$.

To evaluate $\alpha$, we first note that the action of $\bZ_2^3$ on $\Fer(28)$ consists of four copies of the action on $\Fer(7)$, and that anomalies add under stacking (aka tensoring) of SVOAs. We claim that the Majorana layer for the action of $\bZ_2^3$ on $\Fer(7)$ vanishes, i.e.\ that the anomaly lives in restricted supercohomology. (In fact, the Gu--Wen layer of that action also vanishes, but we will not use this.) Indeed, for actions on free fermions, the Majorana layer  merely records the determinant of the corresponding map $\bZ_2^3 \to \rO(7)$. Each element in $\bZ_2^3$ switches the signs of four fermions, and so has positive determinant. The restricted supercohomology group $\rSH^3(\bZ_2^3)$ is an extension of the ordinary cohomology groups $\H^2(\bZ_2^3;\bZ_2)$ and $\H^3(\bZ_2^3;\rU(1))$, each of which is $2$-torsion, and so $\rSH^3(\bZ_2^3)$ is $4$-torsion. Thus the anomaly $\alpha$ of the action of $\bZ_2^3$ on $\Fer(28)$ must vanish, being $4$ times the anomaly of the action on $\Fer(7)$.

Therefore the orbifold SVOA $\Fer(28) \sslash \bZ_2^3$ exists. We will now identify it. In general, there can be multiple non-isomorphic results of orbifolding, because to construct the orbifold requires not just that the anomaly vanishes, but also the choice of a trivialization of the anomaly. In the case of $\Fer(28) \sslash \bZ_2^3$, we will see that such choices do not effect the isomorphism type of the resulting orbifold.

To recognize $\Fer(28) \sslash \bZ_2^3$, first note that the $\bZ_2^3$ symmetry \define{screens} all the free fermions, in the sense that the fixed subalgebra $\Fer(28)^{\bZ_2^3}$ has no free fermions. Second, each twist field has conformal dimension $1$. Indeed, it is a general fact that, for an order-2 element $g \in \Aut(\Fer(n)) = \rO(n)$, the corresponding twist field has conformal dimension $\frac1{16} \times\dim(\text{$(-1)$-eigenspace of $g$})$. (This is why, for instance, the orbifold $\Fer(16) \sslash (-1)^f = \rE_{8,1}$ includes so many extra spin-$1$ fields.) It follows that $\Fer(28) \sslash \bZ_2^3$ has no free fermions. Finally, not all of the spin-$1$ fields in $\Fer(28)$ are screened: the spin-$1$ fields in $\Fer(28)$ are bilinears in the free fermions, and the ones which are not screened are those for which both free fermions come from the same four-dimensional eigenspace of $\bZ_2^3$; thus the spin-$1$ fields in the fixed subalgebra $\Fer(28)^{\bZ_2^3}$, a subalgebra of $\Fer(28)\sslash \bZ_2^3$, form the Lie algebra $\mathfrak{so}(4)^7 = \mathfrak{su}(2)^{14}$. This Lie algebra has rank $14$, equal to the central charge of $\Fer(28) \sslash \bZ_2^3$. It is a general fact that in a (unitary) SVOA of central charge $c$, the spin-$1$ fields form a Lie algebra with rank $c$ if and only if the SVOA is a lattice SVOA \cite{MR2097833}, and that the lattice is self-dual if and only if the SVOA is holomorphic, and that the lattice has vectors of length $1$ if and only if the SVOA has free fermions.

All together, we find that $\Fer(28) \sslash \bZ_2^3$ is isomorphic to a lattice SVOA for a self-dual lattice of rank $14$ without any vectors of length $1$. The classification of self-dual lattices of rank $\leq 23$ is available in \cite[Chapter 16]{MR1662447} (and extended to rank $\leq 25$ in \cite{MR1801771}). Inspecting the list, one finds that there is a unique self-dual lattice of rank $14$ without vectors of length $1$: its root system is $\rE_7^2$, extended by a glue vector of length $\sqrt{3}$. This establishes the existence of an isomorphism $\Fer(28) \sslash \bZ_2^3 \cong \rE_{7,1}^2 / \bZ_2$. (One may also establish such an isomorphism directly without appealing to the classification of lattices: it suffices to understand $\Fer(28) \sslash \bZ_2^3$ as a simple current extension of its subalgebra $\SU(2)_1^{14}$; the extension adjoins new roots indexed by the Hamming code $\mathrm{Ham}(7,4)$, and one recovers in this way a description of $\mathfrak{e}_7$ as a module over its subalgebra $\mathfrak{su}(2)^7$.)

Having established an isomorphism $\Fer(28) \sslash \bZ_2^3 \cong \rE_{7,1}^2 / \bZ_2$, to construct an $N{=}1$ structure it suffices to equip $\Fer(28)$ with an $N{=}1$ structure which is invariant under the $\bZ_2^3$-action (as then the superconformal vector will survive the screening to $\Fer(28)^{\bZ_2^3}$). $N{=}1$ structures on $\Fer(28)$ correspond bijectively to $28$-dimensional semisimple Lie algebras \cite{MR791865}, or more precisely to Lie algebra structures on $\bR^{28}$ whose Killing form agrees with (minus) the Euclidean inner product. There are two $28$-dimensional Lie algebras with a fixed-point-free action by $\bZ_2^3$: $\mathfrak{so}(8)$ and~$\mathfrak{g}_2^2$. Each of these equips $\Fer(28) \sslash \bZ_2^3 \cong \rE_{7,1}^2 / \bZ_2$ with an $N{=}1$ structure. The $N{=}1$ structures on $\Fer(28)$ corresponding to~$\mathfrak{so}(8)$ and~$\mathfrak{g}_2^2$ are not isomorphic, but the orbifold involves adding new fields of spin~$1$ and hence new $N{=}0$ automorphisms (compare the construction from \cite{MR996026}: the moonshine module $V^\natural$ has many automorphisms which are not visible from the orbifold construction), and it is not obvious whether the resulting $N{=}1$ structures on $\rE_{7,1}^2 / \bZ_2$ are isomorphic or not.

\subsection{Existence for $\rE_{7,1}^2$ from finite groups}\label{exist.e71.2}

We now construct an $N{=}1$ structure on $V = \rE_{7,1}^2/\bZ_2$ by the employing the strategy from Section~\ref{sec.existence}. It is not obvious whether the result is isomorphic to either of the $N{=}1$ structures constructed in \S\ref{orbifolds}.

The automorphism group of $V_\ev = \rE_{7,1}^2$ is $\rP\rE_7\wr 2 = \rP\rE_7^2{:}2$, where $\rP\rE_7$ denotes the adjoint form of the simply connected group $\rE_7$. Write $\mathbf{56}$ for the smallest irrep of $\rE_7$; it is symplectic. The space of spin-$\frac32$ fields in $V_\odd$ is $\mathbf{56}^{\boxtimes 2}$, and so $\Aut_{N{=}0}(V) = \rE_7^{\circ 2}{:}2$.

The (isomorphism, not conjugacy, classes of) quasisimple subgroups of $\rE_7$ are listed in \cite{MR1653177}. One of them, which following the ATLAS we will call $\rU_3(8)$ (it is recorded as $\PSU_3(8)$ in \cite{MR1653177}), leaves both $\mathbf{56}$ and the adjoint representation $\mathfrak{e}_7 = \mathbf{133}$ simple. 
In fact, $\rU_3(8)$ has only one 56-dimensional irrep, but three $133$-dimensional irreps; these are permuted by the outer automorphism group of $\rU_3(8)$, which is isomorphic to $3 \times S_3$, and so we find in fact three conjugacy classes of embeddings $\rU_3(8) \subset \rE_7$. Choose an outer automorphism of order $2$. It fixes one of the three conjugacy classes of $\rU_3(8) \subset \rE_7$, and exchanges the other two; we will use the two that are exchanged to provide a map $\rU_3(8) \subset \rE_7^{\circ 2}$. 

By uniqueness, the $56$-dimensional irrep doesn't care about this choice. It follows that $\mathbf{56}^{\boxtimes 2}$ has a unique $\rU_3(8)$-fixed point $\tau$. We can extend the inclusion $\rU_3(8) \subset \rE_7^{\circ 2}$ to an inclusion $\rU_3(8){:}2 \subset \rE_7^{\circ 2}{:}2$ in two different ways, differing by a sign. By uniqueness, one of these extensions fixes $\tau$, and the other antifixes it.

Since we chose the outer automorphism to exchange the $133$-dimensional adjoint representations of the two copies of $\rE_7$, we find that the full adjoint $\mathfrak{e}_7^{\boxplus 2}$ of $\rE_7^{\circ 2}{:}2$ remains simple upon restriction to $\rU_3(8){:}2$. Thus the conformal vector $\nu$ is the unique spin-$2$ field fixed by $S = \rU_3(8){:}2$.

One of the order-3 outer automorphisms of $\rU_3(8)$, called ``$3_1$'' in the ATLAS, commutes with the order-$2$ outer automorphism. Choosing this outer automorphism,
the embeddings $\rU_3(8) \subset \rE_7$ and $\rU_3(8){:}2 \subset \rE_7^{\circ 2}{:}2$ extend (uniquely) to the groups $\rU_3(8){:}3$ and $\rU_3(8){:}6$, respectively. Indeed, $\rU_3(8){:}3$ has three $56$-dimensional irreps: one is quaternionic, and the other two are complex and dual to each other, and the embedding into $\rE_7$ chooses the quaternionic representation. It follows that $\tau$ is in fact preserved by $\rU_3(8){:}6$.

\subsection{Doubt for $\rE_{7,2}$} \label{subsec.e72}

Let $W$ denote the level-$1$ $N{=}1$ minimal model, meaning the minimal-central-charge SVOA generated by a superconformal vector $\tau_W$. Its bosonic subalgebra $W_\ev$ is the level-$2$ Virasoro minimal model of central charge $\frac{7}{10}$, also called the \define{tricritical Ising} algebra. Let $V = \rP\rE_{7,2}$ denote the $\bZ_2$-extension of $\rE_{7,2}$, of central charge $13\frac{3}{10}$. There is a conformal embedding $V \otimes W \subset \rE_{7,1}^2/\bZ_2$. In fact, $W$ is the coset of $V$ inside $\rE_{7,1}^2/\bZ_2$. In particular, the image of $\tau_W$ in $\rE_{7,1}^2$ is the unique spin-$\frac32$ field invariant under the diagonal $\rE_7$-action.

If $V$ admits a superconformal vector $\tau_V$, then $\tau = \tau_V + \tau_W$ will be a superconformal vector on $\rE_{7,1}^2$. Conversely, to give $V$ a supersymmetry is equivalent to finding a supersymmetry on $\rE_{7,1}^2$ such that the embedding $W \subset \rE_{7,1}^2$ is supersymmetric. For this, it would suffice to embed $W \subset \Fer(28)$ supersymmetrically and compatibly with the $\bZ_2^3$-orbifold used in \S\ref{orbifolds}.

There are many close relationships between tricritical Ising, the Lie algebra $\mathfrak{g}_2$, and fermions \cite{MR1354601}, stemming from the fact that $W$ also arises as the coset of $G_{2,1}$ inside $\Fer(7) = \SO(7)_1$. Unfortunately, we were unable to find a suitable supersymmetric embedding $W \subset \Fer(28)$.

Another reason to doubt the existence of an $N{=}1$ structure for $V = \rP\rE_{7,2}$ comes from directly analyzing the equations. The space of spin-$\frac32$ fields is the $\rE_7$-irrep $\mathbf{1463}$, appearing inside $\Sym^2(\mathbf{56}) = \mathfrak{sp}(2{\times}28)$ as
$$ \Sym^2(\mathbf{56}) = \mathbf{1463} \oplus \mathbf{133}.$$
The space of spin-$2$ fields is
$$ V_2 = \mathbf{1}\oplus \mathbf{1539} \oplus \mathbf{7371} \oplus \partial(V_1),$$
where $\partial(V_1)$ means the fields of the form $\frac\partial{\partial z} X(z)$ for $X \in V_1 = \mathbf{133}$.
The OPE $\mathbf{1463} \otimes \mathbf{1463} \to V_2$ includes the map $\mathbf{1463} \otimes \mathbf{1463} \to \mathbf{1}\oplus \mathbf{1539} = \Alt^2(\mathbf{56})$ that sends
$$ A_{ij} \otimes B_{kl} \mapsto A_{ij} \omega^{jk} B_{kl} + B_{ij}\omega^{jk} A_{jl},$$
where $A,B \in \mathbf{1463} \subset \Sym^2(\mathbf{56})$ are thought of as symmetric (complex-valued) matrices, and $\omega^{jk}$ is the symplectic form on the $56$-dimensional representation. Thus if $\tau \in \mathbf{1463} \subset \mathfrak{sp}(2{\times}28)$ is to be a superconformal vector, then it must solve $(\tau\omega)^2 = \mathrm{id}$ when thought of as a $56 \times 56$ matrix. This forces $\tau$ to have spectrum $1^{28} (-1)^{28}$, and its centralizer under the $\Sp(2{\times}28)$-action is necessarily a copy of $\mathrm{GL}(28) \subset \Sp(2{\times}28)$. Conversely, any embedding $\mathrm{GL}(28) \subset \Sp(2{\times}28)$ will provide a candidate superconformal vector (namely the image of the centre of $\mathfrak{gl}(28)$) inside $\Sym^2(\mathbf{56})$. But a necessary condition for it to be a superconformal vector is that it should live inside $\mathbf{1463} \subset \Sym^2(\mathbf{56})$, which is the orthocomplement of the adjoint representation $\mathbf{133} \subset \Sym^2(\mathbf{56})$.

Thus we are faced with an interesting question in Lie theory: choose an element of $\mathfrak{sp}(2{\times}28)$ whose centralizer is a $\mathfrak{gl}(28)$ (these are, up to rescaling, a single conjugacy class); does there exist an $\mathfrak{e}_7$ Lie algebra inside its orthocomplement? More generally, one can ask: choose an element of $\mathfrak{sp}(2{\times}n)$ whose centralizer is a $\mathfrak{gl}(n)$; when does its orthocomplement contain a Lie subalgebra $\mathfrak{g}$ such that the defining representation $\mathbf{2n}$ remains simple upon restriction to $\mathfrak{g}$? It is possible to find various Lie subalgebras $\mathfrak{g}$ of the orthocomplement such that $\mathbf{2n}|_{\mathfrak{g}}$ splits as $\mathbf{n} \oplus \mathbf{n}$ or $\mathbf{n} \oplus \bar{\mathbf{n}}$ (for instance, there is a $\mathfrak{gl}(n)$ in the orthocomplement; cf.~the MathOverflow conversation \href{https://mathoverflow.net/q/337233/}{\texttt{mathoverflow.net/q/337233/}}), but subalgebras for which $\mathbf{2n}$ remains simple are hard to come by, and no examples have yet turned up.

\subsection{Doubt for $\rE_{8,2}$} \label{subsec.failureE}

Let $V = \rE_{8,2} / \bZ_2$. The spin-$\frac32$ fields form the second smallest fundamental representation $\mathbf{3875}$ of $\rE_8$. We record \cite{MR1117679}:
\begin{align*}
  V_1 & = \mathfrak{e}_8 = \mathbf{248} \\
  V_{3/2} & = \mathbf{3875} \\
  V_2 & = \mathbf{1} \oplus \mathbf{3875} \oplus \mathbf{27000} \oplus \partial(V_1) \\
  \Sym^2(V_{3/2}) & = \mathbf{1} \oplus \mathbf{3875} \oplus \mathbf{27000} \oplus \mathbf{147250} \oplus \mathbf{2450240} \oplus \mathbf{4881384}
\end{align*}
Thus to give a superconformal vector for $V$ is to give a real vector $\tau \in \mathbf{3875}$ such that $\tau^2 \in \Sym^2(\mathbf{3875})$ has no components when projected to the $\mathbf{3875}$ and $\mathbf{27000}$ summands.

The $\rE_8$-invariant map $\Sym^2(\mathbf{3875}) \to \mathbf{3875}$ can be thought of as a commutative but nonassociative algebra structure on $\mathbf{3875}$. It is studied in \cite[Section 7]{MR3406824}, where it is shown that $\rE_8$ is exactly the stabilizer of $\star$, and so $(\mathbf{3875},\star)$ can be thought of as a ``Griess algebra'' for the Lie group $\rE_8$.

Any superconformal vector $\tau$ will satisfy in particular $\tau \star \tau = 0$. This equation is easy to solve among complex vectors --- for instance, take $\tau$ to be a highest weight vector --- but not among vectors of nonzero norm and especially not among real vectors. The Conjecture asserting that 
$\rE_{8,2} / \bZ_2$ does not admit an $N{=}1$ structure
is based on the belief that indeed there are no real solutions $\tau \in \mathbf{3875}$ to $\tau \star \tau = 0$.

There is a maximal abelian subgroup $2^5 \subset \rE_8$, with normalizer the Dempwolff group $2^5\cdot \mathrm{SL}_5(\bF_2)$ \cite{MR0399193}. Its action on $\rE_{8,2}$ is nonanomalous, being twice the (bosonic) anomaly of $\rE_{8,1}$, and the orbifold $V \sslash 2^5$ is isomorphic to $\Fer(31)$. One would be able to give $V$ an $N{=}1$ structure following the techniques of \S\ref{orbifolds} if there were a 31-dimensional semisimple Lie algebra with an appropriate $2^5$ grading, but such a Lie algebra does not exist.

\section*{Acknowledgements}

I thank the anonymous referee for their many helpful comments and suggestions.
I also thank M.\ Bischoff, R.\ Derryberry, J.\ Duncan, D.\ Gaiotto, A.\ Henriques, G.\ Moore, and D.\ Treumann for their conversations and suggestions. 
Parts of this work were performed while I was a visitor at the Mathematical Institute at Oxford (supported by the EPSRC grant EP/M024830/1), the International Centre for Mathematical Sciences at Edinburgh, and the Aspen Center for Physics (supported by the NSF grant PHY-1607611).
 Research at the Perimeter Institute is supported by the Government of Canada through Industry Canada and by the Province of Ontario through the Ministry of Economic Development and Innovation.

%\bibliography{ReferencesWithLinks}{}

\begin{thebibliography}{GTVW14}

\bibitem[BB87]{MR867243}
F.~Alexander Bais and Peter~G. Bouwknegt.
\newblock A classification of subgroup truncations of the bosonic string.
\newblock {\em Nuclear Phys. B}, 279(3-4):561--570, 1987.
\newblock \DOI{10.1016/0550-3213(87)90010-1}. \MRnumber{867243}.

\bibitem[Bor00]{MR1801771}
Richard~E. Borcherds.
\newblock Classification of positive definite lattices.
\newblock {\em Duke Math. J.}, 105(3):525--567, 2000.
\newblock \DOI{10.1215/S0012-7094-00-10536-4}. \MRnumber{1801771}.
  \arXiv{math/9912236}.

\bibitem[Bre]{CTblLib}
Thomas Breuer.
\newblock {\em CTblLib - The GAP Character Table Library}.

\bibitem[Car20]{MR4050091}
Scott Carnahan.
\newblock Building {V}ertex {A}lgebras from {P}arts.
\newblock {\em Comm. Math. Phys.}, 373(1):1--43, 2020.
\newblock \MRnumber{4050091}. \DOI{10.1007/s00220-019-03607-0}.
  \arXiv{1408.5215}.

\bibitem[CCN{\etalchar{+}}85]{ATLAS}
J.~H. Conway, R.~T. Curtis, S.~P. Norton, R.~A. Parker, and R.~A. Wilson.
\newblock {\em Atlas of finite groups}.
\newblock Oxford University Press, Eynsham, 1985.
\newblock Maximal subgroups and ordinary characters for simple groups, With
  computational assistance from J. G. Thackray.

\bibitem[CM16]{CarnahanMiyamoto}
Scott Carnahan and Masahiko Miyamoto.
\newblock Regularity of fixed-point vertex operator subalgebras.
\newblock 2016.
\newblock \arXiv{1603.05645}.

\bibitem[CS99]{MR1662447}
J.~H. Conway and N.~J.~A. Sloane.
\newblock {\em Sphere packings, lattices and groups}, volume 290 of {\em
  Grundlehren der Mathematischen Wissenschaften [Fundamental Principles of
  Mathematical Sciences]}.
\newblock Springer-Verlag, New York, third edition, 1999.
\newblock With additional contributions by E. Bannai, R. E. Borcherds, J.
  Leech, S. P. Norton, A. M. Odlyzko, R. A. Parker, L. Queen and B. B. Venkov.

\bibitem[dG]{SLA-GAP}
Willem~Adriaan de~Graaf.
\newblock {\em SLA - Computing with simple Lie algebras - a GAP package}.

\bibitem[DM04]{MR2097833}
Chongying Dong and Geoffrey Mason.
\newblock Rational vertex operator algebras and the effective central charge.
\newblock {\em Int. Math. Res. Not.}, (56):2989--3008, 2004.
\newblock \DOI{10.1155/S1073792804140968}. \MRnumber{2097833}.
  \arXiv{math/0201318}.

\bibitem[DMNO13]{MR3039775}
Alexei Davydov, Michael M{\"u}ger, Dmitri Nikshych, and Victor Ostrik.
\newblock The {W}itt group of non-degenerate braided fusion categories.
\newblock {\em J. Reine Angew. Math.}, 677:135--177, 2013.
\newblock \MRnumber{3039775}. \arXiv{1009.2117}.

\bibitem[Dun06]{JohnDuncan-thesis}
John F.~R. Duncan.
\newblock {\em Vertex operators, and three sporadic groups}.
\newblock PhD thesis, Yale, 2006.

\bibitem[Dun07]{MR2352133}
John~F. Duncan.
\newblock Super-moonshine for {C}onway's largest sporadic group.
\newblock {\em Duke Math. J.}, 139(2):255--315, 2007.
\newblock \DOI{10.1215/S0012-7094-07-13922-X}. \MRnumber{2352133}.
  \arXiv{math/0502267}.

\bibitem[DZ05]{MR2175996}
Chongying Dong and Zhongping Zhao.
\newblock Modularity in orbifold theory for vertex operator superalgebras.
\newblock {\em Comm. Math. Phys.}, 260(1):227--256, 2005.
\newblock \DOI{10.1007/s00220-005-1418-2}. \MRnumber{2175996}.

\bibitem[DZ10]{MR2681777}
Chongying Dong and Zhongping Zhao.
\newblock Modularity of trace functions in orbifold theory for {$\Bbb
  Z$}-graded vertex operator superalgebras.
\newblock In {\em Moonshine: the first quarter century and beyond}, volume 372
  of {\em London Math. Soc. Lecture Note Ser.}, pages 128--143. Cambridge Univ.
  Press, Cambridge, 2010.
\newblock \MRnumber{2681777}.

\bibitem[FBZ04]{MR2082709}
Edward Frenkel and David Ben-Zvi.
\newblock {\em Vertex algebras and algebraic curves}, volume~88 of {\em
  Mathematical Surveys and Monographs}.
\newblock American Mathematical Society, Providence, RI, second edition, 2004.
\newblock \DOI{10.1090/surv/088}. \MRnumber{2082709}.

\bibitem[FLM88]{MR996026}
Igor Frenkel, James Lepowsky, and Arne Meurman.
\newblock {\em Vertex operator algebras and the {M}onster}, volume 134 of {\em
  Pure and Applied Mathematics}.
\newblock Academic Press, Inc., Boston, MA, 1988.
\newblock \MRnumber{996026}.

\bibitem[FQS84]{MR740343}
Daniel Friedan, Zongan Qiu, and Stephen Shenker.
\newblock Conformal invariance, unitarity, and critical exponents in two
  dimensions.
\newblock {\em Phys. Rev. Lett.}, 52(18):1575--1578, 1984.
\newblock \DOI{10.1103/PhysRevLett.52.1575}. \MRnumber{740343}.

\bibitem[Fuc91]{MR1096120}
J\"{u}rgen Fuchs.
\newblock Simple {WZW} currents.
\newblock {\em Comm. Math. Phys.}, 136(2):345--356, 1991.
\newblock \MRnumber{1096120}.

\bibitem[GAP]{GAP}
The GAP~Group.
\newblock {\em {GAP -- Groups, Algorithms, and Programming, Version 4.10.2}}.
\newblock \url{https://www.gap-system.org}.

\bibitem[GG15]{MR3406824}
Skip Garibaldi and Robert~M. Guralnick.
\newblock Simple groups stabilizing polynomials.
\newblock {\em Forum Math. Pi}, 3:e3, 41, 2015.
\newblock \arXiv{1309.6611}. \DOI{10.1017/fmp.2015.3}. \MRnumber{3406824}.

\bibitem[GJF18]{GJFIII}
Davide Gaiotto and Theo Johnson-Freyd.
\newblock Holomorphic {SCFTs} with small index.
\newblock 2018.
\newblock \arXiv{1811.00589}.

\bibitem[GJF19]{MR3978827}
Davide Gaiotto and Theo Johnson-Freyd.
\newblock Symmetry protected topological phases and generalized cohomology.
\newblock {\em J. High Energy Phys.}, (5):007, 34, 2019.
\newblock \DOI{10.1007/JHEP05(2019)007}. \arXiv{1712.07950}.
  \MRnumber{3978827}.

\bibitem[GO85]{MR791865}
P.~Goddard and D.~Olive.
\newblock Kac-{M}oody algebras, conformal symmetry and critical exponents.
\newblock {\em Nuclear Phys. B}, 257(2):226--252, 1985.
\newblock \DOI{10.1016/0550-3213(85)90344-X}. \MRnumber{791865}.

\bibitem[GR99]{MR1653177}
Robert~L. Griess, Jr. and A.~J.~E. Ryba.
\newblock Finite simple groups which projectively embed in an exceptional {L}ie
  group are classified!
\newblock {\em Bull. Amer. Math. Soc. (N.S.)}, 36(1):75--93, 1999.
\newblock \DOI{10.1090/S0273-0979-99-00771-5}. \MRnumber{1653177}.

\bibitem[GS88]{MR968813}
Peter Goddard and Adam Schwimmer.
\newblock Factoring out free fermions and superconformal algebras.
\newblock {\em Phys. Lett. B}, 214(2):209--214, 1988.
\newblock \DOI{10.1016/0370-2693(88)91470-0}. \MRnumber{968813}.

\bibitem[GTVW14]{GTVW}
Matthias~R. Gaberdiel, Anne Taormina, Roberto Volpato, and Katrin Wendland.
\newblock A k3 sigma model with $\mathbb{Z}^8_2 : \mathbb{M}_{20}$ symmetry.
\newblock {\em J. High Energ. Phys.}, 22, 2014.
\newblock \DOI{10.1007/JHEP02(2014)022}. \arXiv{1309.4127}.

\bibitem[H\"{o}h96]{MR1614941}
Gerald H\"{o}hn.
\newblock {\em Selbstduale {V}ertexoperatorsuperalgebren und das
  {B}abymonster}, volume 286 of {\em Bonner Mathematische Schriften [Bonn
  Mathematical Publications]}.
\newblock Universit\"{a}t Bonn, Mathematisches Institut, Bonn, 1996.
\newblock \MRnumber{1614941}.

\bibitem[Hen17]{HenriquesWZW}
Andr\'{e} Henriques.
\newblock The classification of chiral {WZW} models by {$H^4_+(BG,\Bbb{Z})$}.
\newblock In {\em Lie algebras, vertex operator algebras, and related topics},
  volume 695 of {\em Contemp. Math.}, pages 99--121. Amer. Math. Soc.,
  Providence, RI, 2017.
\newblock \DOI{10.1090/conm/695/13998}. \MRnumber{3709708}. \arXiv{1602.02968}.

\bibitem[HK07]{HeluaniKac2007}
Reimundo Heluani and Victor~G. Kac.
\newblock Susy lattice vertex algebras.
\newblock 2007.
\newblock \arXiv{0710.1587}.

\bibitem[HM]{HarveyMoore-unpub}
Jeffrey~A. Harvey and Gregory~W. Moore.
\newblock Unpublished.

\bibitem[Hua08]{MR2468370}
Yi-Zhi Huang.
\newblock Rigidity and modularity of vertex tensor categories.
\newblock {\em Commun. Contemp. Math.}, 10(suppl. 1):871--911, 2008.
\newblock \DOI{10.1142/S0219199708003083}. \MRnumber{2468370}.
  \arXiv{math.QA/0502533}.

\bibitem[JFT18]{JFT}
Theo Johnson-Freyd and David Treumann.
\newblock {$\mathrm{H}^4(\mathrm{Co}_0;\mathbf{Z}) = \mathbf{Z}/24$}.
\newblock {\em Int. Math. Res. Not. IMRN}, 2018.
\newblock \arXiv{1707.07587}.

\bibitem[Kac98]{MR1651389}
Victor Kac.
\newblock {\em Vertex algebras for beginners}, volume~10 of {\em University
  Lecture Series}.
\newblock American Mathematical Society, Providence, RI, second edition, 1998.
\newblock \DOI{10.1090/ulect/010}. \MRnumber{1651389}.

\bibitem[Kir02]{MR1923177}
Alexander Kirillov, Jr.
\newblock Modular categories and orbifold models.
\newblock {\em Comm. Math. Phys.}, 229(2):309--335, 2002.
\newblock \DOI{10.1007/s002200200650}. \MRnumber{1923177}.
  \arXiv{math/0104242}.

\bibitem[KMPS90]{MR1117679}
S.~Kass, R.~V. Moody, J.~Patera, and R.~Slansky.
\newblock {\em Affine {L}ie algebras, weight multiplicities, and branching
  rules. {V}ols. 1, 2}, volume~9 of {\em Los Alamos Series in Basic and Applied
  Sciences}.
\newblock University of California Press, Berkeley, CA, 1990.
\newblock \MRnumber{1117679}.

\bibitem[KO02]{MR1936496}
Alexander Kirillov, Jr. and Viktor Ostrik.
\newblock On a {$q$}-analogue of the {M}c{K}ay correspondence and the {ADE}
  classification of {$\mathfrak{sl}_2$} conformal field theories.
\newblock {\em Adv. Math.}, 171(2):183--227, 2002.
\newblock \DOI{10.1006/aima.2002.2072}. \MRnumber{1936496}.
  \arXiv{math/0101219}.

\bibitem[KS04]{MR2104437}
Yvette Kosmann-Schwarzbach.
\newblock Derived brackets.
\newblock {\em Lett. Math. Phys.}, 69:61--87, 2004.
\newblock \MRnumber{2104437}. \DOI{10.1007/s11005-004-0608-8}.
  \arXiv{math/0312524}.

\bibitem[Li01]{MR1822111}
Haisheng Li.
\newblock Certain extensions of vertex operator algebras of affine type.
\newblock {\em Comm. Math. Phys.}, 217(3):653--696, 2001.
\newblock \DOI{10.1007/s002200100386}. \MRnumber{1822111}.
  \arXiv{math/0003038}.

\bibitem[LL04]{MR2023933}
James Lepowsky and Haisheng Li.
\newblock {\em Introduction to vertex operator algebras and their
  representations}, volume 227 of {\em Progress in Mathematics}.
\newblock Birkh\"{a}user Boston, Inc., Boston, MA, 2004.
\newblock \DOI{10.1007/978-0-8176-8186-9}. \MRnumber{2023933}.

\bibitem[Mooa]{GregLecture-Freedfest}
Gregory~W. Moore.
\newblock K3 surfaces, mathieu moonshine, and quantum codes.
\newblock Lecture given at ``Freed60," University of Austin, January 15, 2019.
  \url{http://www.physics.rutgers.edu/~gmoore/Freed60-Final-2019.pdf}.

\bibitem[Moob]{GregLecture-Readfest}
Gregory~W. Moore.
\newblock Moonshine phenomena, supersymmetry, and quantum codes.
\newblock Lecture given at ``Field theory in condensed matter," Yale, April 12,
  2019. \url{http://www.physics.rutgers.edu/~gmoore/Read60-April10-2019.pdf}.

\bibitem[MT10]{MR2648364}
Geoffrey Mason and Michael Tuite.
\newblock Vertex operators and modular forms.
\newblock In {\em A window into zeta and modular physics}, volume~57 of {\em
  Math. Sci. Res. Inst. Publ.}, pages 183--278. Cambridge Univ. Press,
  Cambridge, 2010.
\newblock \MRnumber{2648364}.

\bibitem[Sch]{Kac-computer}
Bert Schellekens.
\newblock {\em {Kac - Komputations with Algebras and Currents}}.
\newblock \url{https://www.nikhef.nl/~t58/Site/Kac.html}.

\bibitem[SV95]{MR1354601}
S.~L. Shatashvili and C.~Vafa.
\newblock Superstrings and manifolds of exceptional holonomy.
\newblock {\em Selecta Math. (N.S.)}, 1(2):347--381, 1995.
\newblock \DOI{10.1007/BF01671569}. \MRnumber{1354601}. \arXiv{hep-th/9407025}.

\bibitem[SW86]{MR867023}
A.~N. Schellekens and N.~P. Warner.
\newblock Conformal subalgebras of {K}ac-{M}oody algebras.
\newblock {\em Phys. Rev. D (3)}, 34(10):3092--3096, 1986.
\newblock \DOI{10.1103/PhysRevD.34.3092}. \MRnumber{867023}.

\bibitem[Tho76]{MR0399193}
J.~G. Thompson.
\newblock A conjugacy theorem for {$E_{8}$}.
\newblock {\em J. Algebra}, 38(2):525--530, 1976.
\newblock \DOI{10.1016/0021-8693(76)90235-0}. \MRnumber{0399193}.

\bibitem[vE13]{MR3077918}
Jethro van Ekeren.
\newblock Modular invariance for twisted modules over a vertex operator
  superalgebra.
\newblock {\em Comm. Math. Phys.}, 322(2):333--371, 2013.
\newblock \DOI{10.1007/s00220-013-1758-2}. \MRnumber{3077918}.

\bibitem[vE14]{MR3205090}
Jethro van Ekeren.
\newblock Vertex operator superalgebras and odd trace functions.
\newblock In {\em Advances in {L}ie superalgebras}, volume~7 of {\em Springer
  INdAM Ser.}, pages 223--234. Springer, Cham, 2014.
\newblock \DOI{10.1007/978-3-319-02952-8_13}. \MRnumber{3205090}.

\bibitem[Wil83]{MR723071}
Robert~A. Wilson.
\newblock The maximal subgroups of {C}onway's group {${\rm Co}_{1}$}.
\newblock {\em J. Algebra}, 85(1):144--165, 1983.
\newblock \DOI{10.1016/0021-8693(83)90122-9}.

\bibitem[WPN{\etalchar{+}}]{AtlasRep-Gap}
Robert~A. Wilson, Richard~A. Parker, Simon Nickerson, John~N. Bray, and Thomas
  Breuer.
\newblock {\em AtlasRep - A GAP Interface to the Atlas of Group
  Representations}.

\bibitem[Zhu96]{MR1317233}
Yongchang Zhu.
\newblock Modular invariance of characters of vertex operator algebras.
\newblock {\em J. Amer. Math. Soc.}, 9(1):237--302, 1996.
\newblock \DOI{10.1090/S0894-0347-96-00182-8}. \MRnumber{1317233}.

\end{thebibliography}
%\bibliographystyle{alpha}

\newcommand{\etalchar}[1]{$^{#1}$}

\end{document}